\documentclass[oneside, pdflatex, sn-mathphys-num]{sn-jnl}
\usepackage{geometry}

\usepackage{graphicx}%
\usepackage{multirow}%
\usepackage{amsmath,amssymb,amsfonts}%
\usepackage{amsthm}%
\usepackage{mathrsfs}%
\usepackage[title]{appendix}%
\usepackage{xcolor}%
\usepackage{textcomp}%
\usepackage{manyfoot}%
\usepackage{booktabs}%
\usepackage{algorithm}%
\usepackage{algorithmicx}%
\usepackage{algpseudocode}%
\usepackage{listings}%
\usepackage{comment}%
\usepackage{subcaption} 
\usepackage{mathtools}
\usepackage{graphicx}
\usepackage{listings}
\usepackage{xcolor}
\usepackage{titlesec}
\usepackage{float} 
\theoremstyle{thmstyleone}%
\newtheorem{theorem}{Theorem}
\newtheorem{proposition}[theorem]{Proposition}%

\theoremstyle{thmstyletwo}%
\newtheorem{example}{Example}%
\newtheorem{remark}{Remark}%

\theoremstyle{thmstylethree}%
\newtheorem{definition}{Definition}%

\begin{document}

\flushbottom

\title[Article Title]{Optimal Blackjack Betting Strategies Through Dynamic Programming and Expected Utility Theory}

\author*[1,2]{\fnm{Lucas} \sur{Bordeu}}\email{lubordeu@gmail.com}
\author[1,3]{\fnm{Javier} \sur{Castro}}\email{javier.castro@bjtheorem.com}

\affil[1]{\orgname{Blackjack Theorem SpA}, \orgaddress{\city{Santiago}, \postcode{8320000}, \country{Chile}}}
\affil[2]{\orgdiv{Departamento de Ingeniería Industrial}, \orgname{Facultad de Ciencias Físicas y Matemáticas, Universidad de Chile}, \orgaddress{\city{Santiago}, \country{Chile}}}
\affil[3]{\orgdiv{Departamento de Ciencias de la Computación}, \orgname{Facultad de Ciencias Físicas y Matemáticas, Universidad de Chile}, \orgaddress{\city{Santiago}, \country{Chile}}}


\abstract{This study presents a rigorous mathematical approach to the optimization of round and betting policies in Blackjack, using Markov Decision Processes (MDP) and Expected Utility Theory. The analysis considers a direct confrontation between a player and the dealer, simplifying the dynamics of the game. The objective is to develop optimal strategies that maximize expected utility for risk profiles defined by constant (CRRA) and absolute (CARA) aversion utility functions. Dynamic programming algorithms are implemented to estimate optimal gambling and betting policies with different levels of complexity. The evaluation is performed through simulations, analyzing histograms of final returns.
The results indicate that the advantage of applying optimized round policies over the “basic strategy” is slight, highlighting the efficiency of the last one. In addition, betting strategies based on the exact composition of the deck slightly outperform the Hi-Lo counting system, showing its effectiveness. The optimized strategies include versions suitable for mental use in physical environments and more complex ones requiring computational processing. Although the computed strategies approximate the theoretical optimal performance, this study is limited to a specific configuration of rules. As a future challenge, it is proposed to explore strategies under other game configurations, considering additional players or deeper penetration of the deck, which could pose new technical challenges.}




\maketitle

\begin{footnotesize}
\noindent \textbf{Collaboration Note:} This work was carried out in collaboration 
with Prof. Andrés Abeliuk (Department of Computer Science, Universidad de Chile), 
who provided valuable feedback on the structure and academic content of this manuscript.
\end{footnotesize}

\section{Introduction}

The card game “Blackjack” has captured the imagination of millions of people around the world because of its apparent simplicity, but underlying complexity. Over the years, strategies to increase the chances of winning at this card game have evolved from mere tips and tricks, to complex mathematical and computational systems. Motivated by curiosity and passion for Blackjack, this paper studies optimal strategies for the game, modeling the problem as a particular case of the dynamic portfolio optimization problem, highlighting the game's connection to financial theory.

The player's complete strategy is divided into a "round policy" and a "betting policy." The round policy consists of decision-making during each round of play and has been extensively explored by the community. Its optimization is considered partially resolved \cite{peytavi2018matematicas, mgp_analyzer, k_c_cdca}, culminating in good approximations limited by computational costs. In contrast, the optimization of the betting policy remains a less explored area. Although some studies have examined the problem in a simplified manner \cite{vince2013optimal, werthamer2005optimal, jensen2014expected}, none have formally modeled it or addressed it rigorously.

In this paper, a rigorous mathematical approach is proposed for optimizing both policies, employing the theoretical framework of Markov Decision Processes (MDP) and expected utility theory. The developed model considers the relevant variables influencing the player's returns, assuming that the round policy maximizes the expected return of an individual round, while the betting policy seeks to maximize the expected utility of an entire betting session. This approach is versatile enough to capture known strategies, such as the Basic Strategy \cite{shackleford_basic_strategy} and betting strategies based on the Hi-Lo card counting system \cite{wikipedia2024conteo}.

The optimization of the round policy is performed through dynamic programming, using efficient algorithms based on vector structures to minimize the required computational costs. As a result, a semi-optimal round policy is computed, which includes a simplification in the case of the "Splitting" option. Additionally, the probability mass function of the return for each round is estimated, both under the Basic Strategy and the semi-optimal policy. On the other hand, the betting policy is also optimized through dynamic programming, considering a discretization of the state space and a simplification of the transition function. Betting policies based on the Hi-Lo counting system and full knowledge of the deck composition are considered.

The study concludes with the derivation of optimized betting policies for different levels of risk aversion, along with analytical estimates of the return distributions. This is the first work to rigorously address the problem of betting optimization in blackjack, and it is conjectured that the obtained results are close to the theoretical optimum under the assumed simplifications. Therefore, this study represents a formal solution to the game under a specific configuration of rules and conditions, and it opens the possibility of solving it for other configurations and generalizations of the game in future research.

The article is structured to ensure logical continuity and facilitate understanding. First, a comprehensive theoretical framework is presented that mathematically models the optimization problem, including optimality conditions under certain assumptions. Then, the methodology employed to compute the policies and evaluate their performance is detailed. Finally, the obtained results are discussed, and the study is concluded by highlighting its implications and suggesting possible lines of future research.

\subsection{Game Configuration}

There are multiple aspects that fully define a session of Blackjack, so it is essential to be precise when referring to the game as the object of study. This section describes the "game dynamics" and the "rules and conditions" assumed throughout the research. Game dynamics refer to a description of what the game involves and how it unfolds. The rules and conditions refer to a specific set of rules and assumptions considered for this particular study.

\subsubsection{Game Dynamics}

Blackjack is a card game played at a casino table or through online platforms, where multiple players can participate simultaneously, though each competes individually against the house, represented by the dealer. The essence of the game lies in strategy and decision-making based on the value of the dealt cards, with the goal of achieving a score close to or equal to 21, without exceeding it.

\begin{enumerate}
    \item \underline{Start of the Round and Betting:} Before the round begins, each player at the table places an initial bet. This bet is required for the player to participate in the upcoming round. Bets are placed in a designated space in front of each player.

    \item \underline{Dealing the Cards:} Once the bets are placed, the dealer proceeds to deal the cards. Each player receives two face-up cards, while the dealer receives one face-up card and, depending on the game variant, may receive another face-down card (known as the "hole card").

    \item \underline{Card Values:}
    \begin{itemize}
        \item Number cards (2 to 10): Worth their face value.
        \item Face cards (J, Q, K): Worth 10 points.
        \item Aces: Have a dual value, worth either 11 points or 1 point, depending on what is more advantageous for the player's hand. An Ace counts as 11 unless that value would cause the hand to exceed 21 points, in which case it counts as 1.
    \end{itemize}

    \item \underline{Round Progression:} After the initial deal, players evaluate their hands and proceed with their actions, choosing from the following options:
    \begin{itemize}
        \item \textbf{Hit:} Request an additional card to try to get closer to 21.
        \item \textbf{Stand:} Keep the current hand without receiving additional cards.
        \item \textbf{Double Down:} Double the initial bet in exchange for receiving only one more card.
        \item \textbf{Split:} If the two initial cards have the same value, the player can split them into two separate hands, doubling the bet. An additional card is dealt to each new hand, which are then played independently.
        \item \textbf{Surrender:} In some variants, the player can choose to surrender after receiving the first two cards, recovering half of their bet.
    \end{itemize}

    \item \underline{Dealer's Turn:} Once all players have completed their hands, the dealer plays their hand according to fixed rules: the dealer must hit until reaching a minimum of 17 points. In some variants, if the dealer has a "soft 17" (a hand that includes an Ace valued at 11 points), they must take another card.

    \item \underline{Determining the Winner:} Once the dealer has finished taking cards, the round ends, and the winner is determined:
    \begin{itemize}
        \item \textbf{Blackjack:} If a player receives an Ace and a 10-value card as the first two cards, this is called a "Blackjack," which is the best possible hand and usually pays more than a regular win, typically 3 to 2.
        \item If the player beats the dealer without exceeding 21, they win the round and receive a 1-to-1 payout on their bet.
        \item If the dealer busts (exceeds 21) and the player does not, the player wins.
        \item In case of a tie, the player's bet is returned without any additional winnings.
        \item If the dealer has a superior hand without exceeding 21, the player loses their bet.
    \end{itemize}
\end{enumerate}

Blackjack is a game that relies not only on the luck of the dealt cards but also on the player’s strategic decisions. The ability to decide when to hit, stand, double down, or split, considering the dealer’s visible card and the cards that have already been played, can significantly influence the player’s success.

\subsubsection{Rules and Conditions}

For this particular research, the following set of rules was specifically adapted to meet the technical requirements of the study:

\paragraph{}{Set of Rules \textbf{R}:}
\begin{itemize}
    \item \textbf{Initial Deck:} The initial deck from which cards are dealt consists of 8 shoe decks, where each deck contains 52 cards.
    \item \textbf{Penetration Level:} This corresponds to the percentage of cards in the deck that are played before it is reset to its full state. For this study, a cutoff point of 75\% was considered.
    \item \textbf{Dealer:} The dealer receives only one card at the beginning of the round and always stands if their score is 17 or higher.
    \item \textbf{Double Down:} Doubling down is allowed with any hand consisting of 2 cards. It is not allowed to double down after splitting. After doubling down, the player receives only one card.
    \item \textbf{Splitting:} Splitting is allowed for any hand consisting of 2 cards, as long as both cards are identical. Splitting can only be done once per round. If a pair of Aces is split, the player receives one additional card for each Ace, and if blackjack is obtained on either hand, the payout is 1 to 1.
    \item \textbf{Hit or Stand:} The player is allowed to hit or stand at any time, as long as the other rules are respected.
    \item \textbf{Blackjack:} Pays 3 to 2.
    \item \textbf{Surrender and Insurance:} Not allowed.
\end{itemize}

Based on the mentioned set of rules, the dynamics of each round are fully determined by strategies followed by the players at the table. However, the presence of multiple players at the table introduces additional complexity to the problem. Thus, it will be assumed that the player plays alone against the dealer. In this scenario, the player is in a situation of "Complete Information," as they have all the necessary information to deduce the stochastic dynamics of the game.




The paper is organized as follows. Section 2 presents the theoretical framework for modeling Blackjack optimization as a Markov Decision Process, including definitions of round and betting policies. Section 3 describes the methodology for estimating optimal policies and evaluating their performance using dynamic programming and simulations. Section 4 presents the results of the optimization, comparing strategies based on the Basic Strategy and the Hi-Lo card counting system. Finally, Section 5 concludes with a discussion of the implications of the findings and proposes directions for future research.

\section{Theoretical Framework}

\subsection{Round Policy}

In this section, the optimization problem for the round policy is mathematically formulated as a Markov decision process, assuming the absence of additional players and a fixed set of rules \textbf{R} throughout the game. A round policy makes decisions during rounds of play, with the optimal round policy being the one that maximizes the expected return of each round.

\subsubsection{Origin Deck}

If a complete record of the played cards is kept, it is possible to know the exact composition of the deck before each betting round begins. We refer to this composition as the \textit{origin deck} of the round, represented by a vector $\mathbf{d} \in \mathbb{R}^{10}$, where $d_i$ represents the number of cards in the deck with a game value of $i$. Specifically, $d_1$ represents the number of Aces, and $d_{10}$ represents the grouped number of "10," "Jacks," "Queens," and "Kings."

The range of values that $\mathbf{d}$ can take during the game is determined by the initial deck and the penetration level, and this is referred to as the \textit{origin deck space}:

\[
\Omega_{(\tau, l)} = \left\{ \mathbf{d} \in \mathbb{R}^{10} \ \middle|\ 52k(1 - \tau) \leq \sum_{i=1}^{10} d_i \leq 52k, \ d_i \leq 4k \ \forall i \leq 9, \ d_{10} \leq 16k \right\}
\]

Where:
\begin{itemize}
    \item $N_{(\tau, l)} = \text{card}(\Omega_{(\tau,l)}) = |\Omega_{(\tau, l)}| \in \mathbb{N}$: Number of possible origin decks.
    \item $l \in \mathbb{N}$: Number of decks that make up the complete initial deck. A single deck is represented by the vector $\mathbf{d} = [4,4,4,4,4,4,4,4,4,16]$.
    \item $\tau \in [0,1]$: Penetration level set for the game. This is the fraction of the complete initial deck after which a new round begins. Once this point is reached, the deck is reset to its full state $\bar{\mathbf{d}} = [4l,4l,4l,4l,4l,4l,4l,4l,4l,16l]$.
\end{itemize}

\subsubsection{Player's Sequence}

Each possible individual hand of the player can be represented by a sequence 
\(\mathbf{j} \in \mathbb{R}^{10}\), where each component of the sequence indicates the number of cards of each type that make up the hand. There are two sets of particular relevance associated with the player's hands:

\begin{enumerate}
    \item \textbf{Player's Sequences:} \\
    The set of all sequences associated with hands that have not exceeded 21 points (including those with 0 or 1 card) or that have exceeded 21 points by adding exactly one additional card is called \textit{Player's Sequences} and is denoted as \( U^J \). This set is hierarchically organized into levels or \textit{floors}, according to the number of cards that compose the player's hand:
    \[
    U^J = \bigcup_{n=0}^{22} U_n^J, \quad U_n^J = \{\mathbf{j} \in U^J \mid |\mathbf{j}| = n \}.
    \]
    \item \textbf{Relevant Player's Sequences:} \\
    An important subset of \( U^J \) is the set \( U^{(J,r)} \), called \textit{Relevant Player's Sequences}. This set includes only the sequences that do not exceed 21 points. Similar to \( U^J \), these sequences are also grouped into levels:
    \[
    U^{(J,r)} = \bigcup_{n=0}^{21} U_n^{(J,r)}, \quad U_n^{(J,r)} = \{\mathbf{j} \in U^{(J,r)} \mid |\mathbf{j}| = n \}.
    \]
\end{enumerate}

Transitions between sequences can be represented by Directed Acyclic Graphs (DAG), defined by a set of nodes \( U \) (sequences) and a set of directed edges \( A \subseteq U \times U \) (valid transitions). Each edge \( (u,v) \in A \) indicates that the sequence \( v \) is obtained by adding exactly one additional card to the sequence \( u \). In this context, three relevant graphs are identified:

\begin{itemize}
    \item \( \mathbf{G}^J = (U^J, A^J) \): DAG that connects all the player's sequences.
    \item \( \mathbf{G}^{(J,r)} = (U^{(J,r)}, A^{(J,r)}) \): DAG that connects all the relevant player's sequences.
    \item \( \mathbf{G}_T^{(J,r)} = (U^{(J,r)}, A_T^{(J,r)}) \): Directed Spanning Tree (DST) that connects the sequences in \( U^{(J,r)} \) through directed edges, where all nodes have a single parent except for the root (this tree has methodological utility).
\end{itemize}

\subsubsection{Dealer Sequence}

Each possible hand of the dealer can be represented by a sequence \(\mathbf{c} \in \mathbb{R}^{10}\), where each component refers to the number of each type of card that makes up the hand. There are two particularly relevant sets associated with the dealer's hands:

\begin{enumerate}
    \item \textbf{Dealer Sequences:} \\
    The set of all sequences associated with hands that have not exceeded 16 points (including those with 0 or 1 card) or that have surpassed 16 points by adding a single additional card is called \textit{Dealer Sequences} and is denoted by \( U^C \). This set is hierarchically organized into levels or \textit{tiers} based on the number of cards in the dealer's hand:
    \[
    U^C = \bigcup_{n=0}^{12} U_n^C, \quad U_n^C = \{\mathbf{c} \in U^C \mid |\mathbf{c}| = n\}.
    \]
    \item \textbf{Relevant Dealer Sequences:} \\
    An important subset of \( U^C \) is the set \( U^{(C,r)} \), referred to as \textit{Relevant Dealer Sequences}. This set includes only those sequences that exceed 16 points but do not surpass 21. Similar to \( U^C \), these sequences are also grouped into levels:
    \[
    U^{(C,r)} = \bigcup_{n=0}^{12} U_n^{(C,r)}, \quad U_n^{(C,r)} = \{\mathbf{c} \in U^{(C,r)} \mid |\mathbf{c}| = n\}.
    \]
\end{enumerate}

The possible transitions between dealer sequences can be represented using a Directed Acyclic Graph (DAG). Specifically, the following graph is highlighted for this study:

\begin{itemize}
    \item \( \mathbf{G}^C = (U^C, A^C) \): A DAG that connects all the dealer sequences.
\end{itemize}

\subsubsection{State Variables}

The state variables must account for all the information necessary to control the final return of the round. The state, or \textbf{instance}, can be identified through the following variables:

\paragraph{Player Hands} These are represented by sequences \(\mathbf{j}_1, \mathbf{j}_2 \in U^J\). The sequence \(\mathbf{j}_1\) refers to the player's first hand, while \(\mathbf{j}_2\) refers to the player's second hand, in the event of a \textit{Split}.

\paragraph{Dealer Hand} This is represented by a sequence \(\mathbf{c} \in U^C\). When \(|\mathbf{c}| = 1\), the dealer's hand will be denoted as \(\mathbf{c^*}\), indicating it is the dealer's upcard.

\paragraph{Round Deck} 
The round deck \(\mathbf{g} \in \mathbb{R}^{10}\) represents the composition of the deck during a round of play, where each component \(g_i\) indicates the number of cards of value \(i\) available, following the same interpretation as for the origin deck \(\mathbf{d}\). The \textit{round deck space associated with an origin deck} \(\mathbf{d} \in \Omega_{(\tau, k)}\), denoted as \(\mathcal{G}^\mathbf{d}\), is defined as:
\[
\mathcal{G}^\mathbf{d} = \left\{ \mathbf{g} = (\mathbf{d} - \mathbf{j}_1 - \mathbf{j}_2 - \mathbf{c}) \ \middle|\ g_i \geq 0, \ \mathbf{d} \in \Omega_{(\tau, l)}, \ \mathbf{j}_1, \mathbf{j}_2 \in U^J, \ \mathbf{c} \in U^C \right\}.
\]
The general round deck space is \(\mathcal{G} = \bigcup_{\mathbf{d} \in \Omega_{(\tau, l)}} \mathcal{G}^\mathbf{d}\).

\paragraph{Decision Variables} 
The variables \(z \in \{0,1\}\) and \(q \in \{1,2\}\) indicate whether the player has doubled their bet (\(z=1\)) and which of their hands (\(\mathbf{j}_1\) or \(\mathbf{j}_2\)) is active, respectively.

\subsubsection{Game Instances} 
The game instances \(S\) constitute the set of all possible instances that may occur during a Blackjack session:
\[
S = \left\{ s = (\mathbf{g}, \mathbf{j}_1, \mathbf{j}_2, \mathbf{c}, z, q) \ \middle|\ \mathbf{g} \in \mathcal{G}, \ \mathbf{j}_1, \mathbf{j}_2 \in U^J, \ \mathbf{c} \in U^C, \ z \in \{0,1\}, \ q \in \{1,2\} \right\}.
\]

\subsubsection{Round Instances} 
The round instances \(S^\mathbf{d}\) are restricted to a specific origin deck \(\mathbf{d} \in \Omega_{(\tau, k)}\):
\[
S^\mathbf{d} = \left\{ s \in S \ \middle|\ \mathbf{g} \in \mathcal{G}^\mathbf{d} \right\}.
\]

\subsubsection{Decision Instances} 
The decision instances \(S^J\) are game instances prior to the dealer's turn, where the player makes decisions:
\[
S^J = \{ s \in S \ |\ |\mathbf{c}| \leq 1 \}.
\]
We define the following subsets:
\begin{itemize}
    \item \( S^{(J,\mathbf{d})} = S^J \cap S^\mathbf{d} \): Decision instances associated with an origin deck \( \mathbf{d} \).
    \item \( S^{(J, r)} = \{ s \in S^J \ |\ \mathbf{j_1} \in U^{(J,r)} \} \): Decision instances where the player's hand is relevant.
    \item \( S^{(J, r,\mathbf{d})} = S^{(J, r)} \cap S^\mathbf{d} \): Terminal instances associated with an origin deck \( \mathbf{d} \), where the dealer's hand is relevant.
\end{itemize}

\subsubsection{Terminal Instances} 
The terminal instances, denoted as \(S^F\), represent the instances that mark the end of a round. We will assume, without loss of generality, that the dealer always completes their hand, regardless of the player's hands:
\[
S^F = \{ s \in S \mid \mathbf{c} \text{ has more than 16 points} \}.
\]
We define the following subsets:
\begin{itemize}
    \item \( S^{(F,\mathbf{d})} = S^F \cap S^\mathbf{d} \): Terminal instances associated with an origin deck \( \mathbf{d} \).
    \item \( S^{(F, r)} = \{ s \in S^F \ |\ \mathbf{c} \in U^{(C,r)} \} \): Terminal instances where the dealer's hand is relevant.
    \item \( S^{(F, r,\mathbf{d})} = S^{(F, r)} \cap S^\mathbf{d} \): Decision instances associated with an origin deck \( \mathbf{d} \), where the player's hand is relevant.
\end{itemize}

\subsubsection{Stage Instances} 
The stage instances, denoted as \(S_n\), represent the set of possible instances that can occur exactly after \(n\) stages from the beginning of a round. In this context, a stage of the round is defined as each occasion in which the player makes a decision or a new card is dealt from the deck. Formally:
\[
S_n = \{ s \in S \mid n = \lvert \mathbf{c} \rvert + \lvert \mathbf{j_1} \rvert + q - 1 \}.
\]

\subsubsection{Round Policy}

The player's decision-making process during a round is defined by the round policy function \(\theta\), which assigns a decision to each instance \(s \in S^J\):
\[
\theta: S^J \rightarrow \mathcal{A}, \quad s \mapsto \theta(s)
\]

Where:
\begin{itemize}
    \item \(\mathcal{A} = \{\text{Hit, Stand, Double, Split}\}\): Decision space for the round policy.
\end{itemize}

Technically, in certain instances \(s \in S^J\), the player is not required to make decisions. However, the following relaxations are assumed for these instances to simplify modeling and methodology development:

\begin{enumerate}
    \item If the player has busted with their first and only hand, or with both hands, it is assumed that the player makes a decision but can only stand and end their turn.
    \item If the player obtained a Blackjack or doubled their bet in the previous turn, it is assumed that the player makes a decision but can only stand and end their turn.
    \item If the player's hands have fewer than 2 cards, it is assumed that the player is not obligated to take additional cards by default and can also stand and end their turn.
\end{enumerate}

If the player is indeed required to make a decision in an instance \(s \in S^J\), then the policy \(\theta\) must adhere to the previously established set of rules.

\subsubsection{Game Round and Transition Function}

A \textbf{game round} unfolds by dealing cards to the player's and dealer's hands sequentially, according to the player's decisions and the rules of the game. Each round is represented as a sequence of random instances, where the initial instance is determined by an origin deck \(\mathbf{d}\), and the transition probabilities between instances depend on a round policy \(\theta\):

\[
\mathcal{H}_\theta^\mathbf{d} = \{s_t\}_{0 \leq t}
\]

Where:
\begin{itemize}
    \item \(s_t \in S_t^\mathbf{d}\): The instance of the round after \(t\) stages from the origin deck \(\mathbf{d}\).
    \item \(\theta\): The round policy that determines decisions in the instances \(s_t\).
    \item \(s_0 = (\mathbf{g} = \mathbf{d}, \mathbf{j}_1 = \vec{0}, \mathbf{j}_2 = \vec{0}, \mathbf{c} = \vec{0}, z = 0, q = 1)\): The initial instance of the round.
    \item \(\mathbf{d} \in \Omega_{(\tau,l)}\): The origin deck of the round.
\end{itemize}

Two key functions are highlighted regarding the transition between instances during a round: the transition function and the absolute transition probability function.

\paragraph{Transition Function}

This function defines the probability of transitioning from an instance \(s\) to another instance \(s'\) in one stage, given that a decision \(a\) is taken:

\[
T_a(s', s) = \begin{cases} 
\mathbb{P}(s_{t+1} = s' \mid s_t = s, \theta(s) = a), & \text{if } s \in S^J, \\
\mathbb{P}(s_{t+1} = s' \mid s_t = s), & \text{if } s \notin S^J.
\end{cases}
\]

If the variables of instance \(s'\), the variables of instance \(s\), and the decision \(a\) are consistent with the dynamics of the game, the transition function is calculated simply as the probability of occurrence of the card that differentiates the round decks of both instances:

\[
T_a(s', s) = \frac{g - g'}{\sum g_i}, \quad \text{with } s \in H_\theta^\mathbf{d}.
\]

\paragraph{Absolute Transition Probability}

This function defines the probability of transitioning from an instance \(s\) to another instance \(s'\), independent of the number of stages required, under a policy \(\theta\):

\[
Z_\theta(s', s) = \sum_{k > t} \mathbb{P}(s_k = s' \mid s_t = s, \theta).
\]
\subsubsection{Round Return}

The variable of interest is the final return of the round, determined by the terminal instance. We introduce the \textit{instance return} \( Y_\theta^{s} \in R \), defined as the final return of a round \( \mathcal{H}_\theta^\mathbf{d} \) starting from an instance \( s \in S^\mathbf{d} \), under a round policy \( \theta \):

\[
Y_\theta^{s} = (Y_\theta^{s_t} \mid s_t = s) =
\begin{cases}
Y_\theta^{s_{t+1}},  & \text{if } s \notin S^{(F,\mathbf{d})}, \\
f_R(s),  & \text{if } s \in S^{(F,\mathbf{d})}.
\end{cases}
\]

\noindent Where:
\begin{enumerate}
    \item The function \( f_R \) calculates the player's final return relative to their initial bet, based on a terminal instance \( s \in S^F \):
    \[
    f_R: S^F \rightarrow R, \quad s \mapsto f_R(s) = r_i.
    \]
    \item The space of possible values for the final return relative to the initial bet is referred to as the \textit{return space}, and is represented by the set \( {R} \):
    \[
    R = \{r_1, r_2, r_3, r_4, r_5, r_6\} = \{-2, -1, 0, 1, 1.5, 2\}.
    \]
\end{enumerate}

When referring to the return from the origin instance \( Y_\theta^{s_0} \), the instance \( s_0 \) is fully determined by the origin deck of the round \( \mathbf{d} \in \Omega_{(\tau,l)} \). Therefore, this special case is represented as \( X_\theta^\mathbf{d} = Y_\theta^{s_0} \), and is referred to as the \textit{round return}.

\subsubsection{PMF of the Return: Mass Vector}

The return \( Y_\theta^s \) takes discrete values, and its distribution can be represented by the \textit{mass vector} \( M_\theta^s \in \mathbb{R}^6 \). This vector contains the probabilities associated with the possible final returns \( r_i \in R \), starting from a round instance \( s \in S^{\mathbf{d}} \), under a round policy \( \theta \). Each component of the vector is defined as:
\[
(M_\theta^s)_i = \mathbb{P}(Y_\theta^s = r_i).
\]

The calculation of the mass vector \( M_\theta^s \) from an instance \(s\) can be formulated both recursively and directly:

\paragraph{Recursive Formulation}
\[
M_\theta^s = 
\begin{cases} 
\displaystyle\sum\limits_{s' \in S_{t+1}^\mathbf{d}} T_{\theta(s)}(s', s) \cdot M_\theta^{s'}, & \text{if } s \notin S^{(F,\mathbf{d})}, \\[10pt]
\mathbf{e}_i \cdot 1_{\{f_R(s) = x_i\}}, & \text{if } s \in S^{(F,\mathbf{d})}.
\end{cases}
\]

\paragraph{Direct Formulation}
\[
M_\theta^s = 
\displaystyle\sum\limits_{s' \in S^{(F,\mathbf{d})}} Z_\theta(s', s) \cdot \mathbf{e}_i \cdot 1_{\{f_R(s') = x_i\}}.
\]

Where:
\begin{itemize}
    \item \( \mathbf{e}_i \in \mathbb{R}^6 \): Canonical vector representing a standard basis in \(\mathbb{R}^6\).
\end{itemize}

Similar to the round return \( X_\theta^\mathbf{d} \), the mass vector associated with the initial instance \( s_0 \) is denoted as \( M_\theta^\mathbf{d} = M_\theta^{s_0} \).

\subsubsection{Optimization Problem}

The optimization problem for the round policy consists of maximizing the expected return of a game round. This problem can be addressed independently for each origin deck and is formulated as follows:

\[
\theta^* = \arg \max_\theta \{\mathbb{E}[X_\theta^\mathbf{d}]\}, \quad \forall \mathbf{d} \in \Omega_{(\tau,l)}.
\]

The problem can be formulated recursively, from the terminal instances to the initial instance of the round:

\[
\theta^*(s_t = s)
\;=\;
\arg \max_{a \in \mathcal{A}}
\sum\limits_{s' \in S_{t+1}^\mathbf{d}}
T_a(s', s)\,V^*(s'), \quad \forall s \in S^\mathbf{d}.
\]

\noindent Where:
\begin{enumerate}
    \item \(V^*(s)\) represents the \textbf{optimal value} of an instance \(s\), calculated as follows:
    \[
    V^*(s_t = s) \;=\;
    \begin{cases}
    \displaystyle\max_{a \in \mathcal{A}}
    \sum\limits_{s' \in S_{t+1}^\mathbf{d}}
    T_a(s', s)\,V^*(s'),
    & \text{if } s \notin S^{(F,\mathbf{d})}, \\[10pt]
    f_{\mathcal{R}}(s),
    & \text{if } s \in S^{(F,\mathbf{d})}.
    \end{cases}
    \]
\end{enumerate}

\subsubsection{Basic Strategy}

A widely recognized round policy in the community, and relevant to this study, is the \textbf{Basic Strategy}. This policy considers as state variables the player's hand score, an indicator of whether the hand is soft or hard, the possibility of splitting the hand, and the dealer's visible score. The Basic Strategy is derived by maximizing the expected return of a round with a complete origin deck \( \bar{\mathbf{d}} \) and is denoted as \( \theta^{\text{Basic}} \).
\subsection{Betting Policy}

In this section, the optimization problem for the betting policy is formally formulated as a Markov Decision Process. This problem is defined for a betting session consisting of a total of \( H \) rounds and utilizes a utility function \( u \) to model the player's risk profile.

\subsubsection{Evolution of the Origin Deck}

During a session of consecutive game rounds, the origin deck evolves dynamically throughout the rounds. Formally, this evolution is represented as a sequence of random origin decks, where the initial deck corresponds to the full deck \(\bar{\mathbf{d}}\), and the probabilities of transition between decks depend on a round policy \(\theta\):

\[
\mathcal{D}_\theta = \{\mathbf{d}_n\}_{n=0}^{\infty}
\]

Where:
\begin{itemize}
    \item \( \mathbf{d}_n \in \Omega_{(\tau,l)} \): Origin deck of the \(n\)-th round.
    \item \( \mathbf{d}_0 = \bar{\mathbf{d}} \): Initial deck of the session.
\end{itemize}

\subsubsection{Evolution of Wealth}

Let \( \mathcal{D}_\theta = \{\mathbf{d}_n\}_{n=0}^{\infty} \) be the sequence of origin decks for sequential rounds under a round policy \( \theta \). Let \( W_n \) represent the player's wealth before the start of round \( n \), and \( \pi_n \) the bet placed before the round, expressed as a fraction of \( W_n \). Based on the round return \( X_\theta^{\mathbf{d}_n} \), the player's wealth evolves according to the recursive relation:

\[
W_{n+1} = (1 + \pi_n \cdot X_\theta^{\mathbf{d}_n}) \cdot W_n
\]

Where:
\begin{itemize}
    \item \( \mathbf{d}_n \): Origin deck of the \(n\)-th round.
    \item \( \pi_n \): Bet placed before the \(n\)-th round.
    \item \( W_n \): Wealth before the \(n\)-th round.
\end{itemize}

\subsubsection{States}

The state variables during the betting session must contain all the information necessary to control the wealth value over time. From the wealth equation, it follows that the state during a betting session, denoted as \( \psi \), is fully identified by the origin deck of the round, the current wealth, and the number of rounds played since the session began. The state space \( \mathcal{B} \) is defined as:

\[
\mathcal{B} = \{\psi = (\mathbf{d}, w, n) \mid \mathbf{d} \in \Omega_{(\tau,l)}, w \in \mathbb{R}_+, n \in \mathbb{Z}_+\}
\]

\subsubsection{Betting Policy}

The player's decision-making process during a sequence of game rounds is defined by the function \( \pi \), referred to as the \textbf{complete policy}, which assigns a bet to each state \( \psi \in \mathcal{B} \). The possible values for the bet are limited to a maximum of 50\% of the wealth, ensuring the possibility of doubling or splitting in any round:

\[
\pi: \mathcal{B} \rightarrow [0, 0.5], \quad \psi \mapsto \pi(\psi)
\]

\subsubsection{Betting Session and Transition Function}

A \textbf{betting session} takes place as the player participates in sequential game rounds. Formally, it is represented as a sequence of random states \( \{\psi_{n}=(\mathbf{d}_n, W_n, n)\}_{0 \leq n} \), where the initial state is given by \( \psi_{0}=(\bar{\mathbf{d}}, 1, 0) \). The probability of transitioning from a state \( \psi \) to another state \( \psi' \) after placing a bet \( b \) is defined by the transition function \( T_b(\psi', \psi) \):

\[
T_b(\psi', \psi) = \mathbb{P}(\psi_{n+1} = \psi' \mid \psi_n = \psi,   \pi(\psi) = b) = \begin{cases} 
\mathbb{P}(\mathbf{d}_{n+1} = \mathbf{d}', X_\theta^{\mathbf{d}} = x \mid \mathbf{d}_n = \mathbf{d}), & \text{if } \psi' \in \Gamma_{[\psi, b]} \\
0, & \text{if } \psi' \notin \Gamma_{[\psi, b]}
\end{cases}
\]

\noindent Where:
\begin{enumerate}
    \item \( \Gamma_{[\psi, b]} \) represents the set of states \( \psi' \in \mathcal{B} \) to which the player can transition from a state \( \psi \) after placing a bet \( b \), referred to as the \textit{available space}:
    \[
    \Gamma_{[\psi, b]} = \{(\mathbf{d}^*, w^*, n+1) \mid \mathbf{d}^* \in (\mathcal{G}^\mathbf{d} \cup \bar{\mathbf{d}}), w^* = w \cdot (1 + b \cdot r), r \in R \}
    \]

    \item \( \mathbb{P}(\mathbf{d}_{n+1} = \mathbf{d}', X_\theta^{\mathbf{d}} = x \mid \mathbf{d}_n = \mathbf{d}) \) refers to the probability of ending a round with a return \( x \) and a round deck \( \mathbf{d}' \), after following a round policy \( \theta \) from an origin deck \( \mathbf{d} \):
    \[
    \mathbb{P}(\mathbf{d}_{n+1} = \mathbf{d}', X_\theta^{\mathbf{d}} = x \mid \mathbf{d}_n = \mathbf{d}) = 
    \begin{cases} 
    \displaystyle\sum\limits_{s' \in S^{(F, \mathbf{d})} \mid f_R(s') = x, \mathbf{g}' = \mathbf{d}'} Z_\theta(s', s_0), & \text{if } \mathbf{d}' \neq \bar{\mathbf{d}} \\[10pt]
    \displaystyle\sum\limits_{s' \in S^{(F, \mathbf{d})} \mid f_R(s') = x, \mathbf{g}' \notin \Omega_{(\tau,l)}} Z_\theta(s', s_0), & \text{if } \mathbf{d}' = \bar{\mathbf{d}}
    \end{cases}
    \]
\end{enumerate}

Where:
\begin{itemize}
    \item \( \psi = (\mathbf{d}, w, n) \in \mathcal{B} \): Current state.
    \item \( \psi' = (\mathbf{d'}, w', n+1) \in \mathcal{B} \): Transitioned state.
    \item \( x = \frac{w' - w}{w \cdot b} \in R \): Return associated with the transition.
    \item \( s_0 = (\mathbf{g} = \mathbf{d}, \vec{\mathbf{0}}, \vec{\mathbf{0}}, \vec{\mathbf{0}}, 0, 1) \): Origin instance of the round.
    \item \( s' = (\mathbf{g}', \mathbf{j}_1', \mathbf{j}_2', \mathbf{c}', z, q) \in S^{(F, \mathbf{d})} \): Terminal instance of the round.
\end{itemize}

\subsubsection{Optimization Problem}

In a betting session, the player participates in \( H \) game rounds following a round policy \( \theta \), placing bets that influence the evolution of their wealth \( W_n \). At the end of the session, the player retires with a final wealth \( W_H \). We assume that the player participates in all rounds, regardless of whether they decide to place a bet or not.

The optimization problem for the complete policy consists of maximizing the expected utility of the final wealth \( W_H \) at the end of \( H \) rounds, according to a utility function \( u \), and under a round policy \( \theta \):

\[
\pi_{(\theta, u, H)}^* = \arg \max_\pi \left\{ \mathbb{E}[u(W_H)] \right\}
\]

The problem can be formulated recursively from the terminal states to the initial state of the betting session:

\[
\pi^*(\psi_n = \psi)
\;=\;
\arg \max_{b \in [0, 0.5]}
\sum\limits_{\psi' \in \Gamma_{[\psi, b]}}
T_b(\psi', \psi)\,V^*(\psi'), \quad \forall \psi \in \mathcal{B} \mid {n<H}
\]

\noindent Where:
\begin{enumerate}
    \item \(V^*(\psi)\) represents the \textbf{optimal value} of a state \(\psi\), and is calculated as follows:
    \[
    V^*(\psi_{n} = \psi) \;=\;
    \begin{cases}
    \displaystyle\max_{b \in [0, 0.5]}
    \sum\limits_{\psi' \in \Gamma_{[\psi, b]}}
    T_b(\psi', \psi)\,V^*(\psi'),
    & \text{if } n<H
    \\[10pt]
    u(w),
    & \text{if } n=H
    \end{cases}
    \]
\end{enumerate}

Where:
\begin{itemize}
    \item \( \psi = (\mathbf{d}, w, n) \in \mathcal{B} \).
\end{itemize}

In this study, two utility functions are considered:

\paragraph{Utility with CRRA:} Represents a Constant Relative Risk Aversion (CRRA) of \( 1 - \alpha \). Denoted as \( u_1 \), it is defined as:
\[
u_1(w; \alpha) = 
\begin{cases} 
\frac{w^\alpha}{\alpha}, & \text{if } \alpha > 0, \\
\ln(w), & \text{if } \alpha = 0,
\end{cases}
\quad \text{with } \alpha \in \mathbb{R}_+
\]

\paragraph{Utility with CARA:} Represents a Constant Absolute Risk Aversion (CARA) of \( \beta \). Denoted as \( u_2 \), it is defined as:
\[
u_2(w; \beta) = 1 - e^{-\beta w}, \quad \text{with } \beta \in \mathbb{R}
\]
\subsection{Utility with CRRA: Optimality Conditions for Infinite Rounds}

This section establishes the necessary conditions for the betting policy to be optimal in an infinite horizon of rounds (\( H \to \infty \)), considering a utility function with constant relative risk aversion (CRRA). Formally, the optimization problem for the betting policy is defined as:

\[
\pi_{(\theta, u_1, H \to \infty)}^* = \arg \max_\pi \{\mathbb{E}[u_1(W_{H \to \infty}; \alpha)]\}
\]

Where:
\begin{itemize}
    \item \( u_1 \): Utility function with CRRA, parameterized by the coefficient \( \alpha \).
    \item \( W_{H \to \infty} \): Player's wealth in the infinite horizon.
\end{itemize}

\subsubsection{Policy Simplification}

To solve this problem, we rely on Theorems 1  and 2, which are proven in Appendices \ref{teo1} and \ref{teo2}. These results establish the following fundamental properties of the optimal policy \( \pi_{(\theta, u_1, H \to \infty)}^* \):

\begin{itemize}
    \item \textbf{Theorem 1 (\ref{teo1}):} The optimal policy \( \pi_{(\theta, u_1, H \to \infty)}^* \) is independent of the player's wealth \( w \).
    \item \textbf{Theorem 2 (\ref{teo2}):} The optimal policy \( \pi_{(\theta, u_1, H \to \infty)}^* \) is independent of the number of rounds played \( n \).
\end{itemize}

Both theorems allow us to conclude that the optimal policy for an infinite horizon of rounds under CRRA depends exclusively on the origin deck \( \mathbf{d} \) of each round. Mathematically, this is expressed as:

\[
\pi_{(\theta, u_1, H \to \infty)}^*(\mathbf{d}, w, n) = \pi_{(\theta, u_1, H \to \infty)}^*(\mathbf{d}, w', n'), \quad \forall w, w' \in \mathbb{R}_+, \forall n, n' \in \mathbb{Z}_+
\]

\subsubsection{Counting System}

Generalizing the analysis, we assume the player does not know the exact composition of the origin deck \( \mathbf{d} \), but only has partial information about it through a counting system. Formally, a counting system \( \kappa \) is a function that calculates a value or vector \( c \), referred to as the \textit{count}, from the composition of the origin deck \( \mathbf{d} \):

\[
\kappa: \Omega_{(\tau,l)} \rightarrow \text{Im}(\kappa) \subseteq \mathbb{R}^n, \quad \mathbf{d} \mapsto \kappa(\mathbf{d}) = c
\]

Where:
\begin{itemize}
    \item \( 1 \leq n \leq 10 \): Dimension of the count vector.
    \item \( \text{Im}(\kappa) = \{\kappa(\mathbf{d}) \mid \mathbf{d} \in \Omega_{(\tau,l)}\} \): Image of the counting system.
\end{itemize}

There are multiple counting systems. In this study, two relevant systems are considered:

\paragraph{Perfect Count:}
This system manages complete information about the origin deck and is denoted as \( \kappa_0 \):

\[
\kappa_0(\mathbf{d}) = c = \mathbf{d}, \quad \text{Im}(\kappa_0) = \Omega_{(\tau,l)}
\]

\paragraph{True Count:}
A widely used system due to its simplicity and effectiveness. It is denoted as \( \kappa_1 \) and is defined as follows:

\[
\kappa_1(\mathbf{d}) = c = \left\lfloor \frac{52 \cdot \left(\sum_{j=2}^{6} (\bar{\mathbf{d}} - \mathbf{d})_j - (\bar{\mathbf{d}} - \mathbf{d})_1 - (\bar{\mathbf{d}} - \mathbf{d})_{10} \right)}{\sum (\mathbf{d})_i} \right\rfloor, \quad \text{Im}(\kappa_1) = \{n \in \mathbb{Z} \mid -52 \leq n \leq 52\}
\]

\subsubsection{Generalized Problem: Partial Policy}

A generalized policy, referred to as a \textit{partial policy}, is introduced, which uses as a state variable the count associated with the origin deck of the round, derived from a counting system \( \kappa \):

\[
\omega_\kappa: \text{Im}(\kappa) \rightarrow [0, 0.5], \quad c \mapsto \omega_\kappa(c)
\]

The optimization problem for a partial policy is defined in terms of a counting system \( \kappa \), a round policy \( \theta \), a utility function with CRRA \( u_1 \), and an infinite betting horizon \( H \to \infty \):

\[
\omega_{(\kappa, \theta, u_1, H \to \infty)}^* = \arg \max_\omega \{\mathbb{E}[u_1(W_{H \to \infty}; \alpha)]\}
\]

\subsubsection{Optimality Conditions}

Theorem 4, proven in Appendix \ref{teo4}, establishes that the optimality condition for the partial policy can be formulated independently for each count \( c \), as expressed below:

\[
\omega_{(\kappa, \theta, u_1, H \to \infty)}^*(c) = \arg \max_{b \in [0, 0.5]} \{\mathbb{E}[u_1(1 + b \cdot X_\theta^{(\mathbf{D}|c)}; \alpha)]\}, \quad \forall c \in \text{Im}(\kappa)
\]

\noindent Where:
\begin{enumerate}
    \item The variable \( \mathbf{D} = \mathbf{D}_\theta \in \Omega_{(\tau, l)} \) represents the stationary deck of the evolution of the origin deck \( \mathcal{D}_\theta = \{\mathbf{d}_n\}_{n=0}^{\infty} \), and its PMF is defined as follows:
    \[
    \mathbb{P}_\mathbf{D}(\mathbf{d}) = \mathbb{P}(\mathbf{D_\theta} = \mathbf{d}) = \lim_{N \to \infty} \frac{\sum_{n=0}^N 1_{\{\mathbf{d}^{obs}_n = \mathbf{d}\}}}{N+1}
    \]

    \item The variable \( \mathbf{D} \mid c \) represents the stationary deck conditioned on the count associated with that deck being \( c \). Its PMF is expressed as:
    \[
    \mathbb{P}_{(\mathbf{D} | c)}(\mathbf{d})
    = \frac{\mathbb{P}(\mathbf{D} = \mathbf{d}, \, \kappa(\mathbf{D}) = c)}
           {\mathbb{P}(\kappa(\mathbf{D}) = c)}
    \]
\end{enumerate}

Where \( \mathcal{D}_\theta^{obs} = \{\mathbf{d}^{obs}_n\}_{n=0}^{\infty} \): A specific realization of the deck sequence, and each \(\mathbf{d}^{obs}_n\) corresponds to a particular realization of the random deck \(\mathbf{d}_n\).

\subsubsection{Return Distribution}

Theorem 3 (see proof in Appendix \ref{teo3}) provides an approximation for the distribution of \( W_n \) under a partial policy, converging as \( n \to \infty \):

\[
\lim_{n \to \infty} W_{n+1} \sim \text{LogNormal}(n \cdot \mu_\omega, n \cdot \sigma_\omega^2)
\]

\noindent Where:
\begin{enumerate}
    \item The variables \(\mu_\omega=\mathbb{E}[L_\omega]\) and \(\sigma_\omega^2 = \text{Var}[L_\omega]\) are associated with the random variable \(L_\omega\). The latter is conditioned by the partial policy \(\omega\) and is defined as:
    \[
    L_\omega = \ln(1 + \omega(C) \cdot X_\theta^{(\mathbf{D} \mid C)})
    \]
    \item The variable \( C = C_{( \theta, \kappa)} = \kappa(\mathbf{D}_\theta) \in \text{Im}(\kappa) \) represents the stationary count of the evolution of the count \( \mathcal{C}_{(\theta, \kappa)} = \{\kappa({\mathbf{d}_n})\}_{n=0}^{\infty} \), and its PMF is defined as follows:
    \[
    \mathbb{P}_C(c) = \mathbb{P}(C_{( \theta, \kappa)} = c) = \lim_{N \to \infty} \frac{\sum_{n=0}^N 1_{\{\kappa(\mathbf{d}^{obs}_n) = c\}}}{N+1}
    \]
\end{enumerate}

Where  \( \mathcal{C}^{obs}_{(\theta, \kappa)} = \{\kappa(\mathbf{d}^{obs}_n)\}_{n=0}^{\infty} \): A specific realization of the count sequence, and each \(\kappa(\mathbf{d}^{obs}_n)\) corresponds to a particular realization of the random count \(\kappa(\mathbf{d}_n)\).

\subsubsection{Corollary}
\begin{enumerate}
    \item For \( \alpha = 1 \), when the utility equals the value of the wealth, the following holds:
    \[
    \omega_{(\kappa, \theta, u_1, H \to \infty)}^*(c) = \begin{cases} 
    0.5, & \text{if } \mathbb{E}[X_\theta^{(\mathbf{D}|c)}] > 0, \\
    0, & \text{if } \mathbb{E}[X_\theta^{(\mathbf{D}|c)}] \leq 0,
    \end{cases} \quad \forall c \in \text{Im}(\kappa)
    \]
    
    \item The optimal partial policy bets only if the expected return associated with the count is positive:
    \[
    \omega_{(\kappa, \theta, u_1, H \to \infty)}^*(c) = \begin{cases} 
    0.5, & \text{if } \mathbb{E}[X_\theta^{(\mathbf{D}|c)}] > 0, \\
    0, & \text{if } \mathbb{E}[X_\theta^{(\mathbf{D}|c)}] \leq 0,
    \end{cases} \quad \forall c \in \text{Im}(\kappa)
    \]  

    \item For a perfect counting system \( \kappa_0(\mathbf{d}) = \mathbf{d} \), the optimal complete policy is obtained:
    \[
    \pi_{(\theta, u_1, H \to \infty)}^*(\mathbf{d}, w, n) = \omega_{(\kappa_0, \theta, u_1, H \to \infty)}^*(c)
    \]
\end{enumerate}

\section{Methodology}

This section describes the methodology used to address the optimization of the round policy and the betting policy. Both problems are approached approximately, culminating in an estimation of the optimal policies. Subsequently, the performance of the returns of a player following the calculated policies is evaluated using a simulated game environment. Finally, histograms of the wealth obtained at the end of the betting session are generated.

\subsection{Round Policy}

The optimization problem for the round policy can be directly approached using dynamic programming, which would allow for the exact computation of both the optimal round policy and the probability mass function (PMF) of the return for each round \cite{peytavi2018matematicas}. Although conceptually simple, its practical implementation results in extremely high computational costs. This is mainly due to the option of performing a "Split," which leads to exponential growth in the state space and transitions between them, causing the computation to take days or even weeks to complete.

Even excluding the possibility of performing a "Split," the direct implementation of dynamic programming remains inefficient. However, it is possible to identify patterns in both the states and transitions that significantly reduce the required computations. The proposed procedure for efficiently addressing the optimization problem for the round policy is structured into the following steps:

\begin{enumerate}
    \item Mathematically formulate the optimization problem for \textit{restricted round policies} \(\theta_i\), limiting the decision space:
    \begin{itemize}
        \item \(\theta_1 \in \{\text{Hit, Stand, Double}\}\)
        \item \(\theta_2 \in \{\text{Hit, Stand}\}\)
    \end{itemize}
    
    \item Compute the mass vector \(M_{\theta_1}^s\), conditioned on \(\theta_1(s) = \text{Stand}\), for all decision instances of the round.
    
    \item Compute the optimal restricted policies \(\theta_i^*(s)\) and the optimal mass vectors \(M_{\theta_i^*}^s\) using dynamic programming, for all decision instances of the round.
    
    \item Estimate the complete optimal policy \(\theta^*\), referred to as the \textit{semi-optimal policy} and denoted as \(\tilde{\theta}^*\), along with an estimation of the round's mass vector \(\tilde{M}_{\tilde{\theta}^*}^\mathbf{d}\), considering an approximation of the mass vector \(M_\theta^s|\theta(s) = \text{Split}\) for the corresponding decision instances.
\end{enumerate}

\subsubsection{Restricted Problem}

The restricted problem is formulated by limiting the decision space for the round policy. Two possible restrictions are defined:

\begin{itemize}
    \item \(i=1\): Splitting is not allowed.
    \item \(i=2\): Splitting and Doubling are not allowed.
\end{itemize}

Depending on the restriction applied to the problem, the state space of the game is simplified as follows:

\paragraph{Restricted Game Instances}
\[
\prescript{}{i}{S}  = 
\begin{cases} 
\{s = (\mathbf{g}, \mathbf{j_1}, \mathbf{c}, z) \in \mathcal{G}_R \times U^J \times U^C \times \{0, 1\}\}, & \text{if } i = 1, \\
\{s = (\mathbf{g}, \mathbf{j_1}, \mathbf{c}) \in \mathcal{G}_R \times U^J \times U^C\}, & \text{if } i = 2.
\end{cases}
\]
Where:
\begin{itemize}
    \item \(\mathcal{G}_R = \bigcup_{\mathbf{d} \in \Omega_{(\tau, l)}} \mathcal{G}_R^\mathbf{d}  \mid \mathcal{G}_R^\mathbf{d} = \{\mathbf{g} = (\mathbf{d} - \mathbf{j_1} - \mathbf{c}) \mid g_i \geq 0, \mathbf{j}_1 \in U^J, \mathbf{c} \in U^C\}\): Restricted round decks.
\end{itemize}

From the set \(\prescript{}{i}{S}\), it is possible to define, analogously to the definitions presented in the theoretical framework, all subsets of restricted instances. This includes restricted round instances \(\prescript{}{i}{S}^\mathbf{d}\), restricted decision instances \(\prescript{}{i}{S}^J\), restricted terminal instances \(\prescript{}{i}{S}^F\), restricted turn instances \(\prescript{}{i}{S}_n\), as well as analogs of all other subsets described in the theoretical framework.

The restricted policies \(\theta_i\) take as input restricted decision instances \(s \in \prescript{}{i}{S}^J\), and are defined as follows:
\[
\theta_i: \prescript{}{i}{S}^J \rightarrow \prescript{}{i}{\mathcal{A}}=
\begin{cases} 
\{\text{Hit, Stand, Double}\}, & \text{if } i = 1, \\
\{\text{Hit, Stand}\}, & \text{if } i = 2.
\end{cases}
\]

Finally, the restricted optimization problem is formulated for any origin deck \(\mathbf{d} \in \Omega_{(\tau, l)}\):

\begin{equation}
\theta_i^*(s_t = s)
\;=\;
\arg \max_{a \in \prescript{}{i}{\mathcal{A}}}
\sum\limits_{s' \in \prescript{}{i}{S}_{t+1}^\mathbf{d}}
T_a(s', s)\,V_i^*(s'), \quad \forall s \in \prescript{}{i}{S}^\mathbf{d}.
\label{eq:optimizacion_restringida}
\end{equation}

\noindent Where:
\begin{enumerate}
    \item \(V_i^*(s)\) represents the \textbf{optimal value function} of an instance \(s\), calculated as follows:
    \[
    V_i^*(s_t=s) \;=\;
    \begin{cases}
    \displaystyle\max_{a \in \prescript{}{i}{\mathcal{A}}}
    \sum\limits_{s' \in \prescript{}{i}{S}_{t+1}^\mathbf{d}}
    T_a(s', s)\,V_i^*(s'),
    & \text{if } s \notin \prescript{}{i}{S}^{(F,\mathbf{d})},
    \\[10pt]
    f_{\mathcal{R}}(s),
    & \text{if } s \in \prescript{}{i}{S}^{(F,\mathbf{d})}.
    \end{cases}
    \]
\end{enumerate}
\subsubsection{Player and Dealer Sequences}

The first step involves obtaining the data that defines the graphs connecting the player’s and dealer’s sequences. To this end, an algorithm is used to iteratively deal cards, allowing the construction of the following structures:

\[
\mathbf{U}^J, \, \mathbf{A}^J, \, \mathbf{U}^{(J,r)}, \, \mathbf{A}^{(J,r)}, \, \mathbf{A}_T^{(J,r)}, \, \mathbf{U}^C, \, \mathbf{A}^C, \, \mathbf{U}^{(C,r)}, \, \mathbf{A}^{(C,r)}.
\]

In particular, the associated matrices, \( \mathbf{U}^{(J,r)} \in \mathbb{R}^{|U^{(J,r)}| \times 10} \) and \( \mathbf{U}^{(C,r)} \in \mathbb{R}^{|U^{(C,r)}| \times 10} \), structure and organize the sequences in the sets \( U^{(J,r)} \) and \( U^{(C,r)} \), respectively. 

Since the dealer plays their hand according to predetermined rules, the DAG \( G^{C} = (U^{C}, A^{C}) \) follows the same trajectories in each round and ends in the same terminal sequences. The multiplicity of a relevant sequence \( \mathbf{c} \in U^{(C,r)} \) from an upcard \( \mathbf{c}^* \) is defined as the total number of trajectories in the graph \( G^C \) that connect the initial node associated with \( \mathbf{c}^* \) to the terminal node associated with \( \mathbf{c} \):
\[
m(\mathbf{c}, \mathbf{c}^*) =
\begin{cases}
1, & \text{if } \mathbf{c} = \mathbf{c}^*, \\
\sum_{\mathbf{c}' \in \text{Pred}(\mathbf{c})} m(\mathbf{c}', \mathbf{c}^*), & \text{if } \mathbf{c} \neq \mathbf{c}^*,
\end{cases}
\quad \text{where } \text{Pred}(\mathbf{c}) = \{\mathbf{c}' \in U^{(C,r)} : (\mathbf{c}', \mathbf{c}) \in A^{(C,r)}\}.
\]

For more details on the player’s and dealer’s structures, as well as on multiplicities, see Appendix~\ref{sec:secuencias_jugador} (\textit{Player Sequences}) and \ref{sec:secuencias_crupier} (\textit{Dealer Sequences}).

\subsubsection{Sequence Extraction}

A relevant aspect of the methodology for optimizing the round policy is the extraction of a card sequence from a given deck. Specifically, the probability associated with this extraction is of particular importance. Therefore, the probability of extracting a sequence \( c \in \mathbb{R}^{10} \) from a deck \( b \in \mathbb{R}^{10} \) is represented by the function $\zeta$:

\[
\zeta(c, b) =  \prod_{i=1}^{10} \frac{\binom{b_i}{c_i}}{\binom{\sum b_i}{\sum c_i}}
\]

Using Bayes' Theorem, a key property of the function $\zeta$ is derived, which is fundamental for efficiently optimizing the round policy. In particular, this property allows the probability of extracting a sequence \( c \) from a deck \( b - e_i \), where \( \mathbf{e}_i \in \mathbb{R}^{10} \) is a canonical vector representing the extraction of a single card, to be expressed in terms of $\zeta(c, b)$ as follows:

\[
\zeta(c, b - e_i) = \zeta(c, b) \cdot \frac{\zeta(e_i, b - c)}{\zeta(e_i, b)} 
= \zeta(c, b) \cdot \frac{(b - c)_i}{\sum (b - c)_i} \cdot \frac{\sum b_i}{b_i}
\]

\subsubsection{Stand: Mass Vectors}
Under the restricted policy \( \theta_1 \), the player’s turn may opt to stand (\textit{Stand}) at any decision instance \( s \in \prescript{}{1}{S}^{(J, \mathbf{d})} \). When this decision is made, the player’s turn ends, and the remainder of the round proceeds according to predetermined rules. This implies that the mass vector of the return \( (M_{\theta_1}^s \mid \theta_1(s) = \text{Stand}) \) is fully determined once the player decides to stand. Calculating this vector for all decision instances is essential for optimizing \( \theta_1 \), and it constitutes the main computational cost of the problem.

For decision instances where the player exceeds 21 points (\( j_1 \notin U^{(J, r)} \)), the player automatically loses, and the mass vector can be directly deduced. Therefore, the explicit computation of \( (M_{\theta_1}^s \mid \theta_1(s) = \text{Stand}) \) is only necessary for instances \( s \in \prescript{}{1}{S}^{(J, r, \mathbf{d})} \). According to the direct formulation for calculating the mass vector established in the theoretical framework, if \( s \in \prescript{}{1}{S}^{(J, r, \mathbf{d})} \), the calculation of the mass vector can be disaggregated in terms of the relevant terminal instances:

\[
(M_{\theta_1}^s \mid \theta_1(s) = \text{Stand}) = 
\sum_{s' \in \prescript{}{1}{S}^{(F, r, \mathbf{d})}} Z_{\theta_1}(s', s) \cdot \mathbf{e}_i \cdot 1_{\{f_R(s') = x_i\}}
\quad + \quad
\text{DealerBust}_{\theta_1}(s) \cdot \mathbf{e}_{f_{\text{DealerBust}}^s}.
\]

\noindent Where:
\begin{enumerate}
    \item \(\text{DealerBust}_{\theta_1}(s)\) represents the probability that the dealer exceeds 21 points in instance \( s \), ensuring the player’s victory:
    \[
    \text{DealerBust}_{\theta_1}(s) = 1 - \sum_{s' \in \prescript{}{1}{S}^{(F, r, \mathbf{d})}} Z_{\theta_1}(s', s).
    \]

    \item \(Z_{\theta_1}(s',s)\) can be explicitly calculated in terms of the function \(\zeta\) and the multiplicity of the dealer’s terminal sequence:
    \[
    Z_{\theta_1}(s',s) = \frac{\zeta(\mathbf{c'}, \mathbf{d-j_1})}{\zeta(\mathbf{c^*}, \mathbf{d-j_1})} \cdot m(\mathbf{c'}, \mathbf{c^*}),
    \]
    where \(\mathbf{c'} \in U^{(C, r)}\).

    \item \(f_{\text{DealerBust}}^s\) categorizes the return associated with the event that the dealer exceeds 21 points:
    \[
    f_{\text{DealerBust}}^s =
    \begin{cases} 
    4, & \text{if } z = 0 \text{ and } j_1 \text{ is not a Blackjack}, \\
    5, & \text{if } z = 0 \text{ and } j_1 \text{ is a Blackjack}, \\
    6, & \text{if } z = 1.
    \end{cases}
    \footnote{
    It is important to note that the function \(f_{\text{DealerBust}}^s\) depends solely on the decision instance \(s\), meaning its value is exclusively determined by the player’s score in that instance. Specifically, it does not consider the dealer’s exact score beyond the fact that the dealer exceeds 21 points (i.e., a situation of busting). Hence, \(f_{\text{DealerBust}}^s\) is independent of the dealer’s exact score in state \(s'\).}
    \]
\end{enumerate}

\noindent Where:
\begin{itemize}
    \item \( s = (\mathbf{d-c^*-j_1}, \mathbf{j_1}, \mathbf{c^*}, z) \in \prescript{}{1}{S}^{(J, r, \mathbf{d})} \).
    \item \( s' = (\mathbf{d-j_1-c'}, \mathbf{j_1}, \mathbf{c'}, z) \in \prescript{}{1}{S}^{(F, r, \mathbf{d})} \).
\end{itemize}

Based on this formulation, the computation of the mass vector \( (M_{\theta_1}^s \mid \theta_1(s) = \text{Stand}) \) for each decision instance \( s \in \prescript{}{1}{S}^{(J, r, \mathbf{d})} \) primarily relies on calculating the extraction probability \(\zeta(\mathbf{c}, \mathbf{d-j_1})\) and returns \(f_\mathcal{R}(\mathbf{g}, \mathbf{j_1}, \mathbf{c}, z=0)\) for each \(\mathbf{j_1} \in U^{(J,r)}\) and \(\mathbf{c} \in U^{(C,r)}\), leveraging the precomputed structures and multiplicities \(m(\mathbf{c}, \mathbf{c^*})\).

In cases where the bet has been doubled, the final returns \(f_\mathcal{R}(\mathbf{g}, \mathbf{j_1}, \mathbf{c}, z=1)\) are dynamically computed by doubling the precomputed returns for \(z=0\), according to the expression: \( f_\mathcal{R}(\mathbf{g}, \mathbf{j}_1, \mathbf{c}, z = 1) = 2 \cdot f_\mathcal{R}(\mathbf{g}, \mathbf{j}_1, \mathbf{c}, z = 0) \). Additionally, the extraction probabilities associated with the \textit{upcards}, \(\zeta(\mathbf{c^*}, \mathbf{d} - \mathbf{j_1})\), are also required to compute the mass vectors. However, the computation of these probabilities has relatively low complexity, and additional details on their implementation can be found in Appendix~\ref{sec:calculo_secuenciasj_upcards} (\textit{Calculation of Extraction Probabilities for Player Sequences by Upcards}).

\paragraph{Precomputing Final Returns}

The first step involves precomputing the final return \( f_\mathcal{R}(\mathbf{g}, \mathbf{j}_1, \mathbf{c}, z = 0) \) for all combinations of \( \mathbf{j}_1 \in U^{(J,r)} \) and \( \mathbf{c} \in U^{(C,r)} \).

This process generates a matrix of final returns \( \mathbf{R} \) with dimensions \( |U^{(J,r)}| \times |U^{(C,r)}| \), which organizes the values of \( f_\mathcal{R} \). Using the sequence matrices \( \mathbf{U}^{(J,r)} \) and \( \mathbf{U}^{(C,r)} \), the matrix \( \mathbf{R} \) is constructed, defined as:
\(
\mathbf{R}(i,k) = f_\mathcal{R}(\mathbf{g}, \mathbf{U}^{(J,r)}(i), \mathbf{U}^{(C,r)}(k), z = 0),
\) where \( \mathbf{U}^{(J,r)}(i) \) and \( \mathbf{U}^{(C,r)}(k) \) correspond to the \( i \)-th and \( k \)-th sequences of \( U^{(J,r)} \) and \( U^{(C,r)} \), respectively.

For further details on this procedure, including how the return values are defined for each combination \((\mathbf{U}^{(J,r)}(i), \mathbf{U}^{(C,r)}(k))\), refer to Appendix~\ref{sec:retornos_finales} (\textit{Final Return Matrix Calculation}).

\paragraph{Extraction Probability Calculation}

To compute \( Z_{\theta_1}(s', s) \), the process begins by calculating the extraction probability \(\zeta(\mathbf{c}, \mathbf{d} - \mathbf{j}_1)\) for all \(\mathbf{j}_1 \in U^{(J,r)}\) and \(\mathbf{c} \in U^{(C,r)}\). For this purpose, the function \textbf{Probs\_secs()} (Appendix~\ref{sec:anexo3}) and the matrix \(\mathbf{U}^{(C,r)}\), containing the sequences of \(U^{(C,r)}\), are used. This calculation determines the extraction probabilities vectorially from an origin deck \(\mathbf{d}\). As a result, a vector \(\mathbf{Q} \in \mathbb{R}^{|U^{(C,r)}|}\) is generated, defined as:
\(
\mathbf{Q}(k) = \zeta(\mathbf{U}^{(C,r)}(k), \mathbf{d}).
\)

Using the vector \(\mathbf{Q}\), the extraction probability for each sequence \(\mathbf{c} \in U^{(C,r)}\) from adjusted decks of the form \(\mathbf{d} - \mathbf{j}_1\), with \(\mathbf{j}_1 \in U^{(J,r)}\), is calculated recursively. This process involves traversing the transitions of the DAG \(G_T^{(J,r)} = (U^{(J,r)}, A_T^{(J,r)})\) recursively, starting from the matrix \( \mathbf{Q} \) and leveraging the properties of \(\zeta\) to ensure efficiency. This process, implemented using the function \textbf{Matriz\_Probs()} (Appendix~\ref{sec:anexo4}), produces a matrix \(\mathbf{P} \in \mathbb{R}^{|U^{(J,r)}| \times |U^{(C,r)}|}\), defined as:
\(
\mathbf{P}(i,k) = \zeta(\mathbf{U}^{(C,r)}(k), \mathbf{d} - \mathbf{U}^{(J,r)}(i)).
\)

For a detailed description of how this vectorized calculation is implemented and how the transitions of the DAG are traversed, refer to Appendix, Section~\ref{sec:calculo_secuenciasj_upcards} (\textit{Extraction Probability Calculation}).

\paragraph{Mass Vector Calculation for the Stand Decision}

After calculating the objects \(\mathbf{R}\), \(\mathbf{P}\), and the multiplicities \(m(\mathbf{c}, \mathbf{c^*})\), the mass vectors associated with the stand (\textit{Stand}) decision are determined for all decision instances \(s \in \prescript{}{1}{S}^{(J, r, \mathbf{d})}\). To achieve this, instances are grouped by the dealer’s \textit{upcard}, allowing for efficient processing of the relevant combinations of player and dealer sequences. Each group of mass vectors \((M_{\theta_1}^s \mid \theta_1(s) = \text{Stand})\) is calculated simultaneously using Expression 3.1.4, employing vector operations and leveraging the structures \(\mathbf{R}\), \(\mathbf{P}\), and the associated multiplicities \(m(\mathbf{c}, \mathbf{c^*})\). For a detailed description of the procedure, refer to Appendix, Section~\ref{sec:calculo_vectores_plantarse} (\textit{Mass Vector Calculation for the Stand Decision}).
\subsubsection{Hit, Stand, and Double: Exact Solution}

Once the mass vector \( (M_{\theta_1}^s \mid \theta_1(s) = \text{Stand}) \) has been calculated for every decision instance \( s \in \prescript{}{1}{S}^{(J, \mathbf{d})} \), the optimization of the restricted policies \( \theta_1 \) and \( \theta_2 \) is carried out dynamically. At this point, the expected return for the "Stand" action is already available for all decision instances. 

The optimal policy \( \theta_i^*(s) \) is then obtained recursively, as presented earlier in this section (see Equation~\eqref{eq:optimizacion_restringida}):
\[
\theta_i^*(s_t = s)
\;=\;
\arg \max_{a \in \prescript{}{i}{\mathcal{A}}}
\sum\limits_{s' \in \prescript{}{i}{S}_{t+1}^\mathbf{d}}
T_a(s', s)\,V^*(s'), \quad \forall s \in \prescript{}{i}{S}^\mathbf{d}
\]

\noindent The computation is organized according to \( \mathbf{G}^{(J,r)} \), the Directed Acyclic Graph (DAG), where each level groups the player’s sequences \( \mathbf{j}_1 \in U^{(J,r)} \), classified by the number of cards played (\(n\)). These relevant player sequences are traversed within the DAG, combining each sequence with all possible dealer \textit{upcards} \( \mathbf{c^*} \in U^{C} \). This procedure generates all decision instances \( s \), for which the policies \( \theta_1 \) and \( \theta_2 \) are optimized at each stage of the process.

\paragraph{Terminal Level (\(n = N = 21\)):}
At \(n = N\), the player’s sequences correspond to those at the deepest level of the DAG \( \mathbf{G}^{(J,r)} \). These sequences are combined with all dealer upcards \( \mathbf{c^*} \in U^{C} \) to form the decision instances \( s \in \prescript{}{i}{S}_{22}^{(J,r, \mathbf{d})} \). Since the player has exceeded 21 points, they automatically lose in all these instances. Therefore, the optimal decision is \( \theta^*(s) = \text{Stand} \), and the associated mass vectors are directly defined as:
\[
M_{\theta^*}^s = \mathbf{e}_2,
\]

\paragraph{Intermediate Levels (\(n < N\)):}
The recursive computation proceeds from \(n = N - 1\) down to \(n = 0\). At each level, the following are simultaneously calculated for all instances \( s \in \prescript{}{i}{S}_n^{(J,r, \mathbf{d})} \): the value function \( V^*(s) \), the optimal decisions \( \theta_i^*(s) \), and the mass vectors \( M_{\theta_i}^s \). The value function determines \( \theta_i^*(s) \), and these decisions allow the calculation of \( M_{\theta_i}^s \), propagating information from level \( n+1 \).

These calculations are performed vectorially for all instances at a given level, leveraging the structure of the DAG, while the levels are processed sequentially. This procedure ensures consistency and efficiency in propagating information from \( n = N \) back to the initial levels. For a more detailed explanation of the calculation of the mass vectors \( M_{\theta_i}^s \) and the optimal decisions \( \theta_i^*(s) \), the reader is referred to Appendix~\ref{sec:hit_stand_double} (\textit{Hit, Stand, and Double: Exact Solution}).
\subsubsection{Hit, Stand, Double, and Split: Estimated Solution}

We have computed the optimal restricted policy $\theta_2^*$ and the optimal mass vector $M_{\theta_2^*}^s$ for each instance $s \in \prescript{}{2}S_2^{(J, \mathbf{d})}$. When addressing the general problem with the complete policy $\theta$, the possibility of performing a Split involves a high computational cost. Given the need for efficiency in this research, we propose considering the Split option in an approximate manner. The first approximation consists of estimating the state return for the Split action $Y_\theta^s | \theta(s) = \text{Split}$ based on the optimal state return under the restricted policy $Y_{\theta_2^*}^{(s')}$, for instances $s'$ where the player has only one card in hand:

\[
(Y_\theta^s | \theta(s) = \text{Split}) \approx (\tilde{Y}_\theta^s | \theta(s) = \text{Split}) = (Y_1)_{\theta_2^*}^{(s')} + (Y_2)_{\theta_2^*}^{(s')}
\]

\noindent Where:
\begin{enumerate}
    \item The state returns $(Y_i)_{\theta_2^*}^{(s')}$ are independent and identically distributed, with distribution $Y_{\theta_2^*}^{(s')}$, such that the mass vector is known:\\
    \[
    M_{\theta_2^*}^{(s')} = [0, a_{-1}, a_0, a_1, 0, 0]
    \]
    \item The estimated mass vector is calculated as follows:
    \[
    (M_\theta^s | \theta(s) = \text{Split}) \approx (\tilde{M}_\theta^s | \theta(s) = \text{Split}) = [a_{-1}^2, 2a_{-1} a_0, 2a_{-1} a_1 + a_0^2, 2a_1 a_0, 0, a_1^2]
    \]

\end{enumerate}
Where:
\begin{itemize}
    \item $s' = \left( \mathbf{d}- \frac{\mathbf{j_1}}{2}-\mathbf{c^*}, \frac{\mathbf{j_1}}{2}, \mathbf{c^*} \right) \in \prescript{}{2}S_2^{(J, \mathbf{d})}$: Restricted instance associated with dealing the first of the two cards in the player's hand.
    \item $s = (\mathbf{d}- \mathbf{j_1}-\mathbf{c^*}, \mathbf{j_1}, \mathbf{j_2}= \mathbf{0}, \mathbf{c^*}, z, q) \in S^{(J, \mathbf{d})}$: Complete instance, where a Split is possible $j_1 \in \{0, 2\}^{10}, |j_1| = 2$.
\end{itemize}

\textit{Note:} If the split hand is a pair of tens, winning by blackjack is considered impossible. If the split hand is a pair of aces, winning by blackjack and receiving more than one card for each split hand are considered impossible.

From the estimated mass vector $\tilde{M}_\theta^s | \theta(s) = \text{Split}$, an approximation for the optimal policy $\theta^*$ is derived, referred to as the semi-optimal policy and denoted as $\tilde{\theta}^*$:

\begin{itemize}
    \item If $j_1 \in \{0, 2\}^{10}, |j_1| = 2, j_2 = \mathbf{0}$ and $\mathbb{E}[\tilde{Y}_\theta^s | \theta(s) = \text{Split}] = 2\mathbb{E}[Y_{\theta_2^*}^{(s')}] > \mathbb{E}[Y_{\theta_1^*}^s]$:
    \[
    \tilde{\theta}^*(s) = \text{Split}
    \]
    \item If $j_1 \in \{0, 2\}^{10}, |j_1| = 2, j_2 = \mathbf{0}$ and $\mathbb{E}[\tilde{Y}_\theta^s | \theta(s) = \text{Split}] = 2\mathbb{E}[Y_{\theta_2^*}^{(s')}] \leq \mathbb{E}[Y_{\theta_1^*}^s]$:
    \[
    \tilde{\theta}^*(s) = \theta_1^*(q, j_1, c, z)
    \]
    \item If $j_1 \notin \{0, 2\}^{10}, j_2 = \mathbf{0}$:
    \[
    \tilde{\theta}^*(s) = \theta_1^*(g, j_1, c, z)
    \]
    \item If $q = 1$, and $|j_2| > 0$:
    \[
    \tilde{\theta}^*(s) = \theta_2^*(g + j_2, j_1, c)
    \]
    \item If $q = 2$, and $|j_2| > 0$:
    \[
    \tilde{\theta}^*(s) = \theta_2^*(g + j_1, j_2, c)
    \]
\end{itemize}

Where $s = (q, j_1, j_2, c, z, u) \in S^{(J, \mathbf{d})}$ and $s' = \left( q + \frac{j_1}{2}, \frac{j_1}{2}, c \right) \in S_2^{(J, \mathbf{d})}$.

Finally, using the initial mass vectors ${\tilde{M}_{\tilde{\theta}^*}^s |(|j_1| = 2, j_2 = \mathbf{0})}$ along with the occurrence probabilities of the initial instances from the origin deck $\mathbf{d}$, the estimated mass vector for the round under the semi-optimal policy $\tilde{M}_{\tilde{\theta}^*}^\mathbf{d}$ and the basic strategy $\tilde{M}_{\theta^{Basic}}^\mathbf{d}$ is deduced:

\[
\tilde{M}_{\tilde{\theta}^*}^{\mathbf{d}} = \sum_{s \in S_{2}^{(J, \mathbf{d})}} Z_\theta(s, s_0) \cdot \tilde{M}_{\tilde{\theta}^*}^s, \quad 
\tilde{M}_{\theta^{Basic}}^{\mathbf{d}} = \sum_{s \in S_{2}^{(J, \mathbf{d})}} Z_\theta(s, s_0) \cdot \tilde{M}_{\theta^{Basic}}^s.
\]

\subsubsection{Basic Strategy}

The \textbf{Basic Strategy} has been studied for decades \cite{baldwin1956optimum, tejada1985estudio} and is generally presented in tabular form. For this work, we considered the basic strategy provided by the Wizard of Odds website \cite{shackleford_basic_strategy}, adjusted according to the rules established for this research. To compute the mass vector for the round under \textbf{basic strategy} $M_{(\theta^{Basic})}^\mathbf{d}$, the state return $Y_{\theta}^s | \theta^{Basic}(s) = \text{Split}$ was estimated in the same manner as for the semi-optimal round policy, resulting in an estimated mass vector $\tilde{M}_{(\theta^{Basic})}^s$.
\subsection{Betting Policy: Partial Policy}

In the theoretical framework, the optimality condition for the partial optimal policy was established for every count value \(c \in \text{Im}(\kappa)\) (see proof of Theorem 4):

\[
\omega_{(\kappa, \theta, u_1, H \to \infty)}^*(c) = \begin{cases} 
\displaystyle\arg \max_{b \in [0, 0.5]} \sum_{x \in R} \ln(1 + b \cdot x) \cdot \mathbb{E}[\mathbb{P}(X_{\theta}^{(\mathbf{D}|c)} = x)], & \text{if } \alpha = 0, \\
\displaystyle\arg \max_{b \in [0, 0.5]} \sum_{x \in R} (1 + b \cdot x)^\alpha \cdot \mathbb{E}[\mathbb{P}(X_{\theta}^{(\mathbf{D}|c)} = x)], & \text{if } \alpha > 0.
\end{cases}
\]

This section describes the methodology used to estimate two specific policies:

\begin{enumerate}
    \item $\omega_{(\kappa_1, \tilde{\theta}^*, u_1, H \to \infty)}^*$: 
    \newline Optimal Partial Policy for infinite rounds, under the True Count system $\kappa_1$, semi-optimal round policy $\tilde{\theta}^*$, and a CRRA utility function $u_1$ with constant relative risk aversion $1-\alpha$.

    \item $\omega_{(\kappa_1, \theta^{Basic}, u_1, H \to \infty)}^*$:
    \newline Optimal Partial Policy for infinite rounds, under the True Count system $\kappa_1$, “Basic Strategy” round policy $\theta^{Basic}$, and a CRRA utility function $u_1$ with constant relative risk aversion $1-\alpha$.
\end{enumerate}

\subsubsection{Variable Estimation}

Solving this problem involves calculating the expectation within the summations for every count \(c \in \text{Im}(\kappa)\) and every return \(x \in R\):

\[
\mathbb{E}[\mathbb{P}(X_{\theta}^{(\mathbf{D}|c)} = x)] = \sum_{\mathbf{d} \in \Omega_{(\tau,l)} \mid \kappa(\mathbf{d}) = c} \mathbb{P}(\mathbf{D}_{\theta} = \mathbf{d}) \cdot \mathbb{P}(X_{\theta}^\mathbf{d} = x)
\]

\paragraph{Variable Estimation}

\begin{enumerate}
    \item The calculation of \(\mathbb{P}(\mathbf{D}_{\theta} = \mathbf{d})\) depends on the evolution of the origin deck \(\mathcal{D}_{\theta}\). The complexity and size of the space to which this process belongs make it computationally unmanageable. To estimate \(\mathbb{P}(\mathbf{D}_{\theta} = \mathbf{d})\), we propose using the same variable \(\mathbf{D}_z \in \tilde{\Omega}_{(\tau,l)}\) for both round policies \(\tilde{\theta}^*\) and \(\theta^{Basic}\). The PMF of \(\mathbf{D}_z\) is defined from the evolution of the deck under the semi-optimal policy \(\mathcal{D}_{\tilde{\theta}^*}\), using a sampling of the first \(z\) origin decks:

\[
\mathbb{P}(\mathbf{D}_{\tilde{\theta}^*} = \mathbf{d}) \approx \mathbb{P}(\mathbf{D}_{\theta^{Basic}} = \mathbf{d}) \approx \mathbb{P}(\mathbf{D}_z = \mathbf{d}) = \frac{\sum_{n=0}^{z-1} 1_{\{\mathbf{d}^{obs}_n = \mathbf{d}\}}}{z }
\]

Where:
\begin{itemize}
    \item \(\mathcal{D}_{\tilde{\theta}^*}^{z} = \{\mathbf{d}^{obs}_n\}_{n=0}^{z-1}\): Partial sampling of the first \(z\) decks generated in the sequence \(\mathcal{D}_{\tilde{\theta}^*}\).
    \item \(\tilde{\Omega}_{(\tau,l)} \subseteq \Omega_{(\tau,l)}\): Set of decks generated in the sampling \(\mathcal{D}_{\tilde{\theta}^*}^{z}\).
\end{itemize}

    \item To estimate \(\mathbb{P}(X_{\theta}^\mathbf{d} = x)\), we use the approximate expression \(\tilde{\mathbb{P}}(X_{\theta}^\mathbf{d} = x)\) detailed in Section 3.1.
\end{enumerate}

\subsubsection{Estimated Partial Policy}

Based on the variable estimation, an approximation of the optimal partial policy is proposed using the following formulation, which we define as the estimated optimal partial policy:

\[
\omega_{(\kappa, \theta, u_1, H \to \infty)}^*(c) \approx\tilde{\omega}_{(\kappa, \theta, u_1, H \to \infty)}^*(c) = \begin{cases} 
\displaystyle\arg \max_{b \in [0, 0.5]} \sum_{x \in R} \ln(1 + b \cdot x) \cdot \mathbb{E}[\tilde{\mathbb{P}}(X_{\theta}^{(\mathbf{D}_z|c)} = x)], & \text{if } \alpha = 0, \\
\displaystyle\arg \max_{b \in [0, 0.5]} \sum_{x \in R} (1 + b \cdot x)^\alpha \cdot \mathbb{E}[\tilde{\mathbb{P}}(X_{\theta}^{(\mathbf{D}_z|c)} = x)], & \text{if } \alpha > 0.
\end{cases}
\]

\noindent Where:
\begin{enumerate}
    \item The expectation within the summation is calculated as follows:
    \[
    \mathbb{E}[\tilde{\mathbb{P}}(X_{\theta}^{(\mathbf{D}_z|c)} = x)]=\sum_{\mathbf{d} \in \tilde{\Omega}_{(\tau,l)} \mid \kappa(\mathbf{d}) = c} \mathbb{P}(\mathbf{D}_z = \mathbf{d}) \cdot \tilde{\mathbb{P}}(X_{\theta}^\mathbf{d} = x)
    \]
\end{enumerate}

\subsubsection{Computation of the Estimated Policy}

The steps performed to compute the two partial policies previously described are as follows:

\begin{enumerate}
    \item A single player was simulated against the dealer, following a semi-optimal round policy \(\tilde{\theta}^*\). This simulation considered an origin deck space \(\Omega_{(\tau=75\%, l=8)}\), in accordance with the established rules.
    \item From the simulation, a sample of the first 50,000 origin decks \(\mathcal{D}_{\tilde{\theta}^*}^{50.000}\) was obtained, thus determining \(\mathbb{P}(\mathbf{D}_{50.000} = \mathbf{d})\), \(\tilde{\Omega}_{(0.75,8)}\), and \(\text{Im}(\kappa_1) = \{\kappa(\mathbf{d}) \mid \mathbf{d} \in \tilde{\Omega}_{(0.75,8)}\}\).
    \item Using the algorithm designed for the round policy, the estimated probability \(\tilde{\mathbb{P}}(X_{\theta}^\mathbf{d} = x)\) was calculated for each return \(x \in R\) and each origin deck \(\mathbf{d} \in \tilde{\Omega}_{(0.75,8)}\). This calculation was performed under both the round policy \(\theta^{Basic}\) (Basic Strategy) and \(\tilde{\theta}^*\) (Semi-Optimal Policy).
    \item Finally, using these values, the estimated optimal partial policy was calculated for various values of \(\alpha\). This calculation was carried out using a function minimization algorithm implemented with the \texttt{minimize} function from the SciPy library.
\end{enumerate}

\subsection{Betting Policy: Complete Policy}

Given that a large proportion of game rounds have a negative expected return, a clarification is necessary regarding the policies addressed. While in the general model $H$ refers to the number of rounds in the session, we define $H^+$ as the number of rounds in which the expected return is positive during the session. This approach eliminates the variability induced by the proportion of rounds with a positive expected return.

This section describes the methodology used to estimate two specific complete policies:

\begin{enumerate}
    \item $\pi_{(\tilde{\theta}^*, u_2, H^+ = 100)}^*$:
    \newline Optimal Complete Policy for a total of $H^+ = 100$ rounds with positive expected returns, under the semi-optimal round policy $\tilde{\theta}^*$ and a CARA utility function $u_2$.

    \item $\pi_{(\theta^{Basic}, u_2, H^+ = 100)}^*$:
    \newline Optimal Complete Policy for a total of $H^+ = 100$ rounds with positive expected returns, under the "Basic Strategy" round policy $\theta^{Basic}$ and a CARA utility function $u_2$.
\end{enumerate}

To address the problem, an auxiliary surrogate problem is defined, based on a simplification of the evolution of the origin deck.

\subsubsection{Auxiliary Evolution of the Origin Deck}

The complexity and size of the space to which the evolution of the origin deck \(\mathcal{D}_{\theta}\) belongs make it computationally unmanageable. To address this limitation, an auxiliary evolution of the origin deck is proposed, in which the deck distribution in each round is independent of the others and matches the stationary distribution \(\mathbf{D} = \mathbf{D}_\theta\), conditioned on the round having a positive return:

\[
\mathcal{Y}_\theta = \{\mathbf{y}_n\}_{n=0}^{\infty}
\]

\noindent Where:
\begin{enumerate}
    \item \( \mathbf{y}_n \sim \mathbf{D}^+=\mathbf{D}_\theta^+ \in \Omega_{(\tau,l)}^+ \) represents the origin deck of the \(n\)-th round, with its probability mass function defined as follows:
    \[
    \mathbb{P}(\mathbf{y}_n= \mathbf{y}) = \mathbb{P}(\mathbf{D}^+_\theta = \mathbf{y})= 
    \dfrac{\mathbb{P}(\mathbf{D}_\theta = \mathbf{y})}{\displaystyle \sum_{\mathbf{d} \in \Omega_{(\tau,l)}^+} \mathbb{P}(\mathbf{D}_\theta = \mathbf{d})}
    \]

    \item The mass vector of the return in each round, \(M_\theta^{\mathbf{y}_n}\), is independently and identically distributed, with $M_\theta^{\mathbf{y}_n} \sim M_\theta^{\mathbf{D}^+} \in \Upsilon_{\theta}$. Its probability mass function is expressed as:
    \[
    \mathbb{P}(M_\theta^{\mathbf{y}_n} = \mathbf{m}) =\mathbb{P}(M_\theta^{\mathbf{D}^+} = \mathbf{m}) \sum_{\mathbf{y} \in \Omega_{(\tau, l)}^+ \mid M_\theta^{\mathbf{y}}= \mathbf{m}} \mathbb{P}(\mathbf{D}_\theta^+ = \mathbf{y})
    \]

\end{enumerate}

Where:
\begin{itemize}
    \item \( \Omega_{(\tau,l)}^+ =\{ \mathbf{d} \in \Omega_{(\tau,l)} \mid \mathbb{E}[X_\theta^\mathbf{d}] > 0 \} \): Space of favorable origin decks associated with rounds with positive expected returns.
    \item \( \Upsilon_{\theta} = \{M_{\theta}^{\mathbf{d}} \mid \mathbf{d} \in \Omega_{(\tau,l)}^+\} \): Space of favorable mass vectors associated with rounds with positive expected returns.
\end{itemize}
\subsubsection{Auxiliary Mass Vectors}

Under the auxiliary evolution \(\mathcal{Y}_\theta\), the origin deck of each round exclusively determines the return distribution for that round but does not influence the origin deck distribution in the next round. This implies that the state variable is not the origin deck of the round but rather the round's mass vector, allowing the analysis to be reduced to the evolution of the mass vector during the session.

However, the per-round mass vector \(M_\theta^{\mathbf{D}^+} \in \Upsilon_{\theta}\) is computationally prohibitive for use in dynamic programming due to the size of the space \(\Upsilon_{\theta}\). To address this limitation, the random variable \(M_\theta^{\mathbf{D}^+}\) is approximated with a new random variable, the \textbf{auxiliary mass vector}, denoted as \(M_{\theta, k} \in \Upsilon_{\theta, k}\). The probability mass function of \(M_{\theta, k}\) is defined as:

\[
\mathbb{P}(M_{\theta, k} = \mathbf{\hat{m}}_i) = \sum_{\mathbf{m} \in C_i} \mathbb{P}(M_\theta^{\mathbf{D}^+} = \mathbf{m})
\]

\noindent Where:
\begin{enumerate}
    \item \(\mathbf{\hat{m}}_i \in \Upsilon_{\theta, k}\) is an \textbf{auxiliary mass vector}, representing the centroid of the \(i\)-th cluster and approximating a subset of \(\Upsilon_{\theta}\).
    \item \(C_i\) is the \(i\)-th cluster, obtained by grouping the original mass vectors from \(\Upsilon_{\theta}\) using an intracluster variance minimization criterion.
    \item \(\mathbb{P}(M_\theta^{\mathbf{D}^+} = \mathbf{m})\) is the probability associated with the original mass vector \(\mathbf{m} \in \Upsilon_{\theta}\).
    \item \(\Upsilon_{\theta, k} = \{\mathbf{\hat{m}}_1, \mathbf{\hat{m}}_2, \dots, \mathbf{\hat{m}}_k\}\) is the reduced space containing the auxiliary mass vectors.
\end{enumerate}

This construction ensures that the auxiliary variable \(M_{\theta, k}\) provides a representative approximation of \(M_\theta^{\mathbf{D}^+}\), preserving its essential characteristics while significantly reducing computational complexity and costs, enabling its use in dynamic programming schemes.

\subsubsection{Auxiliary States}

To identify the state variables, we consider that the mass vector of each round is defined by the auxiliary mass vector \(M_{\theta, k}\). Based on the equation for wealth, it follows that the auxiliary state during a betting session, denoted as \(v\), is fully identified by the round’s mass vector, the current wealth, and the number of rounds played since the start. The auxiliary state space \(\mathcal{V}\) is defined as:

\[
\mathcal{V} = \{v = (\mathbf{\hat{m}}, w, n) \mid \mathbf{\hat{m}} \in \Upsilon_{\theta, k}, w \in \mathbb{R}_+, n \in \mathbb{Z}_+\}
\]

\subsubsection{Auxiliary Problem}

Analogous to the original problem, the optimization problem for the \textbf{auxiliary policy} $\lambda$ is defined as the maximization of utility at the end of a betting session. The optimization problem is formulated recursively as follows:

\[
\lambda^*(v_n=v)
\;=\;
\arg \max_{b \in [0, 0.5]}
\sum\limits_{v' \in \Gamma_{[\mathbf{\hat{m}}, b]}}
T_b(v', v)\,V^*(v'), \quad \forall v \in \mathcal{V} \mid {n<H}
\]

\noindent Where:
\begin{enumerate}
    \item \(V^*(v)\) represents the \textbf{optimal value} of a state \(v\), calculated as:
    \[
    V^*(v_{n} = v) \;=\;
    \begin{cases}
    \displaystyle\max_{b \in [0, 0.5]}
    \sum\limits_{v' \in \Gamma_{[\mathbf{\hat{m}}, b]}}
    T_b(v', v)\,V^*(v'),
    & \text{if } n<H
    \\[10pt]
    u_2(w; \beta), & \text{if } n=H
    \end{cases}
    \]
    
    \item The \textbf{auxiliary transition function} $ T_b(v', v) $ is expressed as:
    \[
    T_b(v', v) = \mathbb{P}(v_{n+1} = v' \mid v_n = v, \lambda(v) = b) = \begin{cases} 
    \mathbb{P}(M_{\theta, k} = \mathbf{\hat{m}'}) \cdot (\mathbf{\hat{m}})_i , & \text{if } v' \in \Gamma_{[\mathbf{\hat{m}}, b]} \\
    0, & \text{if } v' \notin \Gamma_{[\mathbf{\hat{m}}, b]}
    \end{cases}
    \]

    \item $\Gamma_{[\mathbf{\hat{m}}, b]}$ represents the set of states \( v' \in \mathcal{V} \) to which the player can transition from state \( v \) after placing a bet \( b \):
    \[
    \Gamma_{[\mathbf{\hat{m}}, b]} = \{(\mathbf{\hat{m}^*}, w^*, n+1) \mid \mathbf{\hat{m}^*} \in \Upsilon_{\theta, k}, w^* = w \cdot (1 + b \cdot r_i), r_i \in R \}
    \]
    
\end{enumerate}

Where:
\begin{itemize}
    \item \( v = (\mathbf{\hat{m}}, w, n) \in \mathcal{V} \): Current state.
    \item \( v' = (\mathbf{\hat{m}'}, w', n+1) \in \mathcal{V} \): Transitioned state.
    \item \(r_i = \frac{w' - w}{w \cdot b} \in R\): Return associated with the state transition.
\end{itemize}

\subsubsection{Variable Estimation}

Solving this problem requires calculating the optimal value for all auxiliary states \(v \in \mathcal{V} \mid n \leq H\). However, performing this calculation exactly is infeasible due to two main reasons: the continuous nature of wealth \(w\) and the inability to determine the transition function exactly. Therefore, approximate methods are necessary to address the problem.

\paragraph{Variable Estimation}

\begin{enumerate}
    \item Analogous to the case of the partial policy (see Section 3.2.1), a variable \(\mathbf{D}_z\) is constructed from a sampling of the first \(z\) origin decks \(\mathcal{D}_{\tilde{\theta}^*}^{z}\). Based on this variable, the following estimations are established:
    \begin{itemize}
        \item \(\mathbb{P}(\mathbf{D}_{\theta} = \mathbf{d}) \approx \mathbb{\tilde{P}}(\mathbf{D}_{\theta} = \mathbf{d}) = \mathbb{P}(\mathbf{D}_z = \mathbf{d})\): Estimation of the stationary deck PMF, identical for both round policies.
        
        \item \(\Omega_{(\tau,l)} \approx \tilde{\Omega}_{(\tau,l)}\): Estimation of the origin deck space based on the generated sample.
    \end{itemize}
    
    \item Based on the estimates \(\mathbb{\tilde{P}}(\mathbf{D}_{\theta} = \mathbf{d})\), \(\tilde{\Omega}_{(\tau,l)}\), and the estimated round return PMF (detailed in Section 3.1), the following approximations of the original variables are derived: \(\tilde{\Omega}_{(\tau,l)}^+\), \(\mathbb{\tilde{P}}(\mathbf{D}_{\theta}^+ = \mathbf{d})\), \(\tilde{\Upsilon}_{\theta}\), \(\mathbb{\tilde{P}}(M_\theta^{\mathbf{D}^+} = \mathbf{m})\), \(\tilde{\Upsilon}_{\theta, k}\), \(\mathbb{\tilde{P}}(M_{\theta, k} = \mathbf{\hat{m}})\), and \( \tilde{\Gamma}_{[\mathbf{\hat{m}}, b]}\). These estimates are constructed according to the definitions of each variable and enable approximating both the auxiliary state space and the auxiliary transition function of the session:
    
    \begin{itemize}
        \item \(\mathcal{\tilde{V}} = \{v = (\mathbf{\hat{m}}, w, n) \mid \mathbf{\hat{m}} \in \tilde{\Upsilon}_{\theta, k}, w \in \mathbb{R}_+, n \in \mathbb{Z}_+\}\): Estimate of the auxiliary state space.
        
        \item \(\tilde{T}_b(v', v) =\mathbb{\tilde{P}}(M_{\theta, k} = \mathbf{\hat{m}'}) \cdot (\mathbf{\hat{m}})_i, \quad \text{if } v' \in \tilde{\Gamma}_{[\mathbf{\hat{m}}, b]}\): Estimate of the auxiliary transition function.

    \end{itemize}

\end{enumerate}
\subsubsection{Estimated Auxiliary Policy}

The process for calculating the estimated optimal auxiliary policy for all auxiliary states \(v \in \mathcal{\tilde{V}}\) is detailed below:

\begin{enumerate}
    \item \textbf{Round \(n=H\)}: The optimal value is directly derived as \( V^*(v_{H} = v)= u_2(w; \beta)\), for all states \(v \in \mathcal{\tilde{V}}\) such that \(n=H\).
    \item \textbf{Round \(n=H-1\)}: The optimal value and the optimal policy are initially calculated only for the states \(v \in \mathcal{\tilde{V}}\) that satisfy the conditions \( w = A \cdot \epsilon\), where \( A \in \mathbb{Z}_+ \text{ and } 0 \leq A \leq 5 \cdot \epsilon^{-1}\), and \(n=H-1\).
    \item \textbf{Interpolation}: From the calculated state-optimal policy and state-optimal value pairs, interpolation is performed to estimate both the optimal policy and the optimal value for all states \(v \in \mathcal{\tilde{V}}\) such that \(n=H-1\).
    \item \textbf{Recursion}: By recursively repeating this process up to round \(n=0\), an estimation of the optimal auxiliary policy \(\tilde{\lambda}_{( \theta, u_2, H^+)}^*\) is obtained for all states \(v \in \mathcal{\tilde{V}}\).
\end{enumerate}

\subsubsection{Estimated Complete Policy}

The estimated optimal complete policy is defined in terms of the optimal auxiliary policy as follows:

\[
\pi_{( \theta, u_2, H^+)}^*(\psi) \approx \tilde{\pi}_{( \theta, u_2, H^+)}^*(\psi)
\]

\noindent Where:

\begin{enumerate}
    \item If \(\mathbb{\tilde{E}}[X_\theta^\mathbf{d}] \leq 0\), the policy prescribes no betting, i.e., \(\tilde{\pi}^*(\psi) = 0\).
    \item If \(\mathbb{\tilde{E}}[X_\theta^\mathbf{d}] > 0\), the policy \(\tilde{\pi}^*(\psi)\) is obtained through interpolation. For this, the pairs \((v, \tilde{\lambda}_{(\theta, u_2, H^+)}^*(v))\), with \(v \in \mathcal{\tilde{V}}\), are used. Considering the mass vector associated with the state \(\psi\), these points allow the interpolation and extension of the optimal auxiliary policy \(\tilde{\lambda}^*\) to any general state \(\psi \in \mathcal{B}\), thereby establishing the estimated optimal complete policy.
\end{enumerate}

\subsubsection{Computation of the Estimated Policy}

The steps taken to compute the two previously defined complete policies are as follows:

\begin{enumerate}
    \item The estimates \(\mathbb{\tilde{P}}(\mathbf{D}_{\theta} = \mathbf{d})\) and \(\tilde{\Omega}_{(0.75,8)}\) are obtained from a sample \(\mathcal{D}_{\tilde{\theta}^*}^{50.000}\) (see Section 3.2.3).
    \item Using these estimates, all corresponding estimated sets and expressions are generated. Specifically, the auxiliary vector \(M_{\theta, k=200}\) is constructed by applying the k-means algorithm, clustering the vectors of \(\tilde{\Upsilon}_{\theta}\) into 200 clusters.
    \item To optimize the estimated auxiliary policy, a conventional dynamic programming algorithm was implemented, structured in a vectorized manner for greater efficiency. At each level \(n\), the states \(v \in \mathcal{\tilde{V}}\) where \(w = A \cdot 0.0005\), with \(A \in \mathbb{Z}_+\) and \(0 \leq A \leq 10,000\), were considered. Subsequently, interpolation was performed over the calculated states, following the previously described procedure. Policy optimization was carried out for \(H=100\), three values of \(\beta\): \(0.15\), \(0.25\), and \(0.5\), and for both round policies.
\end{enumerate}

\subsection{Policy Evaluation}

To analyze the performance of different betting policies, two independent simulations were conducted, each consisting of 30,000 betting sessions. In each simulation, the player participated in all rounds of each session following a specific round policy, obtaining a return \(r_i \in R\) at the end of each round. In the first simulation, the player used the Basic Strategy, while in the second, a semi-optimal policy was employed. Before the start of each round, the estimated mass vector corresponding to that round was calculated based on the origin deck.

From the mass vector and the origin deck, several betting policies were evaluated simultaneously. In both simulations, the player maintained a fixed round policy throughout all sessions but applied multiple betting policies in parallel. Specifically, the estimated optimal partial policy and the estimated optimal complete policy were analyzed, considering various values of \(\alpha\) and \(\beta\) for the utility functions.

Each betting session was considered complete after 100 rounds were effectively wagered. At the end of each session, the final wealth achieved under each evaluated betting policy was recorded. This procedure was repeated for all 30,000 sessions in each simulation, resulting in histograms showing the distribution of final wealth after 100 wagered rounds for each betting policy.
\section{Results}

\subsection{Round Policies}

The set of origin decks \(\tilde{\Omega}_{(0.75, 8)}\) was obtained directly through the sampling \(\mathcal{D}_{\tilde{\theta}^*}^{50.000}\). For each deck in this set, the estimated mass vector \(\tilde{M}_\theta^\mathbf{d}\) was calculated under both the semi-optimal policy and the Basic Strategy. This process took approximately 0.5 seconds per deck, highlighting the efficiency of the proposed algorithm. As a result, histograms representing the distribution of estimated returns and estimated standard deviations associated with following a round policy across a series of game rounds were obtained.

The analysis of the histograms reveals that the semi-optimal policy provides a slight advantage over the Basic Strategy in terms of the expected return per round and the frequency of favorable rounds. However, specifically in favorable rounds, this advantage comes at the cost of greater return volatility. The characteristics of both strategies allow us to conclude that the additional complexity of the semi-optimal policy results in a relatively modest improvement over the Basic Strategy. This reaffirms the efficiency and practicality of the Basic Strategy, given its simplicity and strong performance.

It is worth noting that these results were obtained under a specific configuration: a full deck composed of 8 decks and a penetration of 75\%. In scenarios with origin decks containing fewer cards, the differences between the two policies are expected to be amplified, with the semi-optimal policy outperforming the Basic Strategy even further. However, in such cases, the technique used to approximate the possibility of performing a Split presents a higher margin of error, making the exploration of rounds associated with origin decks containing fewer cards a more challenging technical problem.

\subsubsection{Expected Return per Round}

The analysis of the histograms of the estimated expected return per round shows that the composition of the origin deck directly impacts how favorable the round is for participation. Under the semi-optimal policy, the average estimated expected return for all simulated rounds was -0.56\%, slightly surpassing the -0.79\% associated with the Basic Strategy. The proportion of rounds with a positive estimated expected return was 21.9\% under the semi-optimal policy, averaging an expected return of 0.81\%. Meanwhile, 17.7\% of the rounds were favorable under the Basic Strategy, with an average return of 0.64\%. The distribution of the estimated expected return is similar under both policies, with a symmetric and mesokurtic shape centered around \(\sim -0.65\%\), although the semi-optimal policy exhibits a slightly higher concentration of rounds with positive returns.

\begin{figure}[h]
    \centering
    \begin{subfigure}{0.48\textwidth}
        \centering
        \includegraphics[width=\textwidth]{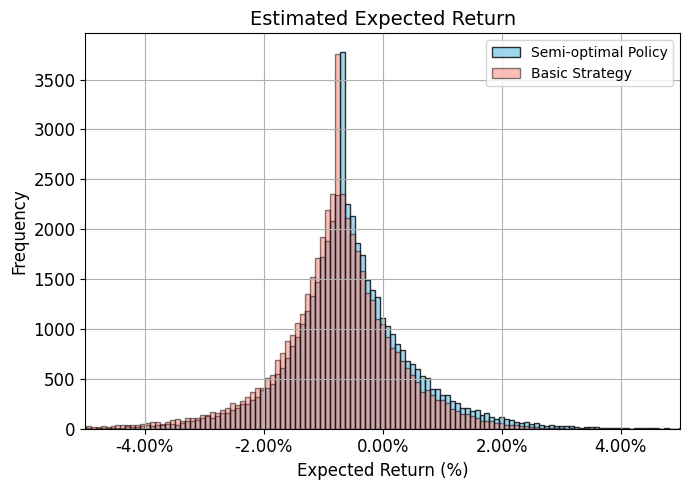}
        \caption{Histogram of the estimated expected return for the 50,000 sampled decks, under the semi-optimal policy and Basic Strategy.}
        \label{fig:Imagen1}
    \end{subfigure}
    \hfill
    \begin{subfigure}{0.48\textwidth}
        \centering
        \includegraphics[width=\textwidth]{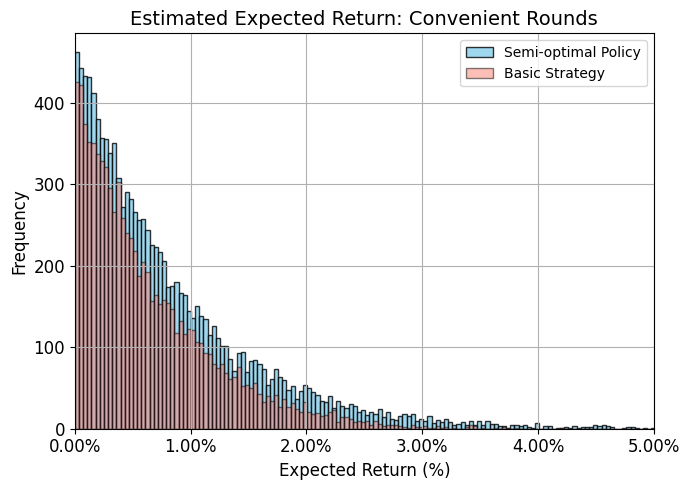}
        \caption{Histogram of the estimated expected return for favorable sampled decks, under the semi-optimal policy and Basic Strategy.}
        \label{fig:Imagen2}
    \end{subfigure}
    \caption{Comparison of the estimated expected return under the semi-optimal policy and Basic Strategy. The first subfigure analyzes the 50,000 sampled decks, while the second focuses on favorable decks.}
    \label{fig:conjunto_histogramas}
\end{figure}

\subsubsection{Standard Deviation per Round}

\begin{figure}[H]
    \centering
    \begin{subfigure}{0.48\textwidth}
        \centering
        \includegraphics[width=\textwidth]{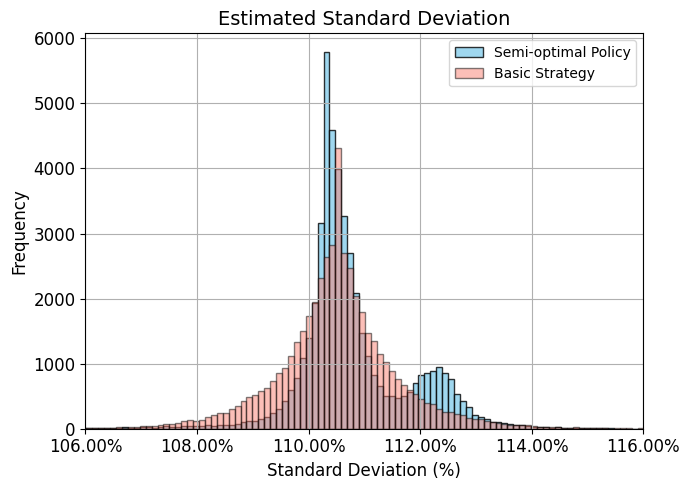}
        \caption{Histogram of the estimated standard deviation of the 50,000 sampled decks under the semi-optimal round policy and Basic Strategy}
        \label{fig:Imagen3}
    \end{subfigure}
    \hfill
    \begin{subfigure}{0.48\textwidth}
        \centering
        \includegraphics[width=\textwidth]{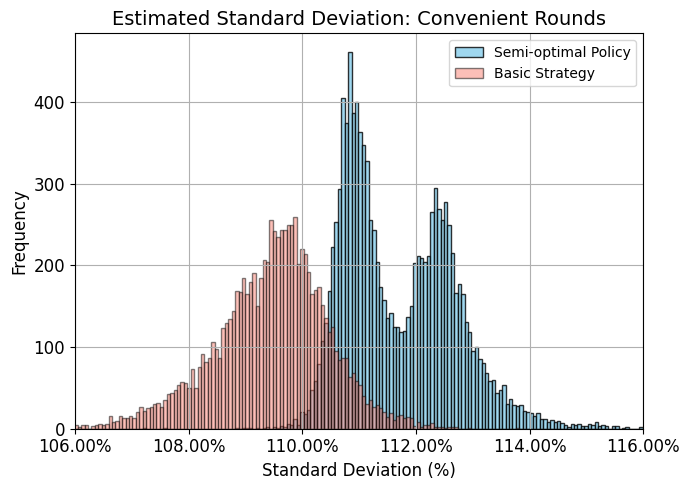}
        \caption{Histogram of the estimated standard deviation of the convenient sampled decks under the semi-optimal round policy and Basic Strategy}
        \label{fig:Imagen4}
    \end{subfigure}
    \caption{Analysis of the estimated standard deviation under the semi-optimal round policy and Basic Strategy. The first subfigure analyzes the 50,000 sampled decks, while the second focuses on the convenient decks.}
    \label{fig:conjunto_desviacion}
\end{figure}

The distribution of the estimated standard deviation of the return per round differs significantly under both round policies. As shown in Figure \ref{fig:Imagen3}, the Basic Strategy exhibits a highly symmetric histogram concentrated around $\sim 110.5\%$. In contrast, the histogram of the semi-optimal policy displays two concentration points: one more pronounced around $\sim 110.3\%$, and another, less intense, around $\sim 112.4\%$.

When analyzing the convenient decks, the differences between both policies become even more noticeable, as illustrated in Figure \ref{fig:Imagen4}. On the one hand, the histogram of the Basic Strategy presents a semi-symmetric shape with a mode around $\sim 109.5\%$ and a higher concentration of cases below this value. On the other hand, the histogram of the semi-optimal policy for these decks exhibits two similar peaks at $\sim 111\%$ and $\sim 112.5\%$, with a longer tail reflecting a larger number of high-deviation cases. Overall, the results indicate that the semi-optimal policy generally exhibits greater volatility than the Basic Strategy, with a particularly notable difference in the convenient rounds.

\subsection{Partial Policies}

From the 50,000 origin decks in the sample \(\tilde{\Omega}_{(0.75, 8)}\), the \textit{true count} \(\kappa_1(\mathbf{d})\) and the estimated mass vector \(\tilde{M}_\theta^\mathbf{d}\) were calculated for each deck under both the semi-optimal policy and the Basic Strategy. These calculations allowed the analysis of the distribution of the \textit{true count}, its correlation with the estimated expected return, and the properties of the optimal partial policy for different risk levels \(\alpha\).

The results show that the \textit{true count} follows an approximately geometric distribution and has an almost perfect linear relationship with the expected return of the round, which is higher under the semi-optimal policy compared to the Basic Strategy. Additionally, the optimal partial policy can be estimated with high precision using a linear expression based on the \textit{true count}, depending on the risk level \(\alpha\). Finally, the necessary parameters were calculated to estimate the lognormal distributions of the returns associated with each policy.

These findings suggest that a player could mentally implement optimal partial policies under CRRA for any value of \(\alpha\), while also being able to anticipate the distributions of returns associated with these strategies.

\subsubsection{True Count: Distribution}

\begin{figure}[h]
    \centering
    \includegraphics[width=0.8\textwidth]{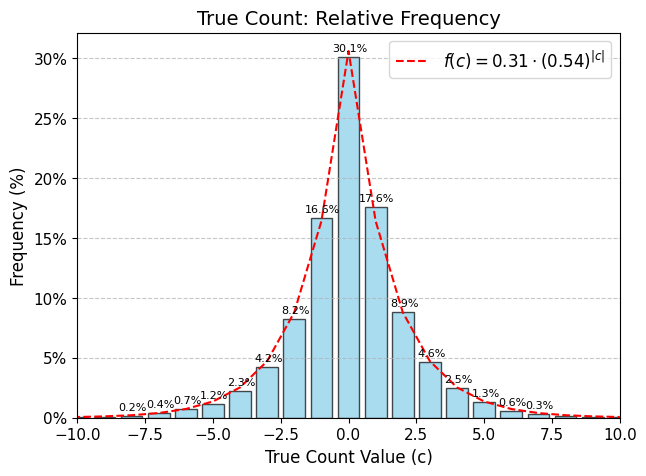}
    \caption{Histogram of relative frequencies of the \textit{true count} for the 50,000 sampled decks.}
    \label{fig:Imagen5}
\end{figure}

From the \textit{true counts} calculated for each of the 50,000 origin decks in the sample \(\tilde{\Omega}_{(0.75, 8)}\), the normalized frequency of each value was plotted. The resulting distribution has a symmetric shape centered at 0, with slight leptokurtosis, and values predominantly between -6 and 6, as illustrated in Figure \ref{fig:Imagen5}. To model this distribution, a geometric distribution was fitted, yielding a good approximation of the PMF of the stationary \textit{true count} value:
\[
\mathbb{\tilde{P}}\left(C_{(\theta, \kappa_1)}=c\right) \approx 0.31 \cdot 0.54^{\left|c\right|}
\]

\subsubsection{True Count and Expected Return}

\begin{figure}[H]
    \centering
    \begin{subfigure}{0.48\textwidth}
        \centering
        \includegraphics[width=\textwidth]{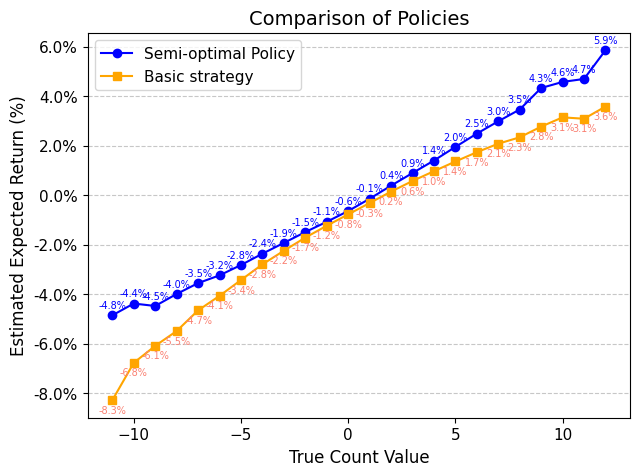}
        \caption{Relationship between the estimated expected return and the \textit{true count} under the semi-optimal policy and Basic Strategy.}
        \label{fig:Imagen6}
    \end{subfigure}
    \hfill
    \begin{subfigure}{0.48\textwidth}
        \centering
        \includegraphics[width=\textwidth]{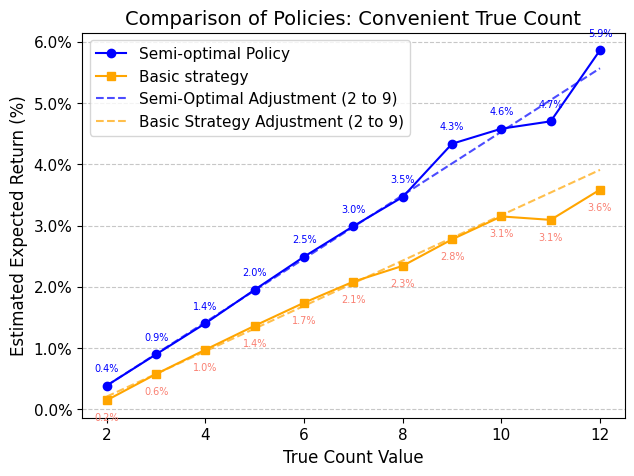}
        \caption{Relationship between the estimated expected return and the \textit{true count} for values greater than or equal to 2, under the semi-optimal policy and Basic Strategy.}
        \label{fig:Imagen7}
    \end{subfigure}
    \caption{Comparison between the estimated expected return and the \textit{true count} under semi-optimal and Basic Strategy round policies. The first subfigure shows the relationship across the entire \textit{true count} range, while the second focuses on values associated with a positive expected return.}
    \label{fig:conjunto_true_count}
\end{figure}

The estimated expected return of a round shows a significant correlation with the value of the \textit{true count} under both round policies, as illustrated in Figures \ref{fig:Imagen6} and \ref{fig:Imagen7}. While both policies perform similarly when the \textit{true count} values are near 0, marked differences are observed at the extremes. As the \textit{true count} becomes more negative, the expected return under the Basic Strategy decreases more rapidly compared to the semi-optimal policy. Conversely, when the \textit{true count} takes on more positive values, the expected return increases more rapidly under the semi-optimal policy than under the Basic Strategy. This reflects the ability of the semi-optimal policy to adjust to the deck's composition, in contrast to the Basic Strategy, which remains fixed and does not respond to changes in the deck.

When analyzing only the \textit{true count} values with positive returns, it is identified that betting becomes favorable for values greater than or equal to 2 under both policies, a condition that occurs in approximately 17\% of rounds. Ignoring outliers and focusing on \textit{true count} values between 2 and 9, a precise estimate of the expected return of the round is obtained depending on the applied round policy. In this range, the expected return shows an average increase of approximately 0.52\% per unit of \textit{true count} under the semi-optimal policy and 0.37\% under the Basic Strategy.

\paragraph{Return vs True Count}

\begin{enumerate}
    \item \textbf{Under Semi-Optimal Policy:}
    \[
    E[X_{\tilde{\theta}^*}^{(V|c)}] \approx 0.00517 \cdot c - 0.00646
    \]

    \item \textbf{Under Basic Strategy:}
    \[
    E[X_{\theta^{Basic}}^{(V|c)}] \approx 0.00369 \cdot c - 0.00529
    \]
\end{enumerate}

\subsubsection{Estimated Partial Policies}

The estimated optimal partial policy was calculated for multiple values of \(\alpha\) in the range from 0 to 1, considering both the semi-optimal policy and the Basic Strategy. Additionally, the parameters \(\mu\) and \(\sigma\) were determined to estimate the lognormal distribution of the returns associated with each policy. As expected, the results revealed that the estimated optimal partial policy is higher under the semi-optimal round policy and leads to more favorable return distributions compared to those obtained under the Basic Strategy.

For values of \(\alpha < 0.9\) and \textit{true count} values less than 9, an almost perfect correlation was observed between the partial policy and the \textit{true count}, allowing the partial policy to be linearly estimated as a function of the \textit{true count}, depending on the value of \(\alpha\) and the round policy employed. This linear relationship suggests that a player could mentally follow estimated partial policies during gameplay by using only the \textit{true count} value.

\begin{itemize}
    \item \textbf{Under Semi-Optimal Policy:}

Figure \ref{fig:Imagen8} shows a clear linear relationship between the estimated optimal partial policy and the \textit{true count} for different values of \(\alpha\), particularly when \(\alpha < 0.9\) and \textit{true count} is less than 9. This relationship can be expressed as follows:
\[
\alpha = \alpha_i \Longrightarrow \tilde{\omega}_{(\tilde{\theta}^\ast, u_1, H\rightarrow\infty)}^* \left(c\right) \approx m_i \cdot c + k_i
\]

The parameters \(m_i\) and \(k_i\) are specific to each value of \(\alpha\), as detailed in Table \ref{tab:parametros_alpha}. The same table also includes the parameters \(\mu_i\) and \(\sigma_i^2\), which estimate the lognormal distribution of wealth over time, reflecting how the value of \(\alpha\) influences the expectation and dispersion of returns.

\begin{figure}[H]
    \centering
    \begin{minipage}{0.48\textwidth}
        \centering
        \includegraphics[width=\textwidth]{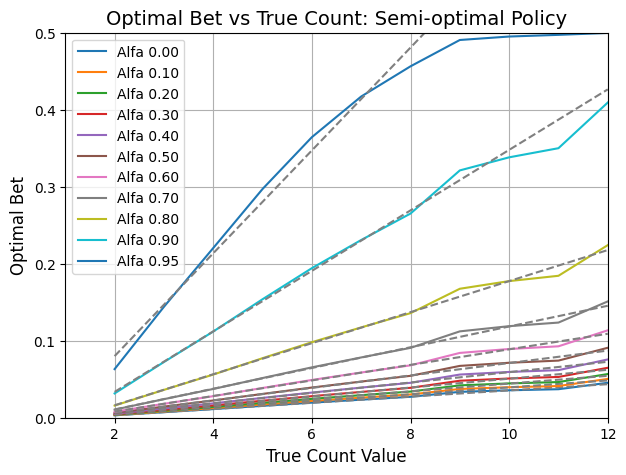}
        \caption{Relationship between the estimated optimal partial policy and the \textit{true count} for different values of \(\alpha\), under the semi-optimal policy.}
        \label{fig:Imagen8}
    \end{minipage}
    \hfill
    \begin{minipage}{0.48\textwidth}
        \centering
        \resizebox{\textwidth}{!}{%
        \begin{tabular}{|c|c|c|c|c|c|}
            \hline
            \(i\) & \(\alpha\) & \(m\) & \(k\) & \(\mu\) & \(\sigma^2\) \\
            \hline
            1  & 0     & 0.0040315  & -0.0048474  & 0.0000097   & 0.0000195   \\
            2  & 0.1   & 0.0044806  & -0.0053894  & 0.0000096   & 0.0000241   \\
            3  & 0.2   & 0.0050420  & -0.0060660  & 0.0000091   & 0.0000305   \\
            4  & 0.3   & 0.0057644  & -0.0069386  & 0.0000080   & 0.0000398   \\
            5  & 0.4   & 0.0067276  & -0.0080999  & 0.0000054   & 0.0000542   \\
            6  & 0.5   & 0.0080767  & -0.0097281  & -0.0000000  & 0.0000780   \\
            7  & 0.6   & 0.0101003  & -0.0121683  & -0.0000122  & 0.0001219   \\
            8  & 0.7   & 0.0134686  & -0.0162177  & -0.0000434  & 0.0002167   \\
            9  & 0.8   & 0.0201597  & -0.0241632  & -0.0001463  & 0.0004882   \\
            10 & 0.9   & 0.0392873  & -0.0450265  & -0.0007768  & 0.0019725   \\
            11 & 0.95  & 0.0667839  & -0.0535707  & -0.0033893  & 0.0085272   \\
            12 & 1     & 0          & 0.5         & -           & -           \\
            \hline
        \end{tabular}
        
        }
        \captionof{table}{Linear policy estimation and parameters of the associated return distribution for each \(\alpha\).}
        \label{tab:parametros_alpha}
    \end{minipage}
\end{figure}

    \item \textbf{Under Basic Strategy:}
    
Figure \ref{fig:Imagen9} shows a clear linear relationship between the estimated optimal partial policy and the \textit{true count} for different values of \(\alpha\), particularly when \(\alpha < 0.9\) and \textit{true count} is less than 9. This relationship can be expressed as follows:
\[
\alpha = \alpha_j \Longrightarrow \tilde{\omega}_{(\theta^{Basic}, u_1, H\rightarrow\infty)}^*\left(c\right) \approx m_j \cdot c + k_j
\]

The parameters \(m_j\) and \(k_j\) are specific to each value of \(\alpha\), as detailed in Table \ref{tab:parametros_alpha_basic_strategy}. The same table also includes the parameters \(\mu_j\) and \(\sigma_j^2\), which estimate the lognormal distribution of wealth over time, reflecting how the value of \(\alpha\) influences the expectation and dispersion of returns.

\begin{figure}[H]
    \centering
    \begin{minipage}{0.48\textwidth}
        \centering
        \includegraphics[width=\textwidth]{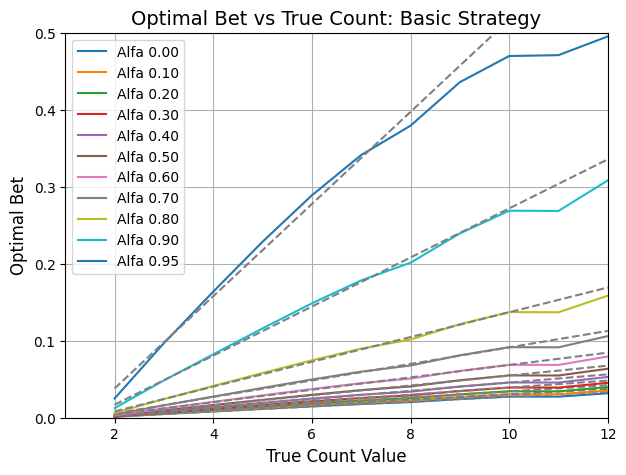}
        \caption{Relationship between the estimated optimal partial policy and the \textit{true count} for different values of \(\alpha\), under the Basic Strategy.}
        \label{fig:Imagen9}
    \end{minipage}
    \hfill
    \begin{minipage}{0.48\textwidth}
        \centering
        \resizebox{\textwidth}{!}{%
        \begin{tabular}{|c|c|c|c|c|c|}
            \hline
            \(j\) & \(\alpha\) & \(m\) & \(k\) & \(\mu\) & \(\sigma^2\) \\
            \hline
            1  & 0     & 0.0032125   & -0.0047781  & 0.0000047   & 0.0000094   \\
            2  & 0.1   & 0.0035702   & -0.0053116  & 0.0000047   & 0.0000116   \\
            3  & 0.2   & 0.0040170   & -0.0059768  & 0.0000044   & 0.0000147   \\
            4  & 0.3   & 0.0045919   & -0.0068328  & 0.0000038   & 0.0000192   \\
            5  & 0.4   & 0.0053587   & -0.0079755  & 0.0000026   & 0.0000262   \\
            6  & 0.5   & 0.0064326   & -0.0095757  & -0.0000000  & 0.0000377   \\
            7  & 0.6   & 0.0080439   & -0.0119765  & -0.0000059  & 0.0000590   \\
            8  & 0.7   & 0.0107292   & -0.0159750  & -0.0000210  & 0.0001048   \\
            9  & 0.8   & 0.0160884   & -0.0239306  & -0.0000708  & 0.0002360   \\
            10 & 0.9   & 0.0318662   & -0.0468364  & -0.0003771  & 0.0009485   \\
            11 & 0.95  & 0.0598706   & -0.0817353  & -0.0016756  & 0.0039103   \\
            12 & 1     & 0           & 0.5         & -           & -           \\
            \hline
        \end{tabular}

        }
        \captionof{table}{Linear policy estimation and parameters of the associated return distribution for each \(\alpha\).}
        \label{tab:parametros_alpha_basic_strategy}
    \end{minipage}
\end{figure}

\end{itemize}

\subsection{Estimated Complete Policies}

The estimated optimal complete policy was derived using the estimated optimal auxiliary policy. This auxiliary policy was calculated under both round policies and for three values of the risk aversion coefficient \(\beta\): 0.15, 0.25, and 0.35. As a result, six matrices of dimensions \(200 \times 10,001 \times 100\) were obtained, corresponding to each combination of \(\beta\) and round policy. The axes of these matrices are associated with the auxiliary state variables: auxiliary mass vector, wealth, and round number. The axis for auxiliary mass vectors was ordered ascendingly based on the expected return associated with each vector.

In general, the estimated optimal auxiliary policy exhibited higher values for lower \(\beta\) values and when considering auxiliary vectors derived from a semi-optimal policy. Additionally, consistent properties of this policy were identified for the three \(\beta\) values and both round policies. In particular, the policy was observed to be increasing with respect to the expected return associated with the auxiliary mass vector and decreasing with respect to wealth. Finally, the variability in betting as a function of the number of rounds was limited, suggesting low sensitivity of the optimal auxiliary policy to this variable.

\begin{figure}[H]
    \centering
    \begin{subfigure}[t]{0.32\textwidth}
        \centering
        \includegraphics[width=\textwidth]{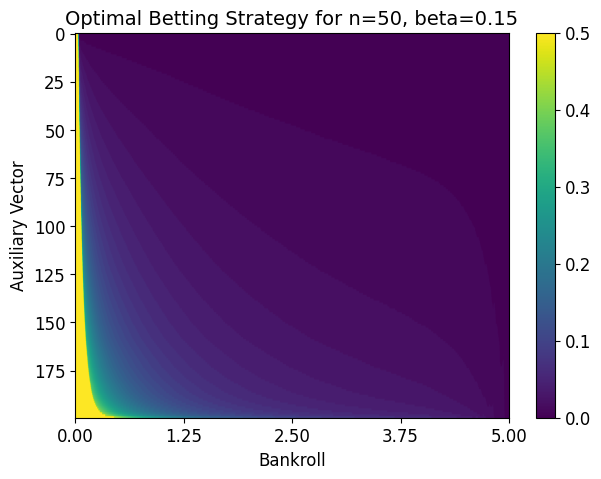}
        \caption{Heatmap of the estimated optimal auxiliary policy for round \(n = 50\) and a risk parameter \(\beta = 0.15\).}
        \label{fig:mapa_calor_patrimonio}
    \end{subfigure}
    \hfill
    \begin{subfigure}[t]{0.32\textwidth}
        \centering
        \includegraphics[width=\textwidth]{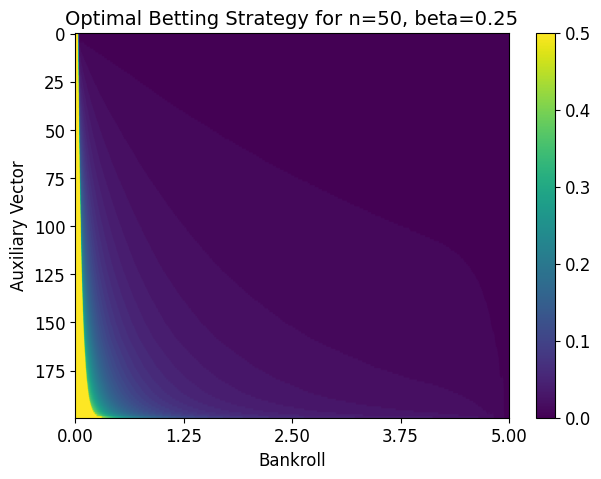}
        \caption{Heatmap of the estimated optimal auxiliary policy for round \(n = 50\) and a risk parameter \(\beta = 0.25\).}
        \label{fig:mapa_calor_vector}
    \end{subfigure}
    \hfill
    \begin{subfigure}[t]{0.32\textwidth}
        \centering
        \includegraphics[width=\textwidth]{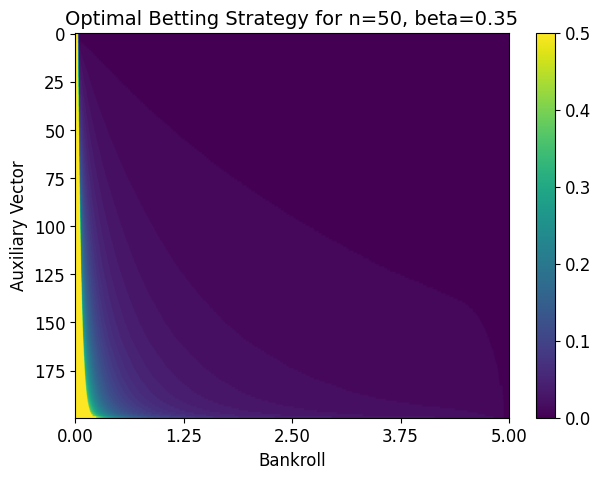}
        \caption{Heatmap of the estimated optimal auxiliary policy for round \(n = 50\) and a risk parameter \(\beta = 0.35\).}
        \label{fig:mapa_calor_patrimonio_ronda}
    \end{subfigure}
    \caption{Heatmaps of the estimated optimal auxiliary policy as a function of wealth and the index of the auxiliary mass vector, for \(n = 50\) rounds and under three risk levels \(\beta\): 0.15, 0.25, and 0.35.}
    \label{fig:mapas_calor_politica}
\end{figure}

\subsection{Evaluation: Return Distribution}

The estimated partial and complete policies were evaluated by simulating a realistic gaming environment, with the player participating only against the dealer. The partial policy was tested for three values of \(\alpha\) (0.3, 0.6, and 0.9), following both round policies. Similarly, the complete policy was evaluated for three values of \(\beta\) (0.15, 0.25, and 0.35), also under the two round policies. In total, 12 different policies were tested, each subjected to 30,000 betting sessions. Each session ended when the policy placed bets in 100 rounds.

The results showed that the return distribution under the estimated optimal partial policy fits well to a lognormal distribution, in line with theoretical expectations. In contrast, although the return distributions under the estimated optimal complete policy exhibited similarities to a normal distribution, this similarity is expected to diminish with a higher number of rounds or lower risk aversion levels. Both partial and complete policies using the semi-optimal round policy achieved higher expected returns and better Sharpe ratios compared to those following the Basic Strategy. Overall, the findings indicate that it is feasible to implement desirable betting policies in practice, capable of reasonably balancing expected return and return dispersion during a betting session.

The formalization of the betting policy's optimality criterion through a utility function proved to be an adequate methodology, allowing precise adjustment of the risk profile of the policies. Future studies could explore and compare return distributions under different rule configurations and utility functions, evaluating their influence on optimal return distributions. In particular, it would be interesting to analyze how key game variables, such as the number of decks in the complete deck and the penetration point, influence the results.

\subsubsection{Return Distribution: Partial Policies}

\begin{figure}[H]
    \centering
    \begin{minipage}{0.48\textwidth}
        \centering
        \includegraphics[width=\textwidth]{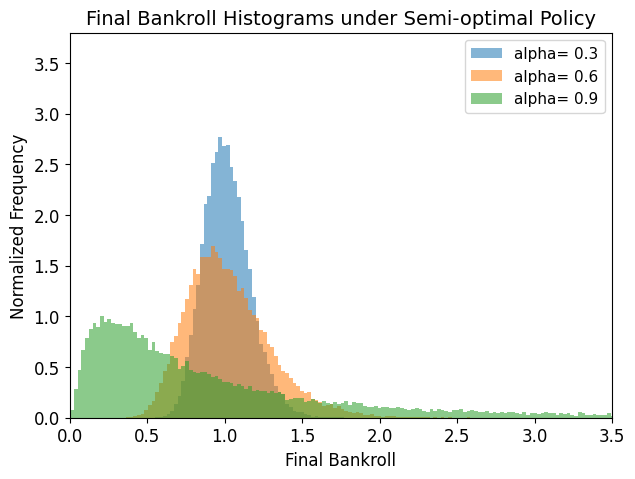}
        \caption{Distributions of final wealth under the estimated optimal partial policy and semi-optimal round policy for \(\alpha = 0.3\), \(\alpha = 0.6\), and \(\alpha = 0.9\).}
        \label{fig:patrimonio_semi_optima}
    \end{minipage}
    \hfill
    \begin{minipage}{0.48\textwidth}
        \centering
        \includegraphics[width=\textwidth]{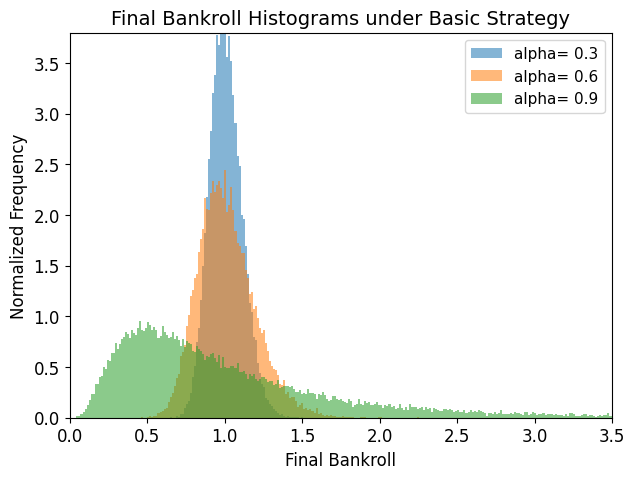}
        \caption{Distributions of final wealth under the estimated optimal partial policy and Basic Strategy for \(\alpha = 0.3\), \(\alpha = 0.6\), and \(\alpha = 0.9\).}
        \label{fig:patrimonio_basic_strategy}
    \end{minipage}

    \begin{minipage}{0.48\textwidth}
        \centering
        \resizebox{\textwidth}{!}{%
        \begin{tabular}{|c|c|c|c|c|c|}
            \hline
            \(\alpha\) & ROR & Median & Expected Return & Standard Deviation & Sharpe Ratio \\
            \hline
            0.3 & 0.00\%  & 0.35\%   & 1.51\%   & 15.11\%   & 10.01\%  \\
            0.6 & 0.00\%  & -0.86\%  & 2.64\%   & 27.09\%   & 9.74\%   \\
            0.9 & 0.00\%  & -32.79\% & 12.49\%  & 153.79\%  & 8.12\%   \\
            \hline
        \end{tabular}
        }
        \captionof{table}{Performance indicators under the estimated optimal partial policy and semi-optimal round policy.}
        \label{tab:resultados_semi_optima}
    \end{minipage}
    \hfill
    \begin{minipage}{0.48\textwidth}
        \centering
        \resizebox{\textwidth}{!}{%
        \begin{tabular}{|c|c|c|c|c|c|}
            \hline
            \(\alpha\) & ROR & Median & Expected Return & Standard Deviation & Sharpe Ratio \\
            \hline
            0.3 & 0.00\%  & 0.21\%   & 0.74\%   & 10.32\%  & 7.15\%  \\
            0.6 & 0.00\%  & -0.37\%  & 1.30\%   & 18.34\%  & 7.08\%  \\
            0.9 & 0.00\%  & -18.53\% & 5.24\%   & 86.33\%  & 6.07\%  \\
            \hline
        \end{tabular}
        }
        \captionof{table}{Performance indicators under the estimated optimal partial policy and Basic Strategy.}
        \label{tab:resultados_basic_strategy}
    \end{minipage}
\end{figure}

The return distribution under the estimated optimal partial policy aligns with a lognormal distribution, as shown in Figures \ref{fig:patrimonio_semi_optima} and \ref{fig:patrimonio_basic_strategy}. Tables \ref{tab:resultados_semi_optima} and \ref{tab:resultados_basic_strategy} indicate that both the volatility and the expected return of the returns increase as the value of \(\alpha\) rises, reflecting a lower relative risk aversion. Additionally, higher expected returns and better Sharpe ratios are observed when following the semi-optimal round policy compared to the Basic Strategy, highlighting the advantage of following this policy over multiple game rounds.

\subsubsection{Return Distribution: Complete Policies}

\begin{figure}[H]
    \centering
    \begin{minipage}{0.48\textwidth}
        \centering
        \includegraphics[width=\textwidth]{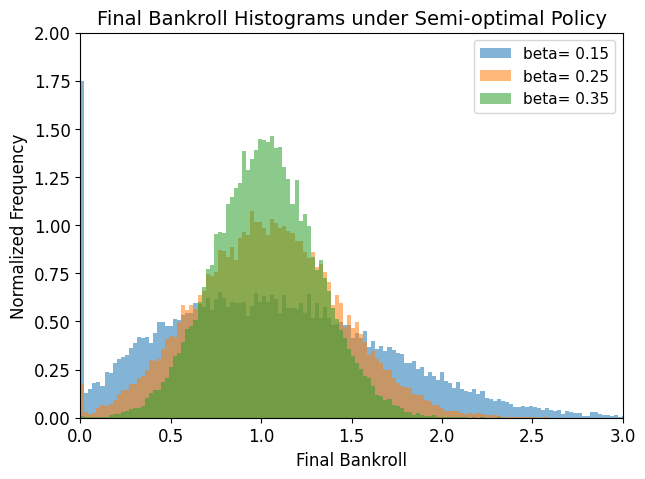}
        \caption{Distributions of final wealth under the estimated optimal complete policy and semi-optimal round policy for \(\beta = 0.15\), \(\beta = 0.25\), and \(\beta = 0.35\).}
        \label{fig:patrimonio_beta_semi_optima}
    \end{minipage}
    \hfill
    \begin{minipage}{0.48\textwidth}
        \centering
        \includegraphics[width=\textwidth]{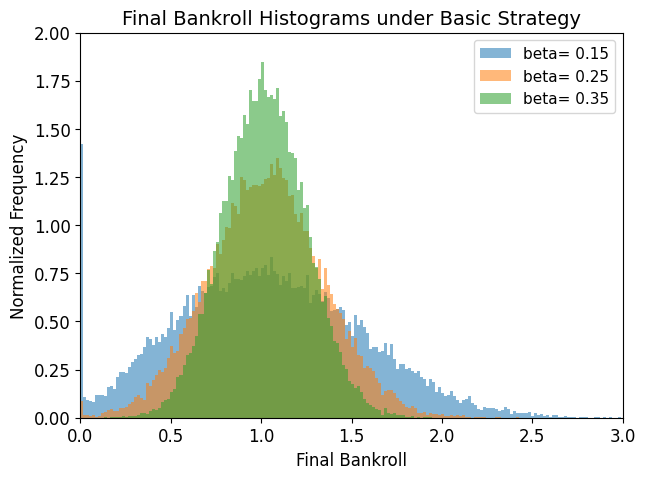}
        \caption{Distributions of final wealth under the estimated optimal complete policy and Basic Strategy for \(\beta = 0.15\), \(\beta = 0.25\), and \(\beta = 0.35\).}
        \label{fig:patrimonio_beta_basic_strategy}
    \end{minipage}

    \begin{minipage}{0.48\textwidth}
        \centering
        \resizebox{\textwidth}{!}{%
        \begin{tabular}{|c|c|c|c|c|c|}
            \hline
            \(\beta\) & ROR & Median & Expected Return & Standard Deviation & Sharpe Ratio \\
            \hline
            0.15 & 3.22\% & 2.85\% & 6.72\%  & 61.34\% & 10.95\% \\
            0.25 & 0.32\% & 3.23\% & 3.68\%  & 39.29\% &  9.38\% \\
            0.35 & 0.01\% & 2.38\% & 2.64\%  & 28.47\% &  9.28\% \\
            \hline
        \end{tabular}
        }
        \captionof{table}{Performance indicators under the estimated optimal complete policy and semi-optimal round policy.}
        \label{tab:resultados_beta_semi_optima}
    \end{minipage}
    \hfill
    \begin{minipage}{0.48\textwidth}
        \centering
        \resizebox{\textwidth}{!}{%
        \begin{tabular}{|c|c|c|c|c|c|}
            \hline
            \(\beta\) & ROR & Median & Expected Return & Standard Deviation & Sharpe Ratio \\
            \hline
            0.15 & 1.96\% & 2.36\% & 3.95\%  & 51.50\% & 7.66\% \\
            0.25 & 0.11\% & 2.17\% & 2.47\%  & 32.18\% & 7.68\% \\
            0.35 & 0.00\% & 1.51\% & 1.69\%  & 22.91\% & 7.36\% \\
            \hline
        \end{tabular}
        }
        \captionof{table}{Performance indicators under the estimated optimal complete policy and Basic Strategy.}
        \label{tab:resultados_beta_basic_strategy}
    \end{minipage}
\end{figure}

The return distribution under the estimated optimal complete policy exhibits high symmetry and a bell-shaped curve similar to a normal distribution, as shown in Figures \ref{fig:patrimonio_beta_semi_optima} and \ref{fig:patrimonio_beta_basic_strategy}. However, this similarity is expected to diminish for lower risk aversion levels or longer betting sessions, partly due to the irregularity induced by the possibility of ruin. Tables \ref{tab:resultados_beta_semi_optima} and \ref{tab:resultados_beta_basic_strategy} highlight that the semi-optimal round policy offers higher expected returns, a higher median, and a superior Sharpe ratio compared to the Basic Strategy, demonstrating the advantage of following this policy over multiple game rounds.

\section{Conclusion}
This study aims to formalize and solve the optimization challenges associated with two critical decision components in Blackjack: the round policy, which governs in-game actions, and the betting policy, which determines how much to bet. Using a combination of Markov Decision Processes (MDPs) and Expected Utility Theory, we developed a framework capable of generating strategies that approach theoretical optimality while remaining computationally efficient.

On the technical front, round policy optimization demonstrated how dynamic programming techniques can compute strategies that outperform traditional heuristics such as Basic Strategy. On the betting side, the integration of constant relative risk aversion (CRRA) utility functions not only improved the adaptability of the betting policy to different player profiles, but also extended the results beyond blackjack, providing insights relevant to investment and risk management scenarios.

Perhaps the most compelling result of this study is the confirmation that Blackjack is effectively solved under the defined rules and assumptions. The optimized strategies presented here achieve return distributions that leave little room for practical improvement, providing a clear roadmap for players seeking to maximize their results. Beyond the confines of this specific game, however, the methods developed have broader implications. The principles of decision optimization, utility maximization, and computational efficiency explored in this study have the potential to impact other domains, from finance to resource allocation to dynamic system control.

Despite these significant achievements, the research also highlights the complexities of Blackjack optimization and points to areas for further research. One notable challenge is reducing the residual error introduced by approximations, such as those applied to the "split" action. As variations of Including additional player interactions and alternatives introduce new levels of complexity, future work could extend this framework to account for such dynamics. In addition, the inclusion of multiplayer scenarios or simulations involving multiple hands per player could provide insights into the role of strategic diversity in decision-making.

\newpage
\begin{appendices}
\section{Appendix}

\subsection{Round Policy Algorithm}
This section provides a detailed explanation of the algorithm used to compute the semi-optimal round policy and the PMF of the return per round.

\subsubsection{Player Sequences}
\label{sec:secuencias_jugador}
Using a code (Appendix~\ref{sec:anexo2}) that deals cards sequentially to a player, from an empty sequence until exceeding 21 points by at most one card, the DAG \( G^J = (U^J, A^J) \) is traversed, computing the following data structures:

\begin{itemize}
    \item \(\mathbf{U}^J\): A matrix of size \(\mathbb{R}^{|U^J| \times 10}\) representing the set \( U^J \). Each row \(\mathbf{U}^J[i]\) corresponds to a player sequence \(\mathbf{j} \in U^J\), where \(\mathbf{U}^J[0]\) is a sequence of all zeros.
    
    \item \(\mathbf{A}^J\): A dictionary representing the directed edges \( A^J \). \(\mathbf{A}^J[i]\) returns a list of 10 indices associated with rows of \(\mathbf{U}^J\), representing the 10 "child" sequences to which the sequence \(\mathbf{U}^J[i]\) can transition when a new card is dealt to the player.
\end{itemize}

By adding the limitation to the algorithm that cards are only added to the previous sequence if the total does not exceed 21 points, the DAG \( G^{(J,r)} = (U^{(J,r)}, A^{(J,r)}) \) is traversed, and the following data structures are computed and stored:

\begin{itemize}
    \item \(\mathbf{U}^{(J,r)}\): A matrix of size \(\mathbb{R}^{|U^{(J,r)}| \times 10}\), representing the set \( U^{(J,r)} \). Each row \(\mathbf{U}^{(J,r)}[i]\) corresponds to a relevant player sequence \(\mathbf{j} \in U^{(J,r)}\), where \(\mathbf{U}^{(J,r)}[0]\) is a sequence of all zeros.
    
    \item \(\mathbf{A}^{(J,r)}\): A dictionary representing the directed edges \( A^{(J,r)} \). \(\mathbf{A}^{(J,r)}[i]\) returns a list of indices associated with rows of \(\mathbf{U}^{(J,r)}\), representing the "child" sequences to which the sequence \(\mathbf{U}^{(J,r)}[i]\) can transition.

    \item \(\mathbf{U}_n^{(J,r)}\): A dictionary representing the sequences \( U_n^{(J,r)} \). \(\mathbf{U}_n^{(J,r)}[n]\) returns a matrix of size \(\mathbb{R}^{|U_n^{(J,r)}| \times 10}\), where each row is a sequence \( \mathbf{j} \in U_n^{(J,r)} \).
\end{itemize}

Finally, based on the sequences \(\mathbf{U}^{(J,r)}\), an algorithm generates a directed spanning tree (AGD) of the DAG \( G^{(J,r)} \). To achieve this, all sequences in the DAG must be traversed, ensuring that each sequence has a single parent, except for the root. As a result, the following data structure is computed:

\begin{itemize}
    \item \(\mathbf{A}_T^{(J,r)}\): A dictionary representing the directed edges \( A_T^{(J,r)} \). \(\mathbf{A}_T^{(J,r)}[i]\) returns a list of indices associated with rows of \( \mathbf{A}^{(J,r)} \), representing the "child" sequences to which the sequence \( \mathbf{A}^{(J,r)}[i] \) can transition, ensuring that all sequences have a single parent.
\end{itemize}

\subsubsection{Dealer Sequences}
\label{sec:secuencias_crupier}
Using a code (Appendix~\ref{sec:anexo2}) that deals cards sequentially to a dealer, from an empty sequence until exceeding 16 points by at most one card, the DAG \( G^C = (U^C, A^C) \) is traversed, computing the following data structures:

\begin{itemize}
    \item \(\mathbf{U}^C\): A matrix of size \(\mathbb{R}^{|U^C| \times 10}\) representing the set \( U^C \). Each row \(\mathbf{U}^C[i]\) corresponds to a dealer sequence \(\mathbf{c} \in U^C\).
    
    \item \(\mathbf{A}^C\): A dictionary representing the directed edges \( A^C \). \(\mathbf{A}^C[i]\) returns a list of 10 indices associated with rows of \(\mathbf{U}^C\), representing the 10 "child" sequences to which the sequence \(\mathbf{U}^C[i]\) can transition when a new card is dealt to the dealer.
\end{itemize}

From the DAG \( G^C \) where the total does not exceed 21 points, the following additional data structures are derived:

\begin{itemize}
    \item \(\mathbf{U}^{(C,r)}\): A matrix of size \(\mathbb{R}^{|U^{(C,r)}| \times 10}\) representing the set \( U^{(C,r)} \). Each row \(\mathbf{U}^{(C,r)}[i]\) corresponds to a relevant dealer sequence \(\mathbf{c} \in U^{(C,r)}\).
    
    \item \(\mathbf{Mult\_up}\): A dictionary where \(\mathbf{Mult\_up[\mathbf{c^*}]}\) returns a two-column matrix for each upcard \( \mathbf{c^*} \):
    \begin{itemize}
        \item The first column contains the indices of the \(\mathbf{U}^{(C,r)}\) matrix corresponding to terminal sequences reachable from the upcard \( \mathbf{c^*} \).
        \item The second column indicates the multiplicities of these terminal sequences from \( \mathbf{c^*} \), calculated as \( m(\mathbf{c}, \mathbf{c}^*) \).
    \end{itemize}
\end{itemize}
\subsubsection{Final Return Matrix Calculation}
\label{sec:retornos_finales}

Using the matrices \(\mathbf{U}^{(J,r)}\) and \(\mathbf{U}^{(C,r)}\), the scores associated with the relevant sequences of the player and dealer are computed using the code described in Appendix~\ref{sec:anexo1}. These scores are then compared between the sequences \(\mathbf{j_1} \in U^{(J,r)}\) and \(\mathbf{c} \in U^{(C,r)}\), yielding the final returns \( f_\mathcal{R}(\mathbf{g}, \mathbf{j_1}, \mathbf{c}, z=0) \). As a result, a matrix \( \mathbf{R} \in \mathbb{R}^{|U^{(J,r)}| \times |U^{(C,r)}|} \) is generated, such that:

\[
\mathbf{R}[i, k] =
\begin{cases} 
    1.5, & \text{if } \mathbf{U}^{(J,r)}[i] \text{ against } \mathbf{U}^{(C,r)}[k] \text{ is a Blackjack}, \\
    1, & \text{if } \mathbf{U}^{(J,r)}[i] \text{ beats } \mathbf{U}^{(C,r)}[k], \\
    0, & \text{if } \mathbf{U}^{(J,r)}[i] \text{ ties with } \mathbf{U}^{(C,r)}[k], \\
    -1, & \text{if } \mathbf{U}^{(C,r)}[k] \text{ beats } \mathbf{U}^{(J,r)}[i].
\end{cases}
\]

\subsubsection{Extraction Probability Calculation}
\label{sec:probabilidades_extraccion}

The calculation of extraction probabilities is implemented using the \textbf{Probs\_secs()} function (Appendix~\ref{sec:anexo3}) and the matrix \(\mathbf{U}^{(C,r)}\). These tools allow for the vectorized calculation of the probabilities of extracting sequences \(\mathbf{c} \in U^{(C,r)}\) from an adjusted initial deck \(\mathbf{d}\), where:
\(
\mathbf{g} = \mathbf{d} - \mathbf{U}^{(J,r)}_0.
\). As a result, a probability vector \(\mathbf{Q} \in \mathbb{R}^{|U^{(C,r)}|}\) is generated, defined using Einstein notation as:

\[
\mathbf{Q}_k = \zeta(\mathbf{U}^{(C,r)}_k, \mathbf{g}),
\]

where \(\zeta(\mathbf{U}^{(C,r)}_k, \mathbf{d})\) represents the probability of extracting the sequence \(\mathbf{U}^{(C,r)}_k\) from the adjusted initial deck \(\mathbf{d}\).

With this initial vector \(\mathbf{Q}\) and the transition dictionary \(\mathbf{A}_T^{(J,r)}\), the Directed Spanning Tree (AGD) \(G_T^{(J,r)} = (U^{(J,r)}, A_T^{(J,r)})\) is iteratively traversed. This traversal calculates the probabilities of extracting the sequences \(\mathbf{c} \in U^{(C,r)}\) from adjusted decks of the form \(\mathbf{d} = \mathbf{d} - \mathbf{j}_1\), where \(\mathbf{j}_1 \in U^{(J,r)}\). The result of this calculation is a matrix \(\mathbf{P} \in \mathbb{R}^{|U^{(J,r)}| \times |U^{(C,r)}|}\), where each entry is defined using Einstein notation as:

\[
\mathbf{P}_{i k} = \zeta(\mathbf{U}^{(C,r)}_k, \mathbf{d} - \mathbf{U}^{(J,r)}_i),
\]

where \(\mathbf{P}_{i k}\) represents the probability of extracting the sequence \(\mathbf{U}^{(C,r)}_k\) from the adjusted deck \(\mathbf{d} - \mathbf{U}^{(J,r)}_i\).

\paragraph{Iterative Process}

\textbf{Initialization:}  
At the root node of the AGD, associated with the initial deck \(\mathbf{g} = \mathbf{d}\), the initial probabilities are defined as:  
\[
\mathbf{P}_{0,k} = \mathbf{Q}_k = \zeta(\mathbf{U}^{(C,r)}_{k}, \mathbf{d}),
\]
where:
\begin{itemize}
    \item \(\mathbf{P}_{0,k}\) represents the initial probabilities at the root node, corresponding to the relevant sequences \( \mathbf{U}^{(C,r)}_k \).
    \item \(\mathbf{Q}_k\) is the initial probability vector calculated from the function \(\zeta\), which depends on the relevant sequences \( \mathbf{U}^{(C,r)}_{k} \) and the initial deck \(\mathbf{d}\).
\end{itemize}
This calculation is performed simultaneously for all indices \(k\) of the relevant sequences \( \mathbf{U}^{(C,r)} \).

\textbf{Tree Traversal:}  
For each node \(j_1 \in G_T^{(J,r)}\), associated with an adjusted deck \(\mathbf{g} = \mathbf{d} - j_1\), the extraction probabilities are calculated for each child \(j_1'\) of the node \(j_1\), where \(j_1'\) corresponds to a sequence obtained by adding a single card to \(j_1\), represented as a canonical sequence \(e_x \in \mathbb{R}^{10}\). The adjusted deck of the child is \(\mathbf{g}' = \mathbf{g} - e_x\). The probabilities are updated simultaneously for all relevant sequences \( \mathbf{U}^{(C,r)} \), corresponding to the dealer's sequences, using Einstein notation as follows:

\[
\mathbf{P}_{i,k} = \mathbf{P}_{p,k} \cdot \left( \frac{(\mathbf{g}_{x} - \mathbf{U}^{(C,r)}_{kx})}{\sum_{y} (\mathbf{g}_{y} - \mathbf{U}^{(C,r)}_{ky})} \cdot \frac{\sum_{y} \mathbf{g}_{y}}{\mathbf{g}_{x}} \right),
\]

where:
\begin{itemize}
    \item The calculation of \(\mathbf{P}_{i,k}\) is performed simultaneously for all indices \(k\) of the relevant sequences \( \mathbf{U}^{(C,r)} \).
    \item Formally, this implies that the values of \(\mathbf{P}_{i,k}\) are defined for each \(k\) as the result of the vector operation over all sequences \(k \in \{1, \dots, |U^{(C,r)}|\}\).
    \item \(\mathbf{P}_{i,k}\) is the probability associated with the node \(j_1'\) for the dealer sequence \(k\).
    \item \(\mathbf{P}_{p,k}\) is the probability associated with the parent node \(j_1\) for the dealer sequence \(k\).
\end{itemize}

\textbf{Matrix \(\mathbf{P}\) Update:}  
The calculation of \(\mathbf{P}_{i,k}\) for all indices \(k\) simultaneously updates the rows of the probability matrix \(\mathbf{P}\), defining the value of \(\mathbf{P}_i\) as a vector whose dimension corresponds to all relevant sequences of \( \mathbf{U}^{(C,r)} \).

\textbf{Complete Iteration:}  
The iterative process continues until all nodes of the tree \(G_T^{(J,r)}\) have been traversed, completing the calculation of \(\mathbf{P}\).

\paragraph{Results}  
This procedure ensures that each entry \(\mathbf{P}[i, k]\) correctly captures the probability of extracting \(\mathbf{U}^{(C,r)}[k]\) from any adjusted round deck \(\mathbf{g}=\mathbf{d} - \mathbf{U}^{(J,r)}[i]\), thus achieving an efficient and structured calculation. For additional details on the implementation, see Appendix~\ref{sec:anexo4}.

\subsubsection{Calculation of Extraction Probabilities for Player Sequences by Upcard}
\label{sec:calculo_secuenciasj_upcards}

Additionally, the extraction probability \(\zeta(\mathbf{c^*}, \mathbf{d} - \mathbf{j}_1)\) is calculated for every player sequence \(\mathbf{j}_1 \in U^{(J,r)}\) and each dealer upcard \(\mathbf{c^*}\) (canonical vector). This is generalized in the matrix \(\zeta(\mathbf{C}, \mathbf{G}) \in \mathbb{R}^{|U^{(J,r)}| \times 10}\), defined as:

\[
\zeta(\mathbf{C}, \mathbf{G})_{ij} = \frac{\mathbf{G}_{ij}}{\sum_{k=1}^{10} \mathbf{G}_{ik}},
\]

where:

\begin{itemize}
    \item \(\mathbf{C} \in \mathbb{R}^{10 \times 10}\): Canonical matrix representing all possible dealer upcards.
    \item \(\mathbf{G} = \mathbf{d} - \mathbf{U}^{(J,r)}\): Matrix of available cards after adjusting the initial deck with the player sequences.
    \item \(\mathbf{G}_{ij}\): Element corresponding to row \(i\) (player sequence) and column \(j\) (specific upcard).
    \item \(\sum_{k=1}^{10} \mathbf{G}_{ik}\): Total sum of cards in row \(i\), ensuring normalization.
\end{itemize}

This result compactly describes the extraction probabilities for all combinations of player sequences and dealer upcards.

\subsubsection{Calculation of Mass Vectors for the Stand Decision}
\label{sec:calculo_vectores_plantarse}

Using the final return matrix \(\mathbf{R} \in \mathbb{R}^{|U^{(J,r)}| \times |U^{(C,r)}|}\), the probability matrix \(\mathbf{P} \in \mathbb{R}^{|U^{(J,r)}| \times |U^{(C,r)}|}\), and the multiplicity dictionary \(\mathbf{M}_{\text{up}}\), the mass vectors \( (M_{\theta_1}^s \mid \theta_1(s) = \text{Stand}) \) are calculated for all decision instances \( s \in \prescript{}{1}{S}^{(J, r, \mathbf{d})} \). The procedure is carried out by the \textbf{Vector\_Masa\_Stand()} function (Appendix \ref{sec:anexo5}), described as follows:

\begin{enumerate}
    \item \textbf{Identification of relevant sequences:} The indices in the matrix \(\mathbf{U}^{(J,r)}\) associated with the sequences \(\mathbf{j_1} \in U_3^{(J,r)}\) are identified and labeled as \(\mathbf{index}_3\).

    \item \textbf{Initialization of mass vector matrices:} Two matrices initialized with zeros are created:
    \[
    \mathbf{MV}_{\text{Stand}, 0} \in \mathbb{R}^{|U^{(J,r)}| \times 10 \times 6}, \quad
    \mathbf{MV}_{\text{Stand}, 1} \in \mathbb{R}^{|U_3^{(J,r)}| \times 10 \times 6}.
    \]

    \item \textbf{Filtering for the first upcard:} Starting with the ace as the upcard (\(\mathbf{c^*} = \mathbf{e_1}\)), the columns of the matrices \(\mathbf{R}\) and \(\mathbf{P}\) are filtered using the indices from the first column of the dictionary \(\mathbf{M}_{\text{up}}[\mathbf{e_1}]\). The resulting matrices are \(\mathbf{R}_1\) and \(\mathbf{P}_1\).

    \item \textbf{Adjustment of probabilities by multiplicities:} The columns of the matrix \(\mathbf{P}_1\) are multiplied by the values in the second column of \(\mathbf{M}_{\text{up}}[\mathbf{e_1}]\), which correspond to the multiplicities of the sequences associated with the filtered columns.

    \item \textbf{Normalization by extraction probabilities:} The rows of the matrix \(\mathbf{P}_1\) are divided vectorially by the probabilities \(\zeta(\mathbf{c^*}, \mathbf{d} - \mathbf{j_1}) = \zeta(\mathbf{C}, \mathbf{G})_{ij}\), according to the sequence \(\mathbf{j_1}\) associated with each row of \(\mathbf{P}_1\) and the upcard \(\mathbf{c^*}\). After this step, the values in the matrix \(\mathbf{P}_1\) represent:
    \[
    D_{\theta_1}(s', s) \mid \theta_1(s) = \text{Stand}.
    \]

    \item \textbf{Calculation of mass vectors (\(z = 0\)):} By iterating element by element through the values in \(\mathbf{P}_1\) and \(\mathbf{R}_1\), the mass vector is calculated for all decision instances \( s \in \prescript{}{1}{S}^{(J, r, \mathbf{d})} \) with \( z = 0 \). The results are stored in the matrix \(\mathbf{MV}_{\text{Stand}, 0}\).

    \item \textbf{Calculation of mass vectors (\(z = 1\)):} Filtering the rows of \(\mathbf{R}_1\) and \(\mathbf{P}_1\) using \(\mathbf{index}_3\), the mass vector is calculated for instances with \( z = 1 \), which is added to the matrix \(\mathbf{MV}_{\text{Stand}, 1}\).

    \item \textbf{Iteration for the 10 upcards:} The described process is repeated for all 10 possible upcards \(\mathbf{c^*}\), completing the matrices \(\mathbf{MV}_{\text{Stand}, 0}\) and \(\mathbf{MV}_{\text{Stand}, 1}\).

    \item \textbf{Organization into a dictionary:} Finally, the matrix \(\mathbf{MV}_{\text{Stand}, 0}\) is organized into a dictionary \(\mathbf{MV}_{\text{Stand}, 0, n}\), where \(\mathbf{MV}_{\text{Stand}, 0, n}[n]\) returns a matrix containing the mass vectors associated with decision instances \( s \in \prescript{}{1}{S}^{(J, r, \mathbf{d})} \), such that \(\mathbf{j_1} \in U_n^{(J,r)}\).
\end{enumerate}
\subsubsection{Hit, Stand, and Double: Exact Solution}
\label{sec:hit_stand_double}

The calculation of the mass vectors \( M_{\theta_i^*}^s \) and the optimal decisions \( \theta_i^*(s) \) is performed using a dynamic programming procedure. This approach propagates information recursively from the terminal levels (\( n = N \)) to the initial levels (\( n = 0 \)) of the Directed Acyclic Graph (DAG) \(\mathbf{G}^{(J,r)} = (U^{(J,r)}, A^{(J,r)})\):

\begin{enumerate}
    \item \textbf{Terminal Level (\( n = N \))}: The mass vectors are directly defined as:
    \[
    M_{\theta^*}^s = \mathbf{e}_2,
    \]

    \item \textbf{Intermediate Levels (\( n < N \))}: For each node \( s \in U_n^{(J, \text{rel})} \), the mass vector is calculated by propagating information from the nodes at the next level \( U_{n+1}^{(J, \text{rel})} \) and considering the transitions \( T_\theta(s', s) \). The optimal policy \( \theta_i^*(s) \) is obtained by maximizing the expected return as defined earlier in the theoretical framework:
    \[
    \theta_i^*(s) = \arg\max_{a \in \{\text{Hit}, \text{Stand}, \text{Double}\}} \mathbb{E}[Y_a^s],
    \]
    where the expected return under each action \( a \) is calculated as:
    \[
    V_i(s \mid \theta_i(s) = \text{Hit}) = \mathbb{E}[Y_{\theta_i}^s \mid \theta_i(s) = \text{Hit}] = \sum_{s' \in S_i^{(J, \mathbf{d}, n+1)}} T_{\theta_1}(s', s) \cdot M_{\theta_i^*}^{s'},
    \]
    for \( a = \text{Hit} \), and for the action \( \text{Double} \) (with \( z = 1 \)):
    \[
    \mathbb{E}[Y_{\theta_1}^s \mid \theta_1(s) = \text{Double}] = \sum_{\substack{s' \in S_1^{(J, \mathbf{d}, 3)} \\ z = 1}} T_{\theta_1}(s', s) \cdot M_{\theta_1^*}^{s'}.
    \]
\end{enumerate}

The procedure progresses level by level. First, the mass vectors \( M_{\theta_2^*}^s \) and decisions \( \theta_2^*(s) \) are calculated simultaneously for all instances \( s \in S_2^{(J, \mathbf{d}, n)} \) at each level \( n \). Subsequently, at \( n = 2 \), the mass vectors \( M_{\theta_1^*}^s \) and decisions \( \theta_1^*(s) \) are calculated for \( s \in S_1^{(J, \mathbf{d}, n)} \), including the \textit{Double} action. For levels \( n > 2 \) in \( \theta_1 \), where the \textit{Double} action is not allowed, the optimal policy \( \theta_1^* \) coincides with \( \theta_2^* \). Therefore, at these levels, the calculations for \( M_{\theta_1^*}^s \) and \( \theta_1^*(s) \) are equivalent to those performed for \( M_{\theta_2^*}^s \) and \( \theta_2^*(s) \).

This model generates the mass vectors \( M_{\theta_i^*}^s \) for all instances \( s \in S_i^{(J, \mathbf{d}, n)} \) and the decisions \( \theta_i^*(s) \), maximizing the expected return under the problem's constraints.

The complete pseudocode of the implementation is presented in the \textbf{Optimizacion\_Problema\_Restringido()} function (Appendix \ref{sec:anexo6}). This function details the step-by-step procedure for optimizing the policies \( \theta_1 \) and \( \theta_2 \), applying dynamic programming to efficiently compute the value functions \( V^*(s) \), the optimal decisions \( \theta_i^*(s) \), and the mass vectors \( M_{\theta_i^*}^s \) at each level of the DAG.
\subsection{Code Appendices}
\subsubsection{Appendix 1: Sequence Score Calculation}

\label{sec:anexo1}
\begin{lstlisting}[
    language=Python,
    caption={Function calculate\_score},
    label={lst:calculate_score},
    texcl=true,                   % Enables LaTeX comments
    basicstyle=\small\sf,         % Base font style
    commentstyle=\small\rm,       % Font style for comments
    mathescape=true,              % Allows math expressions
    escapeinside={(*)}{*)},        % Configures escape for LaTeX
    breaklines=true,              % Enables automatic line breaks
    breakatwhitespace=true        % Breaks at whitespace
]
def calculate_score(sequence):

    indices = np.arange(1, 11)   # [1, 2, 3, ..., 10]
    base_score = np.dot(sequence, indices)

    # Blackjack rule
    if sequence[0] > 0 and sequence[9] > 0 and base_score == 11:
        return 100  # Blackjack

    # Adjustment for Aces
    if sequence[0] > 0 and base_score + 10 <= 21:
        base_score += 10

    return base_score
\end{lstlisting}

\subsubsection{Appendix 2: Sequence and Edge Calculation}
\label{sec:anexo2}
\begin{lstlisting}[
    language=Python,
    caption={Generation of sequences and edges},
    label={lst:generate_sequences_edges},
    texcl=true,                   % Enables LaTeX comments
    basicstyle=\small\sf,         % Base font style
    commentstyle=\small\rm,       % Font style for comments
    mathescape=true,              % Allows math expressions
    escapeinside={(*)}{*)},        % Configures escape for LaTeX
    breaklines=true,              % Enables automatic line breaks
    breakatwhitespace=true        % Breaks at whitespace
]
def traverse_graph_sequences(player_or_dealer, state, state_list, edge_list, visited, expanded):
    """
    Performs depth-first search (DFS) from 'state'.
    - state: np.array of 10 positions.
    - state_list: accumulation of generated states (sequences).
    - edge_list: accumulation of edges [parent_idx, child_idx].
    - visited: dict {tuple(state): idx}, to avoid repeating states.
    - expanded: set of indices already expanded to avoid re-expansion.

    Returns the index (in state_list) of the current 'state'.
    """
    # Convert to tuple to use as a key in the dictionary
    t_state = tuple(state)
    if t_state not in visited:
        idx_current = len(state_list)
        visited[t_state] = idx_current
        state_list.append(state.copy())
    else:
        idx_current = visited[t_state]

    if idx_current in expanded:
        return idx_current

    expanded.add(idx_current)

    current_score = calculate_score(state)
    max_point = 21 if player_or_dealer == 'player' else 16

    if current_score > max_point or current_score == 100:
        return idx_current

    for i in range(10):
        child = state.copy()
        child[i] += 1
        child_score = calculate_score(child)

        t_child = tuple(child)
        if t_child not in visited:
            idx_child = len(state_list)
            visited[t_child] = idx_child
            state_list.append(child)
        else:
            idx_child = visited[t_child]

        edge_list.append([idx_current, idx_child])

        if child_score <= max_point and child_score != 100:
            traverse_graph_sequences(player_or_dealer, child, state_list, edge_list, visited, expanded)

    return idx_current

def generate_sequences_and_edges(player_or_dealer):

    initial_state = np.zeros(10, dtype=int)

    state_list = []
    edge_list = []
    visited = {}
    expanded = set()

    traverse_graph_sequences(player_or_dealer, initial_state, state_list, edge_list, visited, expanded)

    sequence_matrix = np.array(state_list, dtype=int)
    edge_matrix = np.array(edge_list, dtype=int)

    return sequence_matrix, edge_matrix
\end{lstlisting}

\subsubsection{Appendix 3: Calculation of Card Sequence Probabilities}
This code calculates the extraction probabilities of card sequences from a given deck.
\label{sec:anexo3}
\begin{lstlisting}[
    language=Python,
    caption={Function probs\_sequences},
    label={lst:probs_sequences},
    texcl=true,
    basicstyle=\small\sf,
    commentstyle=\small\rm,
    breaklines=true,
    breakatwhitespace=true
]
def probs_sequences(sequence_matrix, deck):
    """
    Calculates the extraction probabilities of card sequences from a deck.
    :param sequence_matrix: Nx10 matrix containing the extracted card sequences.
    :param deck: Vector of size 10 representing the number of cards available for each type in the deck.
    :return: Vector of size N with the calculated probabilities for each sequence.
    """
    numerators = np.ones((10, 20))
    denominators = np.ones(20)

    total_cards = np.sum(deck)
    divisor = 1.0

    for i in range(1, 20):
        if total_cards - i + 1 > 0:
            divisor *= (total_cards - i + 1)
            denominators[i] = divisor
        else:
            denominators[i] = 1

    for card in range(10):
        remaining_cards = deck[card]
        numerator = 1.0
        for i in range(1, 20):
            if remaining_cards - i + 1 > 0:
                numerator *= (remaining_cards - i + 1)
                numerators[card, i] = numerator
            else:
                numerators[card, i] = 0

    sequence_numerators = numerators[np.arange(10), sequence_matrix.T].T
    sequence_denominators = denominators[np.sum(sequence_matrix, axis=1).astype(int)]

    sequence_probabilities = np.prod(sequence_numerators, axis=1) / sequence_denominators

    return sequence_probabilities
\end{lstlisting}

\subsubsection{Appendix 4: Updating Probabilities in Sequence Trees}
\label{sec:anexo4}
This recursive function adjusts the probability vector in the player's sequence tree.

\begin{lstlisting}[
    language=Python,
    caption={Function update\_probs\_matrix},
    label={lst:update_probs_matrix},
    texcl=true,
    basicstyle=\small\sf,
    commentstyle=\small\rm,
    mathescape=true,
    escapeinside={(*)}{*)},
    breaklines=true,
    breakatwhitespace=true
]
def update_probs_matrix(origin_deck, PLAYER_SAFE, PlayerChildren_SAFE_GEN, DEALER_TERMINAL, sequence_index=root,
                         Probs_Matrix=None, Probs_Vector=None, Bayes_decks=None):
    """
    Recursive function that updates and records the probability vector for extracting the dealer's terminal sequences,
    traversing the player's sequence DAG (Directed Acyclic Graph).

    :param root: Index of the root sequence in the player's DAG
    :param PLAYER_SAFE: Matrix containing the player's sequences that do not exceed 21
    :param PlayerChildren_SAFE_GEN: Dictionary containing the child sequences for each sequence
    :param DEALER_TERMINAL: Matrix of the dealer's terminal sequences
    :param origin_deck: origin deck of the round
    :param sequence_index: Index of the current sequence in the player's tree
    :param Probs_Matrix: Matrix containing the probabilities of extracting the dealer's terminal sequences 
                         from all sequences in PLAYER_SAFE
    :param Probs_Vector: Vector with extraction probabilities of the dealer's terminal sequences for the 
                         current player's sequence
    :param Bayes_decks: Auxiliary deck matrix used to adjust the probability vector using Bayes' theorem
    :return: Probs_Matrix
    """
    if sequence_index == root:
        Probs_Matrix = np.zeros((PLAYER_SAFE.shape[0], DEALER_TERMINAL.shape[0]))
        Bayes_decks = origin_deck - DEALER_TERMINAL
        Probs_Vector = probs_sequences(DEALER_TERMINAL, origin_deck)
        Probs_Matrix[root] = Probs_Vector

    for child in PlayerChildren_SAFE_GEN[sequence_index]:
        Bayes_decks_copy = Bayes_decks.copy()
        origin_deck_copy = origin_deck.copy()

        card = np.where((PLAYER_SAFE[child] - PLAYER_SAFE[sequence_index]) == 1)

        numerator_adj = (Bayes_decks_copy[:, card] / np.sum(Bayes_decks_copy, axis=1, keepdims=True)).flatten()
        denominator_adj = (origin_deck_copy[card] / np.sum(origin_deck_copy, keepdims=True)).flatten()
        adjusted_vector = Probs_Vector * numerator_adj / denominator_adj

        Probs_Matrix[child] = adjusted_vector

        origin_deck_copy[card] += -1
        Bayes_decks_copy[:, card] += -1

        update_probs_matrix(origin_deck_copy, PLAYER_SAFE, PlayerChildren_SAFE_GEN, DEALER_TERMINAL, child, Probs_Matrix, Probs_Vector, Bayes_decks)

    return Probs_Matrix
\end{lstlisting}

\subsubsection{Appendix 5: Calculation of the Mass Vector for the Stand State}
\label{sec:anexo5}
This function calculates the probabilities of winning, tying, and losing for the Stand state.

\begin{lstlisting}[
    language=Python,
    caption={Function Stand\_Mass\_Vector},
    label={lst:Stand_Mass_Vector},
    texcl=true,
    basicstyle=\small\sf,
    commentstyle=\small\rm,
    mathescape=true,
    escapeinside={(*)}{*)},
    breaklines=true,
    breakatwhitespace=true
]
def Stand_Mass_Vector(origin_deck, PLAYER_SAFE, Win_Loss_Matrix, MULT_INDICES, Probs_Matrix):
    Stand_Mass_Vector = np.zeros((Probs_Matrix.shape[0], 10, 6))
    adjusted_deck = origin_deck - PLAYER_SAFE
    adjusted_deck_sum = np.sum(adjusted_deck, axis=1)

    BJ_indices = np.where(PLAYER_SAFE is Bj)
    noBJ_indices = np.where(PLAYER_SAFE is not Bj)

    for upcard in range(1, 11): 
        upcard_indices = MULT_INDICES[upcard][:,0]
        upcard_multiplicities = MULT_INDICES[upcard][:,1]

        WinLoss_upcard = Win_Loss_Matrix[:, upcard_indices]
        Probs_upcard = Probs_Matrix[:, upcard_indices] * upcard_multiplicities

        Probs_upcard /= (adjusted_deck[:, upcard-1] / adjusted_deck_sum)[:, np.newaxis]

        prob_win = np.sum(Probs_upcard * (WinLoss_upcard == 1), axis=1)
        prob_tie = np.sum(Probs_upcard * (WinLoss_upcard == 0), axis=1)
        prob_lose = np.sum(Probs_upcard * (WinLoss_upcard == -1), axis=1)

        prob_dealer_bust = 1 - (prob_win + prob_tie + prob_lose)

        Stand_Mass_Vector[noBJ_indices, upcard-1, 3] = prob_win[noBJ_indices] + prob_dealer_bust[noBJ_indices]
        Stand_Mass_Vector[noBJ_indices, upcard-1, 2] = prob_tie[noBJ_indices]
        Stand_Mass_Vector[noBJ_indices, upcard-1, 1] = prob_lose[noBJ_indices]

        Stand_Mass_Vector[BJ_indices, upcard-1, 4] = prob_win[BJ_indices] + prob_dealer_bust[BJ_indices]
        Stand_Mass_Vector[BJ_indices, upcard-1, 2] = prob_tie[BJ_indices]

    Mass_Vector_Dict = [Stand_Mass_Vector[np.sum(PLAYER_SAFE, 1) == i] for i in range(22)]

    return Mass_Vector_Dict
\end{lstlisting}

\subsubsection{Appendix 6: Optimization of Restricted Round Policies}
This code implements a procedure for optimizing round policies using a decision tree.
\label{sec:anexo6}
\begin{lstlisting}[
    language=Python,
    caption={Function Optimize\_Restricted\_Policy},
    label={lst:Optimize_Restricted_Policy},
    basicstyle=\small\ttfamily,
    mathescape=true,
    escapeinside={(*)}{*)},
    breaklines=true
]
def Optimize_Restricted_Policy(origin_deck, PLAYER_SAFE_BY_SIZE, Stand_Mass_Vector_Dict, 
                                n=0, Optimal_Mass_Vector_Dict={}):
    # Terminal level
    if n == 22:
        return [0,1,0,0,0,0], Optimal_Mass_Vector_Dict

    # Retrieve player sequences of size n
    player_sequences_n = PLAYER_SAFE_BY_SIZE[n]

    # Define mass vectors for hit and stand decisions
    mass_vector_hit = np.zeros([len(player_sequences_n),10,6])
    mass_vector_stand = Stand_Mass_Vector_Dict[n]

    # Recursively retrieve optimal mass vector for n+1
    future_mass_vector, Optimal_Mass_Vector_Dict = Optimize_Restricted_Policy(
        origin_deck, PLAYER_SAFE_BY_SIZE, Stand_Mass_Vector_Dict, n+1, Optimal_Mass_Vector_Dict)

    # Calculate transition probabilities to 10 possible children per sequence
    transition_probs = transition_probabilities(origin_deck, player_sequences_n)

    # Compute hit mass vector using future optimal mass vector
    mass_vector_hit = weight_children_by_transition(transition_probs, future_mass_vector)

    optimal_mass_vector, optimal_decisions = max_expectation_filter(mass_vector_hit, mass_vector_stand)

    # Save optimal values for this n
    Optimal_Mass_Vector_Dict[n] = optimal_mass_vector, optimal_decisions

    return optimal_mass_vector, Optimal_Mass_Vector_Dict
\end{lstlisting}

\section{Proofs}
\subsection{Theorem 1}
\label{teo1}
\textit{The optimal policy under a utility function with constant relative risk aversion is independent of the state variable “wealth”}:
\[
\pi_{\theta, u_1, H}^*(\mathbf{d}, w, n) = \pi_{\theta, u_1, H}^*(\mathbf{d}, w', n), \quad \forall w, w' \in \mathbb{R}_+
\]

\noindent{Where:}
\begin{enumerate}
    \item The utility function \(u_1\) represents constant relative risk aversion:
    \[
    u_1(w; \alpha) = \begin{cases}
    \frac{w^\alpha}{\alpha} &, \text{if } \alpha > 0 \\
    \ln(w) &, \text{if } \alpha = 0
    \end{cases}
    \]
\end{enumerate}

With:
\begin{itemize}
    \item \(\alpha \in \mathbb{R}_+\).
\end{itemize}

\subsubsection{Proof of Theorem 1}
To prove the theorem, we will recursively establish the independence of the optimal policy with respect to wealth, starting from terminal states and progressing to the initial state. For this, we use the recursive formulation of the optimization problem:

\[
\pi^*(\psi_n=\psi)
\;=\;
\arg \max_{b \in [0, 0.5]}
\sum\limits_{\psi' \in \Gamma_{[\psi, b]}}
T_b(\psi', \psi)\,V^*(\psi'), \quad \forall \psi \in \mathcal{B} \mid {n<H}
\]

\noindent{Where:}
\begin{enumerate}
    \item \(V^*(\psi)\) represents the \textbf{optimal value} of a state $\psi$, and is calculated as follows:
    \[
    V^*(\psi_{n} = \psi) \;=\;
    \begin{cases}
    \displaystyle\max_{b \in [0, 0.5]}
    \sum\limits_{\psi' \in \Gamma_{[\psi, b]}}
    T_b(\psi', \psi)\,V^*(\psi'),
    & \text{if } n<H
    \\[10pt]
    u_1(w; \alpha)    & \text{if } n=H
    \end{cases}
    \]
    
    \item The transition function $ T_b(\psi', \psi) $ is expressed as:
    \[
    T_b(\psi', \psi) = \mathbb{P}(\psi_{n+1} = \psi' \mid \psi_n = \psi, \pi(\psi) = b) = \begin{cases} 
    \mathbb{P}(\mathbf{d}_{n+1} = \mathbf{d}', X_\theta^{\mathbf{d}} = x \mid \mathbf{d}_n = \mathbf{d}), & \text{if } \psi' \in \Gamma_{[\psi, b]} \\
    0, & \text{if } \psi' \notin \Gamma_{[\psi, b]}
    \end{cases}
    \]

    \item $\Gamma_{[\psi, b]}$ represents the set of states \( \psi' \in \mathcal{B} \) to which the player can transition from a state \( \psi \) after making a bet \( b \), referred to as the \textit{available space}:
    \[
    \Gamma_{[\psi, b]} = \{(\mathbf{d}^*, w^*, n+1) \mid \mathbf{d}^* \in (\mathcal{G}^\mathbf{d} \cup \bar{\mathbf{d}}), w^* = w \cdot (1 + b \cdot r), r \in \mathcal{R} \}
    \]
    
\end{enumerate}

With:
\begin{itemize}
    \item \( \psi = (\mathbf{d}, w, n) \in \mathcal{B} \): Current state.
    \item \( \psi' = (\mathbf{d'}, w', n+1) \in \mathcal{B} \): Transition state.
    \item \(x = \frac{w' - w}{w \cdot b} \in \mathcal{R}\): Return associated with the state transition.
\end{itemize}

We proceed by recursively developing the calculation of $V^*( \psi)$. Subsequently, we conclude that the policy is independent of wealth for every state $\psi \in \mathcal{B} \mid {n<H}$.

\paragraph{Explicit Transition Function}

First, we note that the calculation of the value \(V^*(\psi)\) only considers transitions to states \( \psi' \in \Gamma_{[\psi,b]} \). For these states, the transition function is independent of wealth $w$, and simplifies to the following expression:

\[
T_b(\psi', \psi)  =
\mathbb{P}(\mathbf{d}_{n+1} = \mathbf{d}', X_\theta^{\mathbf{d}} = x \mid \mathbf{d}_n = \mathbf{d})
\]

\paragraph{Value Function: Simplification}

Considering this expression for the transition function and decomposing the summation in the value equation based on the space \( \Gamma_{[\psi, b]} \), the formula for calculating the optimal state value becomes:

\[
V^*(\psi_{n} = \psi) \;=\;
\begin{cases}
\displaystyle\max_{b \in [0, 0.5]}
\sum\limits_{x \in \mathcal{R}} \sum\limits_{\mathbf{d}' \in (Q^{\mathbf{d}} \cup \bar{\mathbf{b}})}
\mathbb{P}(\mathbf{d}_{n+1} = \mathbf{d}', X_\theta^{\mathbf{d}} = x \mid \mathbf{d}_n = \mathbf{d})\,V^*( w \cdot (1 + b \cdot x)),
& \text{if } n<H
\\[10pt]
u_1(w; \alpha)    & \text{if } n=H
\end{cases}
\]

We proceed with the proof for \( \alpha > 0 \), followed by \( \alpha = 0 \).
\paragraph{\textbf{Case \( \alpha > 0 \)}: Proof}

\paragraph{Period \( H-1 \)}

We proceed to calculate the optimal value for terminal decision states \( \psi_{H-1} = \psi = (\mathbf{d}, w, H-1) \). The value function for terminal states \( V^*(\psi_{H-1} = \psi) \) is expressed as:

\[
V^*(\psi_{H-1} = \psi) = \frac{w^\alpha}{\alpha} \cdot U^*(\psi_{H-1} = \psi),
\]

\noindent{Where:}
\begin{enumerate}
    \item \(U^*(\psi)\) represents the \textbf{auxiliary optimal value} of a state $\psi$, calculated as:
    \[
    U^*(\psi_{H-1} = \psi) = \max_{b \in [0, 0.5]} \sum\limits_{x \in \mathcal{R}} \sum\limits_{\mathbf{d}' \in (Q^{\mathbf{d}} \cup \bar{\mathbf{b}})}
    \mathbb{P}(\mathbf{d}_H = \mathbf{d}', X_\theta^{\mathbf{d}} = x \mid \mathbf{d}_{H-1} = \mathbf{d})\,(1 + b \cdot x)^\alpha.
    \]
\end{enumerate}

\paragraph{Period \( H-2 \)}

We now calculate the value for states \( \psi_{H-2} = \psi = (\mathbf{d}, w, H-2) \). Analogous to the \( H-1 \) period, the state value is expressed in terms of the auxiliary state value as follows:

\[
V^*(\psi_{H-2} = \psi) = \frac{w^\alpha}{\alpha} \cdot U^*(\psi_{H-2} = \psi),
\]

\noindent{Where:}
\begin{enumerate}
    \item The auxiliary value for the \( H-2 \) period is given by:
    \[
    U^*(\psi_{H-2} = \psi) = \max_{b \in [0, 0.5]} \sum\limits_{x \in \mathcal{R}} \sum\limits_{\mathbf{d}' \in (Q^{\mathbf{d}} \cup \bar{\mathbf{b}})}
    \mathbb{P}(\mathbf{d}_{H-1} = \mathbf{d}', X_\theta^{\mathbf{d}} = x \mid \mathbf{d}_{H-2} = \mathbf{d})\,(1 + b \cdot x)^\alpha \cdot U^*(\psi').
    \]
\end{enumerate}

\paragraph{Generalization}

The process described above can be iteratively repeated until the initial decision state \( \psi_0 \), resulting in a general recursive formulation for the state value calculation:

\[
V^*(\psi_n = \psi) = \begin{cases}
\frac{w^\alpha}{\alpha}, & \text{if } n = H, \\
\frac{w^\alpha}{\alpha} \cdot U^*(\psi), & \text{if } n < H,
\end{cases}
\]

\noindent{Where:}
\begin{enumerate}
    \item The auxiliary value for the \( n \)-th period is recursively calculated as:
    \[
    U^*(\psi_n = \psi) = \begin{cases}
    \max_{b \in [0, 0.5]} \sum\limits_{x \in \mathcal{R}} \sum\limits_{\mathbf{d}' \in (Q^{\mathbf{d}} \cup \bar{\mathbf{b}})}
    \mathbb{P}(\mathbf{d}_{n+1} = \mathbf{d}', X_\theta^{\mathbf{d}} = x \mid \mathbf{d}_n = \mathbf{d})\,(1 + b \cdot x)^\alpha \cdot U^*(\psi'), & \text{if } n < H, \\
    1, & \text{if } n = H.
    \end{cases}
    \]
\end{enumerate}

With:
\begin{itemize}
    \item \( \psi = (\mathbf{d}, w, n) \in \mathcal{B} \),
    \item \( \psi' = (\mathbf{d'}, w \cdot (1 + b \cdot x), n+1) \in \mathcal{B} \).
\end{itemize}

\paragraph{Policy Conclusion}

Evaluating the expression for the optimal value \( V^*(\psi_n = \psi) \) within the equation defining the optimal policy, we conclude that the optimal policy is given by:

\[
\pi^*(\psi_n = \psi) =
\begin{cases}
\arg\max_{b \in [0, 0.5]}
\sum\limits_{x \in \mathcal{R}} \sum\limits_{\mathbf{d}' \in (Q^{\mathbf{d}} \cup \bar{\mathbf{b}})}
\mathbb{P}(\mathbf{d}_{n+1} = \mathbf{d}', X_\theta^{\mathbf{d}} = x \mid \mathbf{d}_n = \mathbf{d})\,(1 + b \cdot x)^\alpha \cdot U^*(\psi'), & \text{if } n < H-1, \\
\arg\max_{b \in [0, 0.5]}
\sum\limits_{x \in \mathcal{R}} \sum\limits_{\mathbf{d}' \in (Q^{\mathbf{d}} \cup \bar{\mathbf{b}})}
\mathbb{P}(\mathbf{d}_{n+1} = \mathbf{d}', X_\theta^{\mathbf{d}} = x \mid \mathbf{d}_n = \mathbf{d})\,(1 + b \cdot x)^\alpha, & \text{if } n = H-1.
\end{cases}
\]

Observing that the auxiliary value \( U^*(\psi_n = \psi) \) is independent of wealth \( w \) for all \( n \leq H \), we conclude that the optimal policy is also independent of wealth, completing the proof for \( \alpha > 0 \):

\[
\pi^*(\mathbf{d}, w, n) = \pi_{\theta, u_1, H}^*(\mathbf{d}, w', n), \quad \forall w, w' \in \mathbb{R}_+, \alpha > 0.
\]

\paragraph{\textbf{Case \( \alpha = 0 \)}: Proof}

\paragraph{Period \( H-1 \)}

We calculate the optimal value for terminal decision states \( \psi_{H-1} = \psi = (\mathbf{d}, w, H-1) \). The value function for terminal states \( V^*(\psi_{H-1} = \psi) \) is expressed as:

\[
V^*(\psi_{H-1} = \psi) = \ln(w) + U^*(\psi_{H-1} = \psi),
\]

\noindent{Where:}
\begin{enumerate}
    \item For this case, the auxiliary value function is:
    \[
    U^*(\psi_{H-1} = \psi) = \max_{b \in [0, 0.5]}
    \sum\limits_{x \in \mathcal{R}} \sum\limits_{\mathbf{d}' \in (Q^{\mathbf{d}} \cup \bar{\mathbf{b}})}
    \mathbb{P}(\mathbf{d}_H = \mathbf{d}', X_\theta^{\mathbf{d}} = x \mid \mathbf{d}_{H-1} = \mathbf{d})\,\ln(1 + b \cdot x).
    \]
\end{enumerate}

\paragraph{Generalization}

This process can be identically repeated until the initial decision state \( \psi_0 \), yielding a general recursive formulation for the state value:

\[
V^*(\psi_n = \psi) = \begin{cases}
\ln(w), & \text{if } n = H, \\
\ln(w) + U^*(\psi), & \text{if } n < H.
\end{cases}
\]
\noindent{Where:}
\begin{enumerate}
    \item The auxiliary value for the \(n\)-th period is recursively calculated as:
    \[
    U^*(\psi_n = \psi) = \begin{cases}
    \displaystyle \max_{b \in [0, 0.5]}
    \sum\limits_{x \in \mathcal{R}} \sum\limits_{\mathbf{d}' \in (Q^{\mathbf{d}} \cup \bar{\mathbf{b}})}
    \mathbb{P}(\mathbf{d}_{n+1} = \mathbf{d}', X_\theta^{\mathbf{d}} = x \mid \mathbf{d}_n = \mathbf{d})\,\big(\ln(1 + b \cdot x) + U^*(\psi')\big), & \text{if } n < H, \\
    0, & \text{if } n = H.
    \end{cases}
    \]
\end{enumerate}

With:
\begin{itemize}
    \item \( \psi = (\mathbf{d}, w, n) \in \mathcal{B} \),
    \item \( \psi' = (\mathbf{d'}, w \cdot (1 + b \cdot x), n+1) \in \mathcal{B} \).
\end{itemize}

\paragraph{Policy Conclusion}

Evaluating the expression obtained for the optimal value \(V^*(\psi_n = \psi)\) within the equation defining the optimal policy, we conclude that the optimal policy is given by:

\[
\pi^*(\psi_n = \psi) =
\begin{cases}
\displaystyle
\arg\max_{b \in [0, 0.5]}
\sum\limits_{x \in \mathcal{R}} \sum\limits_{\mathbf{d}' \in (Q^{\mathbf{d}} \cup \bar{\mathbf{b}})}
\mathbb{P}(\mathbf{d}_{n+1} = \mathbf{d}', X_\theta^{\mathbf{d}} = x \mid \mathbf{d}_n = \mathbf{d})\,\big(\ln(1 + b \cdot x) + U^*(\psi')\big), & \text{if } n < H-1, \\[10pt]
\displaystyle
\arg\max_{b \in [0, 0.5]}
\sum\limits_{x \in \mathcal{R}} \sum\limits_{\mathbf{d}' \in (Q^{\mathbf{d}} \cup \bar{\mathbf{b}})}
\mathbb{P}(\mathbf{d}_{n+1} = \mathbf{d}', X_\theta^{\mathbf{d}} = x \mid \mathbf{d}_n = \mathbf{d})\,\ln(1 + b \cdot x), & \text{if } n = H-1.
\end{cases}
\]

Observing that the auxiliary value \(U^*(\psi_n = \psi)\) is independent of wealth \(w\) for all \(n \leq H\), we conclude that the optimal policy is also independent of wealth. Thus, the proof for the case \(\alpha = 0\) is completed, finalizing the proof of the theorem:

\[
\pi_{\theta, u_1, H}^*(\mathbf{d}, w, n) = \pi_{\theta, u_1, H}^*(\mathbf{d}, w', n), \quad \forall w, w' \in \mathbb{R}_+, \alpha = 0.
\]

\paragraph{Final Conclusion for Theorem 1}

Combining the results for both cases \(\alpha > 0\) and \(\alpha = 0\), we conclude that the optimal policy under a utility function with constant relative risk aversion is independent of wealth:

\[
\pi_{\theta, u_1, H}^*(\mathbf{d}, w, n) = \pi_{\theta, u_1, H}^*(\mathbf{d}, w', n), \quad \forall w, w' \in \mathbb{R}_+.
\]
\subsection{Theorem 2}
\label{teo2}
\textit{The optimal policy under a utility function with constant relative risk aversion, and for an infinite horizon of betting rounds, is independent of the state variable “number of rounds played”}:

\[
\pi_{\theta, u_1, H \to \infty}^*(\mathbf{d}, w, n) = \pi_{\theta, u_1, H \to \infty}^*(\mathbf{d}, w, n'), \quad \forall n, n' \in \mathbb{Z}_+
\]

\noindent{Where:}
\begin{enumerate}
    \item The utility function \(u_1\) represents constant relative risk aversion:
    \[
    u_1(w; \alpha) = \begin{cases}
    \frac{w^\alpha}{\alpha} &, \text{if } \alpha > 0 \\
    \ln(w) &, \text{if } \alpha = 0
    \end{cases}
    \]
\end{enumerate}

With:
\begin{itemize}
    \item \(\alpha \in \mathbb{R}_+\).
\end{itemize}

\subsubsection{Proof of Theorem 2}
The proof is based on the formulations obtained for the value function in the proof of Theorem 1. In particular, by taking the limit \(H \to \infty\), it is deduced that the value function satisfies the following relation:

\[
V^*(\psi_n = \psi) = \begin{cases}
\frac{w^\alpha}{\alpha} \cdot U^*(\psi), & \text{if } \alpha >0 \\
\ln(w) + U^*(\psi), & \text{if } \alpha = 0 
\end{cases}
\]

\noindent{Where:}
\begin{enumerate}
    \item The auxiliary value for the \(n\)-th period is recursively calculated as:
    \[
    U^*(\psi_n = \psi) = \begin{cases}
    \displaystyle \max_{b \in [0, 0.5]}
    \sum\limits_{x \in \mathcal{R}} \sum\limits_{\mathbf{d}' \in (Q^{\mathbf{d}} \cup \bar{\mathbf{b}})}
    \mathbb{P}(\mathbf{d}_{n+1} = \mathbf{d}', X_\theta^{\mathbf{d}} = x \mid \mathbf{d}_{n} = \mathbf{d})\,(1 + b \cdot x)^\alpha \cdot U^*(\psi'), & \text{if } \alpha > 0 \\
    \max_{b \in [0, 0.5]}
    \sum\limits_{x \in \mathcal{R}} \sum\limits_{\mathbf{d}' \in (Q^{\mathbf{d}} \cup \bar{\mathbf{b}})}
    \mathbb{P}(\mathbf{d}_{n+1} = \mathbf{d}', X_\theta^{\mathbf{d}} = x \mid \mathbf{d}_{n} = \mathbf{d})\,(\ln(1 + b \cdot x) + U^*(\psi')), & \text{if } \alpha = 0
    \end{cases}
    \]
\end{enumerate}

Evaluating the expression obtained for the optimal value \(V^*(\psi_n = \psi)\) within the equation defining the optimal policy, the optimal policy is determined as follows:

\[
\pi^*(\psi_n = \psi) =
\begin{cases}
\displaystyle
\arg\max_{b \in [0, 0.5]}
\sum\limits_{x \in \mathcal{R}} \sum\limits_{\mathbf{d}' \in (Q^{\mathbf{d}} \cup \bar{\mathbf{b}})}
\mathbb{P}(\mathbf{d}_{n+1} = \mathbf{d}', X_\theta^{\mathbf{d}} = x \mid \mathbf{d}_{n} = \mathbf{d})\,(1 + b \cdot x)^\alpha \cdot U^*(\psi'),
& \text{if } \alpha > 0 \\[10pt]
\displaystyle
\arg\max_{b \in [0, 0.5]}
\sum\limits_{x \in \mathcal{R}} \sum\limits_{\mathbf{d}' \in (Q^{\mathbf{d}} \cup \bar{\mathbf{b}})}
\mathbb{P}(\mathbf{d}_{n+1} = \mathbf{d}', X_\theta^{\mathbf{d}} = x \mid \mathbf{d}_{n} = \mathbf{d})\,(\ln(1 + b \cdot x) + U^*(\psi')),
& \text{if } \alpha = 0
\end{cases}
\]

Observing that the probability \(\mathbb{P}(\mathbf{d}_{n+1} = \mathbf{d}', X_\theta^{\mathbf{d}} = x \mid \mathbf{d}_{n} = \mathbf{d})\) does not depend on the round \(n\), it follows that the auxiliary value \(U^*(\psi_n = \psi)\) is also independent of the number of rounds played. Consequently, it is concluded that the optimal policy \(\pi^*(\psi_n = \psi)\) is independent of the variable \(n\), thus completing the proof of the theorem.
\subsection{Theorem 3}
\textit{Under a partial policy \(\omega\) and a counting system \( \kappa\), the long-term distribution of utility under a CRRA utility function converges as follows:}

\[
\lim_{N \to \infty} u_1(W_{N+1}; \alpha) \sim 
\begin{cases} 
\text{Normal}(N \cdot \mu_\omega, N \cdot \sigma_\omega^2) &, \text{if } \alpha = 0 \\
\text{LogNormal}(N \cdot \alpha \cdot \mu_\omega - \ln(\alpha), N \cdot \alpha^2 \cdot \sigma_\omega^2) &, \text{if } \alpha > 0
\end{cases}
\]

\noindent{Where:}
\begin{enumerate}
    \item The utility function \(u_1\) represents constant relative risk aversion:
    \[
    u_1(w; \alpha) = \begin{cases}
    \frac{w^\alpha}{\alpha} &, \text{if } \alpha > 0 \\
    \ln(w) &, \text{if } \alpha = 0
    \end{cases}
    \]
\end{enumerate}

With:
\begin{itemize}
    \item \( L_\omega = \ln(1 + \omega(\kappa(\mathbf{D})) \cdot X_\theta^{(\mathbf{D}|C)}). \)
    \item \( \mu_\omega = \mathbb{E}[L_\omega]. \)
    \item \( \sigma_\omega^2 = \text{Var}[L_\omega]. \)
\end{itemize}

\subsubsection{Proof of Theorem 3}
\label{teo3}
The proof is based on the stationary variables associated with the evolution of the deck \( \mathcal{D}_\theta = \{\mathbf{d}_n\}_{n=0}^\infty \) and the evolution of the count \( \mathcal{C}_{(\theta, \kappa)} = \{c_n= \kappa(\mathbf{d}_n)\}_{n=0}^\infty \).

\paragraph{Stationary Distribution of the Deck Evolution}

The deck evolution \(\mathcal{D}_\theta\) returns to the initial full deck \(\bar{\mathbf{d}}\) every time the penetration point is reached, and thus is a positive recurrent process. This implies the existence of a unique stationary distribution associated with the evolution \(\mathcal{D}_\theta\). We define the \textit{stationary deck} variable \( \mathbf{D} = \mathbf{D}_\theta \in \Omega_{(\tau,l)} \), which models the frequency of occurrence of decks \( \mathbf{d} \in \Omega_{(\tau,l)} \) in the process \( \mathcal{D}_\theta = \{\mathbf{d}_n\}_{n=0}^\infty \) in the long run:

\[
\mathbb{P}_\mathbf{D}(\mathbf{d}) = \mathbb{P}(\mathbf{D_\theta} = \mathbf{d}) = \lim_{N \to \infty} \frac{\sum_{n=0}^N 1_{\{\mathbf{d}^{obs}_n = \mathbf{d}\}}}{N+1}
\]

With:
\begin{itemize}
    \item \( \mathcal{D}_\theta^{z=\infty} = \{\mathbf{d}^{obs}_n\}_{n=0}^{z=\infty} \): A specific realization of the deck evolution up to round \(z=\infty\), where each \(\mathbf{d}^{obs}_n\) corresponds to a particular realization of the random deck \(\mathbf{d}_n\).
\end{itemize}

\paragraph{Stationary Distribution of the Count Evolution}

Similarly, we define the stationary count \( C = C_{(\theta, \kappa)} \in \text{Im}(\kappa) \), which models the frequency of occurrence of counts \( c \in \text{Im}(\kappa) \) in the process \(\mathcal{C}_{(\theta, \kappa)} = \{c_n= \kappa(\mathbf{d}_n)\}_{n=0}^\infty\) over an infinite horizon:

\[
\mathbb{P}_C(c) = \mathbb{P}(C_{(\theta, \kappa)} = c) = \lim_{N \to \infty} \frac{\sum_{n=0}^N 1_{\{c_n^{obs} = c\}}}{N+1}
\]

With:
\begin{itemize}
    \item \( \mathcal{C}_{(\theta, \kappa)}^{z=\infty} = \{c_n^{obs}= \kappa(\mathbf{d}_n^{obs})\}_{n=0}^{z=\infty} \): A specific realization of the count evolution up to round \(z=\infty\), where each \(c^{obs}_n\) corresponds to a particular realization of the random count \(c_n\).
\end{itemize}

\paragraph{Equation for Wealth}

By virtue of Theorems 1 and 2, we know that under a CRRA utility function and for infinite rounds, the state variable is solely the initial deck of each round. Assuming that the player has partial information about the deck through a counting system \(\kappa\), we derive an expression for wealth in the long run:

\[
\lim_{N \to \infty} W_{N+1} = \lim_{N \to \infty} \prod_{n=0}^N (1 + \omega(c_n) \cdot X_\theta^{\mathbf{d}_n}),
\]

With:
\begin{itemize}
    \item \( \mathbf{d}_n \in \mathcal{D}_\theta = \{\mathbf{d}_n\}_{n=0}^\infty \): Deck of origin for the \(n\)-th round.
    \item \( c_n \in \mathcal{C}_{(\theta, \kappa)} = \{c_n= \kappa(\mathbf{d}_n)\}_{n=0}^\infty \): Count value for the \(n\)-th round under a counting system \(\kappa\).
\end{itemize}

\paragraph{Wealth Based on the Stationary Deck}

In the limit \(N \to \infty\), each possible factor \((1 + \omega(\kappa(\mathbf{d})) \cdot X_\theta^{\mathbf{d}})\) that contributed to wealth evolution has occurred with a frequency given by \(\mathbb{P}(\mathbf{D_\theta} = \mathbf{d})\). Consequently, we can reformulate the wealth calculation in terms of the stationary deck:

\[
\lim_{N \to \infty} W_{N+1} = \lim_{N \to \infty} \prod_{n=0}^N (1 + \omega(c_n) \cdot X_\theta^{\mathbf{d}_n}) = \lim_{N \to \infty} \prod_{i=0}^N (1 + \omega(C_i) \cdot X_\theta^{\mathbf{D}_i}),
\]

\begin{itemize}
    \item \( \mathbf{D}_i \) represents a collection of i.i.d. variables with distributions \( \mathbf{D}_\theta \).
    \item \( C_i = \kappa(\mathbf{D}_i) \) represents a collection of i.i.d. variables with distributions \( \kappa(\mathbf{D}_\theta) \).
\end{itemize}

\paragraph{Case \(\alpha = 0\): Utility Distribution}

Taking the logarithm of wealth transforms the iterated product into a sum of infinitely many independent and identically distributed random variables. Thus, by the Central Limit Theorem, the sum converges to a normal distribution, yielding the following result:

\[
\lim_{N \to \infty} u_1(W_{N+1}; 0) = \lim_{N \to \infty} \ln(W_{N+1}) \sim \text{Normal}(N \cdot \mu_\omega, N \cdot \sigma_\omega^2),
\]

With:
\begin{itemize}
    \item \( L_\omega = \ln(1 + \omega(\kappa(\mathbf{D})) \cdot X_\theta^{\mathbf{D}}). \)
    \item \( \mu_\omega = \mathbb{E}[L_\omega]. \)
    \item \( \sigma_\omega^2 = \text{Var}[L_\omega]. \)
\end{itemize}

\paragraph{Case \(\alpha > 0\): Utility Distribution}

Since the logarithm of wealth follows a normal distribution, it follows that wealth is log-normally distributed. Using the properties of the log-normal distribution, we deduce the utility distribution of wealth when \(\alpha > 0\), concluding the theorem's proof:

\[
\lim_{N \to \infty} u_1(W_{N+1}; \alpha) = \lim_{N \to \infty} \frac{W_{N+1}^\alpha}{\alpha} \sim \text{LogNormal}(N \cdot \alpha \cdot \mu_\omega - \ln(\alpha), N \cdot \alpha^2 \cdot \sigma_\omega^2),
\]

With:
\begin{itemize}
    \item \( L_\omega = \ln(1 + \omega(\kappa(\mathbf{D})) \cdot X_\theta^{\mathbf{D}}). \)
    \item \( \mu_\omega = \mathbb{E}[L_\omega]. \)
    \item \( \sigma_\omega^2 = \text{Var}[L_\omega]. \)
\end{itemize}
\subsection{Theorem 4}
\label{teo4}

\textit{The optimal partial policy must solve the following optimization problem for each count value:}

\[
\omega_{(\kappa, \theta, u_1, H \to \infty)}^*(c) = \arg \max_{b \in [0, 0.5]} \left\{ \mathbb{E}\left[u_1(1 + b \cdot X_\theta^{(\mathbf{D}|c)}; \alpha)\right] \right\}, \quad \forall c \in \text{Im}(\kappa)
\]

\subsubsection{Proof of Theorem 4}

The problem to solve is formulated as follows:

\[
\omega_{(\kappa, \theta, u_1, H \to \infty)}^* = \arg \max_\omega \{\mathbb{E}[u_1(W_{H \to \infty}; \alpha)]\}
\]

\textbf{Case \( \alpha = 0 \):}

From Theorem 3, the optimization problem for the policy can be explicitly expressed in terms of the long-term expected utility based on the convergent distribution:

\[
\arg \max_\omega \{\mathbb{E}[u_1(W_{H \to \infty}; 0)]\} = \arg \max_\omega \left\{ \mathbb{E}\left[\lim_{H \to \infty} \ln(W_H)\right] \right\} = \arg \max_{\omega} \left\{ \mu_\omega \right\}
\]

We now develop the problem using the algebraic form of \( \mu_\omega \).

\noindent{1. Expanding \( \mu_\omega \):}

\[
\mu_\omega = \mathbb{E}[L_\omega] = \mathbb{E}\left[\ln(1 + \omega(C) \cdot X_\theta^{\mathbf{D}})\right]
\]

\[
\mu_\omega = \sum_{c \in \text{Im}(\kappa)} \mathbb{P}_C(c) \cdot \mathbb{E}\left[\ln(1 + \omega(c) \cdot X_\theta^{(\mathbf{D}|c)})\right]
\]

\noindent{2. From this expression, it follows that maximizing \( \mu_\omega \) is an independent problem for each count value \( c \in \text{Im}(\kappa) \), concluding this part of the proof:}

\[
\omega_{(\kappa, \theta, u_1, H \to \infty)}^*(c) = \arg \max_{b \in [0, 0.5]} \left\{ \mathbb{E}\left[\ln(1 + b \cdot X_\theta^{(\mathbf{D}|c)})\right] \right\}, \quad \forall c \in \text{Im}(\kappa)
\]

\noindent{3. Expanding the expression to maximize:}

\[
\mathbb{E}\left[\ln(1 + b \cdot X_\theta^{(\mathbf{D}|c)})\right] = \sum_{\mathbf{d} \in \Omega_{(\tau, l)} | \kappa(\mathbf{d}) = c} \frac{\mathbb{P}_{\mathbf{D}}(\mathbf{d})}{\mathbb{P}_C(c)} \cdot \left(\sum_{x \in R} \mathbb{P}(X_\theta^{\mathbf{d}} = x) \cdot \ln(1 + b \cdot x)\right)
\]

\noindent{4. Using summation properties to interchange sums, the following is deduced:}

\[
\mathbb{E}\left[\ln(1 + b \cdot X_\theta^{(\mathbf{D}|c)})\right] = \sum_{x \in R} \ln(1 + b \cdot x) \cdot \left(\sum_{\mathbf{d} \in \Omega_{(\tau, l)} | \kappa(\mathbf{d}) = c} \frac{\mathbb{P}_{\mathbf{D}}(\mathbf{d})}{\mathbb{P}_C(c)} \cdot \mathbb{P}(X_\theta^{\mathbf{d}} = x)\right)
\]

\[
\mathbb{E}\left[\ln(1 + b \cdot X_\theta^{(\mathbf{D}|c)})\right] = \sum_{x \in R} \ln(1 + b \cdot x) \cdot \mathbb{E}[\mathbb{P}(X_\theta^{(\mathbf{D}|c)} = x)]
\]

\noindent{5. Finally, the explicit optimality condition that the policy \( \omega^* \) must satisfy for each count value \( c \in \text{Im}(\kappa) \) independently is concluded:}

\[
\omega_{(\kappa, \theta, u_1, H \to \infty)}^*(c) = \arg \max_{b \in [0, 0.5]} \left\{ \sum_{x \in R} \ln(1 + b \cdot x) \cdot \mathbb{E}[\mathbb{P}(X_\theta^{(\mathbf{D}|c)} = x)] \right\}, \quad \forall c \in \text{Im}(\kappa)
\]

\textbf{Case \( \alpha > 0 \):}

From Theorem 3, the optimization problem for the policy can be formulated in terms of wealth expressed with the stationary deck:

\[
\arg \max_\omega \{\mathbb{E}[u_1(W_{H \to \infty}; \alpha)]\} = \arg \max_\omega \left\{ \mathbb{E}\left[\lim_{H \to \infty} \frac{W_H^\alpha}{\alpha}\right] \right\} = \arg \max_{\omega} \left\{ \mathbb{E}\left[\lim_{H \to \infty} \prod_{i=0}^{H-1} (1 + \omega(C_i) \cdot X_\theta^{\mathbf{D}_i})^\alpha\right] \right\}
\]

Where \( \mathbf{D}_i \) are i.i.d. random variables with distributions \( \mathbf{D} \), and \( C_i = \kappa(\mathbf{D}_i) \). The problem develops as follows:

\noindent{1. By the expectation multiplication theorem, the following equality is deduced:}

\[
\mathbb{E}\left[\lim_{H \to \infty} \prod_{i=0}^{H-1} (1 + \omega(C_i) \cdot X_\theta^{\mathbf{D}_i})^\alpha\right] = \lim_{H \to \infty} \mathbb{E}^{H}[(1 + \omega(C) \cdot X_\theta^{\mathbf{D}})^\alpha]
\]

\noindent{2. Using this equality, the optimization problem can be stated as follows:}

\[
\arg \max_\omega \{\mathbb{E}[u_1(W_{H \to \infty}; \alpha)]\}= \arg \max_{\omega} \left\{ \mathbb{E}[(1 + \omega(C) \cdot X_\theta^{\mathbf{D}})^\alpha] \right\}
\]

\noindent{3. Expanding the expression:}

\[
\mathbb{E}[(1 + \omega(C) \cdot X_\theta^{\mathbf{D}})^\alpha] = \sum_{c \in \text{Im}(\kappa)} \mathbb{P}_C(c) \cdot \mathbb{E}[(1 + \omega(c) \cdot X_\theta^{(\mathbf{D}|c)})^\alpha]
\]

\noindent{4. Similar to the case \( \alpha = 0 \), the optimization problem for the policy is deduced to be independent for each count value \( c \in \text{Im}(\kappa) \), concluding the proof of the theorem:}

\[
\omega_{(\kappa, \theta, u_1, H \to \infty)}^*(c) = \arg \max_{b \in [0, 0.5]} \left\{ \mathbb{E}[(1 + b \cdot X_\theta^{(\mathbf{D}|c)})^\alpha] \right\}, \quad \forall c \in \text{Im}(\kappa)
\]

\noindent{5. Analogously to the case \( \alpha = 0 \), summation properties can be used to deduce the following formulation for the optimization problem:}

\[
\omega_{(\kappa, \theta, u_1, H \to \infty)}^*(c) = \arg \max_{b \in [0, 0.5]} \left\{ \sum_{x \in R} (1 + b \cdot x)^\alpha \cdot \mathbb{E}[\mathbb{P}(X_\theta^{(\mathbf{D}|c)} = x)] \right\}, \quad \forall c \in \text{Im}(\kappa)
\]

\end{appendices}
\bibliography{sn-bibliography}
\end{document}